\documentclass[preprint,3p,12pt]{elsarticle}
\usepackage{amsmath}
\usepackage{stmaryrd}
\usepackage{bbding}
\usepackage{dcolumn}
\usepackage{graphicx}
\usepackage{amsfonts}
\usepackage{amssymb}
\usepackage{psfrag}
\usepackage{wrapfig}
\usepackage{subfigure}
\usepackage{makeidx}
\usepackage{bm}
\usepackage{epsf}
\usepackage{epsfig}
\usepackage{setspace}
\usepackage{color}
\usepackage{setspace}
\usepackage[ruled,linesnumbered]{algorithm2e}

\definecolor{SpringGreen4}{RGB}{0,139,69}

\begin{document}
\title{Implicit high-order gas-kinetic schemes for compressible flows on three-dimensional unstructured meshes I: steady flows}

\author[BNU,HKUST1]{Yaqing Yang}
\ead{yqyangbnu@163.com}

\author[BNU]{Liang Pan\corref{cor}}
\ead{panliang@bnu.edu.cn}

\author[HKUST1,HKUST2]{Kun Xu}
\ead{makxu@ust.hk}

\cortext[cor]{Corresponding author}
\address[BNU]{Laboratory of Mathematics and Complex Systems, School of Mathematical Sciences, Beijing Normal University, Beijing, China}
\address[HKUST1]{Department of Mathematics, Hong Kong University of Science and Technology, Clear Water Bay, Kowloon, Hong Kong}
\address[HKUST2]{Shenzhen Research Institute, Hong Kong University of Science and Technology, Shenzhen, China}

\begin{abstract}
In the previous studies, the high-order gas-kinetic schemes (HGKS)
have achieved successes for unsteady flows on three-dimensional unstructured meshes. 
In this paper, to accelerate the rate of convergence  for steady flows, the implicit non-compact and compact
HGKSs are developed.  For non-compact scheme, the simple weighted essentially non-oscillatory (WENO) reconstruction is 
used to achieve the spatial accuracy, where the stencils for reconstruction contain two levels of neighboring cells. 
Incorporate with the nonlinear generalized minimal residual (GMRES) method,  
the implicit non-compact HGKS is developed. In order to improve the 
resolution and parallelism of non-compact HGKS, the implicit  compact HGKS is developed with
Hermite WENO (HWENO) reconstruction, in which the reconstruction stencils only contain one level of neighboring cells. 
The cell averaged conservative variable is also updated with GMRES method. Simultaneously, a simple strategy is used to update the cell averaged 
gradient  by the time evolution of spatial-temporal coupled  gas distribution function.
To accelerate the computation, the implicit non-compact and compact HGKSs are 
implemented with the graphics processing unit (GPU) using compute unified device architecture (CUDA). 
A variety of numerical examples, from the subsonic to supersonic flows, are presented to validate the 
accuracy, robustness and efficiency of both inviscid and viscous flows.
\end{abstract}


\maketitle

\section{Introduction}
The simulation of compressible flows with complex geometry is an important issue for computational fluid dynamics, and 
the unstructured meshes are widely used due to the flexibility. 
For the spatial reconstruction, various high-order numerical methods on unstructured meshes have been developed in the past decades, 
such as essential non-oscillatory (ENO) \cite{ENO-un} and weighted essential non-oscillatory (WENO) \cite{WENO-un1,WENO-un2},  
Hermite WENO (HWENO) methods \cite{HWENO-1,HWENO-2,HWENO-3}, discontinuous Galerkin (DG) \cite{DG, DDG},  
flux reconstruction (FR) \cite{FR} and correction procedure using reconstruction (CPR) \cite{CPR}, etc.  
For the temporal discretization, early efforts mainly focused on explicit schemes on unstructured meshes, 
and the most widely used method is the Runge-Kutta schemes \cite{TVD-RK}.  However, for steady flows, the rate of convergence slows 
down dramatically. In order to speed up the convergence, the implicit temporal discretization is required. 
In general, the implicit method requires to solve a large system of equation, which arises from the linearization 
of a fully implicit scheme at each time step. Several methods are used to solve the large sparse system 
on unstructured meshes, including approximate factorization methods and iterative solution methods.
As an approximate factorization method, the lower-upper symmetric Gauss-Seidel (LUSGS) 
method on structured meshes was originally developed by Jameson and Yoon \cite{LUSGS-2}, and 
it has been successfully extended to unstructured meshes \cite{LUSGS-3,LUSGS-4}.
The most attractive feature of this method is that it does not require any extra memory 
compared with the explicit methods and is free from any matrix inversion. 
However, LUSGS method is not ideally effective because of slow convergence, 
and requires thousands of time steps to achieve the steady state. 
As an iterative method, the most successful and effective method is the Krylov 
subspace method, such as the nonlinear generalized minimal residual (GMRES) 
method \cite{GMRES-1,GMRES-2}.  The GMRES method always has a faster convergence speed, 
but the drawback is that they require a considerable amount of memory to store the Jacobian matrix, 
which may be prohibitive for large problems. To save the storage, the matrix-free GMRES method has been applied to 
unstructured meshes for steady and unsteady flows \cite{GMRES-3,GMRES-4}. To improve the convergence,
a wide variety of preconditioners \cite{Implicit-2,Implicit-3,Implicit-4,Implicit-5,Implicit-6} are 
also applied to the methods above.

In the past decades, the gas-kinetic scheme (GKS) based on the
Bhatnagar-Gross-Krook (BGK) model \cite{BGK-1,BGK-2} have been
developed systematically for the computations from low speed flows
to supersonic ones \cite{GKS-Xu1,GKS-Xu2}. The gas-kinetic scheme
presents a gas evolution process from the kinetic scale to
hydrodynamic scale, and both inviscid and viscous fluxes can be
calculated in one framework.  With the two-stage temporal discretization, 
which was originally developed for the Lax-Wendroff type flow solvers \cite{GRP-high-1}, 
a reliable framework was provided to construct gas-kinetic scheme with
 fourth-order and even higher-order temporal accuracy \cite{GKS-high-1,GKS-high-2}. 
 With the simple WENO type reconstruction, the third-order gas-kinetic schemes
 on three-dimensional unstructured meshes \cite{GKS-WENO-2,GKS-WENO-3} are developed, 
 in which a simple strategy of selecting stencils for reconstruction is adopted and the 
 topology independent linear weights are used. 
 Based on the spatial and temporal coupled property of GKS solver and HWENO reconstruction, 
 the explicit high-order compact gas-kinetic schemes are also developed \cite{GKS-HWENO-1,GKS-HWENO-2,GKS-HWENO-3}. 
In the compact scheme, the time-dependent gas distribution function at a cell interface is used to calculate the fluxes for the updating the 
cell-averaged flow variables, and to evaluate the cell-averaged gradients of flow variables.
Numerical results demonstrate that the superior robustness in high speed flow computation and the
favorable mesh adaptability for complex geometry of high-order compact schemes.
 For the gas-kinetic scheme, several implicit algorithms have also been developed to simulate from continuum 
and rarefied flows \cite{GKS-implicit-1,GKS-WENO-2,GKS-WENO-3,GKS-implicit-2}. The implicit methods 
provide efficient techniques for speeding up the convergence of steady flows.

In this paper, the implicit non-compact and compact HGKSs are developed on the
three-dimensional unstructured meshes. For non-compact GKS scheme, the third-order WENO reconstruction is used, 
where the stencils are selected from the neighboring cells and the neighboring cells of neighboring cells.  
Incorporate with the GMRES method, the implicit non-compact scheme is developed for steady problems. 
In order to balance the computational efficiency and memory storage, 
the GMRES method is based on numerical Jacobian matrix with Roe's approximation. 
To improve the resolution and parallelism, the implicit compact HGKS is also developed with HWENO reconstruction, 
where the stencils only contain one level of neighboring cells.  
Since the cell averaged conservative variables and gradients need to be updated simultaneously, 
the GMRES method is associated with a suitable update strategy at each time step. 
For steady problems, the update of cell averaged gradient can be driven by time evolution of gas distribution function directly.  
To further accelerate the computation, the Jacobi iteration is chosen as preconditioner for both non-compact and compact schemes. 
Various three-dimensional numerical experiments, from the subsonic to supersonic flows, are presented to to validate the accuracy and robustness of current implicit scheme. 
To accelerate the computation, the current schemes are implemented to run on graphics processing unit (GPU) using compute unified device architecture (CUDA). 
The GPU code is implemented with single Nvidia Quadro RTX 8000 GPU, and the CPU code is run with Intel Xeon Gold 6230R CPU with 16 OpenMP threads.
Compared with the CPU code, 8x speedup is achieved for GPU code.  
In the future, more challenging compressible flow problems will be investigated with multiple GPUs.

This paper is organized as follows. In Section 2, BGK equation and finite volume scheme will be introduced. 
The third-order non-compact and compact gas-kinetic scheme will be presented in Section 3. 
Section 4 includes the implicit method and its parallel computation. 
Numerical examples are included in Section 5. The last section is the conclusion.

\section{BGK equation and finite volume scheme}
The Boltzmann equation expresses the behavior of a many-particle
kinetic system in terms of the evolution equation for a single
particle gas distribution function. The BGK equation
\cite{BGK-1,BGK-2}  is the simplification of Boltzmann equation, and
the three-dimensional BGK equation can be written as
\begin{equation}\label{bgk}
f_t+uf_x+vf_y+wf_z=\frac{g-f}{\tau},
\end{equation}
where $\boldsymbol{u}=(u,v,w)$ is the particle velocity, $\tau$ is
the collision time, $f$ is the gas distribution function. $g$ is the
equilibrium state given by Maxwellian distribution
\begin{equation*}
g=\rho(\frac{\lambda}{\pi})^{(N+3)/2}e^{-\lambda[(u-U)^2+(v-V)^2+(w-W)^2+\xi^2]},
\end{equation*}
where $\rho$ is the density, $\boldsymbol{U}=(U,V,W)$ is the
macroscopic fluid velocity, and $\lambda$ is the inverse of gas
temperature, i.e., $\lambda=m/2kT$. In the BGK model, the collision
operator involves a simple relaxation from $f$ to the local equilibrium
state $g$. The variable $\xi$ accounts for the internal degree of
freedom, $\xi^2=\xi_1^2+\dots+\xi_N^2$, $N=(5-3\gamma)/(\gamma-1)$ is
the internal degree of freedom, and $\gamma$ is the specific heat
ratio. The collision term satisfies the compatibility condition
\begin{equation*}
\int \frac{g-f}{\tau}\psi \text{d}\Xi=0,
\end{equation*}
where $\displaystyle\psi=(1,u,v,w,\frac{1}{2}(u^2+v^2+w^2+\xi^2))^T$ and $\text{d}\Xi=\text{d}u\text{d}v\text{d}w\text{d}\xi_1\dots\text{d}\xi_{N}$.
According to the Chapman-Enskog expansion for BGK equation, the
macroscopic governing equations can be derived. In the continuum
region, the BGK equation can be rearranged and the gas distribution
function can be expanded as
\begin{align*}
f=g-\tau D_{\boldsymbol{u}}g+\tau D_{\boldsymbol{u}}(\tau D_{\boldsymbol{u}})g-\tau D_{\boldsymbol{u}}[\tau D_{\boldsymbol{u}}(\tau D_{\boldsymbol{u}})g]+...,
\end{align*}
where $D_{\boldsymbol{u}}=\displaystyle\frac{\partial}{\partial
t}+\boldsymbol{u}\cdot \nabla$. With the zeroth-order truncation
$f=g$, the Euler equations can be obtained. For the first-order
truncation
\begin{align*}
f=g-\tau (ug_x+vg_y+wg_z+g_t),
\end{align*}
the Navier-Stokes equations can be obtained \cite{GKS-Xu1,GKS-Xu2}.

Taking moments of Eq.\eqref{bgk} and integrating with respect to
space, the semi-discretized finite volume scheme can be expressed as
\begin{align}\label{semi}
|\Omega_{i}|\frac{\text{d} Q_i}{\text{d} t}=\mathcal{L}(Q_{i}),
\end{align}
where $Q_i=(\rho, \rho U,\rho V, \rho W, \rho E)$ is the cell averaged
conservative value of $\Omega_{i}$,  $\rho$ is the density, $U,V,W$
is the flow velocity, $\rho E$ is the total energy density and $|\Omega_{i}|$ is the volume of $\Omega_{i}$. The
operator $\mathcal{L}$ is defined as
\begin{equation}\label{finite}
\mathcal{L}(Q_{i})=-\sum_{i_p\in N(i)}F_{i,i_p}(t)S_{i_p}=-\sum_{i_p\in N(i)}\iint_{\Sigma_{i_p}}\boldsymbol{F}(Q,t)\text{d}\sigma,
\end{equation}
where $\Sigma_{i_p}$ is the common cell interface of $\Omega_{i}$, $S_{i_p}$ is the area of $\Sigma_{i_p}$ and
$N(i)$ is the set of index for neighboring cells of
$\Omega_{i}$. To achieve the expected order of accuracy, the
Gaussian quadrature is used for the flux integration
\begin{align*}
\iint_{\Sigma_{i_p}}\boldsymbol{F}(Q,t)\text{d}\sigma=\sum_{G}\omega_{G}F_{G}(t)S_{i_p},
\end{align*}
where $\omega_{G}$ is the quadrature weights. The numerical flux $F_{G}(t)$ at Gaussian quadrature
point $\boldsymbol{x}_{G}$ can be given by taking moments of gas
distribution function
\begin{align*}
F_{G}(t)=\int\boldsymbol\psi \boldsymbol{u}\cdot\boldsymbol{n}_{G} f(\boldsymbol{x}_{G},t,\boldsymbol{u},\xi)\text{d}\Xi,
\end{align*}
where $F_{G}(t)=(F^{\rho}_{G},F^{\rho U}_{G},F^{\rho V}_{G},F^{\rho W}_{G}, F^{\rho E}_{G})$ and
$\boldsymbol{n}_{G}$ is the local normal direction of cell
interface. With the coordinate transformation, the numerical flux in
the global coordinate can be obtained. Based on the integral
solution of BGK equation Eq.\eqref{bgk}, the gas distribution
function $f(\boldsymbol{x}_{G},t,\boldsymbol{u},\xi)$ in the local
coordinate can be given by
\begin{equation*}
f(\boldsymbol{x}_{G},t,\boldsymbol{u},\xi)=\frac{1}{\tau}\int_0^t
g(\boldsymbol{x}',t',\boldsymbol{u}, \xi)e^{-(t-t')/\tau}\text{d}t'+e^{-t/\tau}f_0(-\boldsymbol{u}t,\xi),
\end{equation*}
where $\boldsymbol{x}'=\boldsymbol{x}_{G}-\boldsymbol{u}(t-t')$ is
the trajectory of particles, $f_0$ is the initial gas distribution
function, and $g$ is the corresponding equilibrium state. With the
first order spatial derivatives, the second-order gas distribution
function at cell interface can be expressed as
\begin{align}\label{flux}
f(\boldsymbol{x}_{G},t,\boldsymbol{u},\xi)=&(1-e^{-t/\tau})g_0+
((t+\tau)e^{-t/\tau}-\tau)(\overline{a}_1u+\overline{a}_2v+\overline{a}_3w)g_0\nonumber\\
+&(t-\tau+\tau e^{-t/\tau}){\bar{A}} g_0\nonumber\\
+&e^{-t/\tau}g_r[1-(\tau+t)(a_{1}^{r}u+a_{2}^{r}v+a_{3}^{r}w)-\tau A^r)](1-H(u))\nonumber\\
+&e^{-t/\tau}g_l[1-(\tau+t)(a_{1}^{l}u+a_{2}^{l}v+a_{3}^{l}w)-\tau A^l)]H(u),
\end{align}
where the equilibrium state $g_{0}$ and the corresponding
conservative variables $Q_{0}$ can be determined by the
compatibility condition
\begin{align*}
\int\psi g_{0}\text{d}\Xi=Q_0=\int_{u>0}\psi
g_{l}\text{d}\Xi+\int_{u<0}\psi g_{r}\text{d}\Xi.
\end{align*}
With the reconstruction of macroscopic variables, the coefficients
in Eq.\eqref{flux} can be fully determined by the reconstructed
derivatives and compatibility condition
\begin{equation*}
\begin{aligned} 
\displaystyle 
\langle a_{1}^{k}\rangle=\frac{\partial Q_{k}}{\partial \boldsymbol{n_x}}, 
\langle a_{2}^{k}\rangle=\frac{\partial Q_{k}}{\partial \boldsymbol{n_y}}, 
\langle a_{3}^{k}\rangle&=\frac{\partial Q_{k}}{\partial\boldsymbol{n_z}}, 
\langle
a_{1}^{k}u+a_{2}^{k}v+a_{3}^{k}w+A^{k}\rangle=0,\\ 
\displaystyle
\langle\overline{a}_1\rangle=\frac{\partial Q_{0}}{\partial \boldsymbol{n_x}}, 
\langle\overline{a}_2\rangle=\frac{\partial Q_{0}}{\partial \boldsymbol{n_y}},
\langle\overline{a}_3\rangle&=\frac{\partial Q_{0}}{\partial \boldsymbol{n_z}},
\langle\overline{a}_1u+\overline{a}_2v+\overline{a}_3w+\overline{A}\rangle=0,
\end{aligned}
\end{equation*}
where $k=l$ and $r$,  $\boldsymbol{n_x}$, $\boldsymbol{n_y}$, $\boldsymbol{n_z}$ are the unit directions of local coordinate at $\boldsymbol{x}_{G}$
 and $\langle...\rangle$ are the moments of the
equilibrium $g$ and defined by
\begin{align*}
\langle...\rangle=\int g (...)\psi \text{d}\Xi.
\end{align*}
More details of the gas-kinetic scheme can be found in \cite{GKS-Xu1,GKS-Xu2}. 

\section{Spatial reconstruction}
To deal with the complex geometry, the three-dimensional
unstructured meshes are considered, and the tetrahedral and
hexahedral meshes are used in this paper for simplicity. In the
previous studies, the high-order gas-kinetic schemes have be
developed with the third-order non-compact WENO reconstruction
\cite{GKS-WENO-2, GKS-WENO-3} and compact HWENO reconstruction
\cite{GKS-HWENO-2, GKS-HWENO-3}. Successes have been achieved for
the unsteady flows from the subsonic to supersonic flow problems. For the
unstructured meshes,  the classical WENO reconstruction, in which
the high-order accuracy is achieved by non-linear combination of
lower order polynomials, becomes extremely complicated for
three-dimensional problems \cite{WENO-un1,WENO-un2}. In this paper,
the idea of simple WENO reconstruction is adopted
\cite{WENO-simple-1,WENO-simple-2}. The high-order accuracy is
achieved by non-linear combination of high-order and lower-order
polynomials, and the constant linear weights are adopted.

\subsection{Selection of stencils}
For the cell $\Omega_i$, the faces are labeled as $F_{p}$, where
$p=1,\dots,4$ for tetrahedral cell, and $p=1,\dots,6$ for hexahedral
cell. The neighboring cell of $\Omega_{i}$, which shares the face
$F_{p}$, is denoted as $\Omega_{i_p}$. Meanwhile, the neighboring
cells of $\Omega_{i_p}$ are denoted as $\Omega_{i_{pm}}$. To achieve
the third-order accuracy, a big stencil for non-compact
reconstruction for cell $\Omega_i$  is selected as follows
\begin{align*}
S_i^{WENO}&=\{\Omega_{i},\Omega_{i_p},\Omega_{i_{pm}}\},
\end{align*}
which is consist of the neighboring cells and neighboring cells of
neighboring cells of $\Omega_{i}$.  Meanwhile, a big stencil for
compact reconstruction for cell $\Omega_i$  is selected as follows
\begin{align*}
S_i^{HWENO}&=\{\Omega_{i},\Omega_{i_p}\},
\end{align*}
which is only consist of the neighboring cells of $\Omega_{i}$.

To deal with the discontinuity, the sub-stencils $S_{i_m}^{WENO}$ in
non-compact WENO reconstruction for cell $\Omega_{i}$ are selected,
where $m=1,\dots, M$ and $M$ is the number of sub-stencils. For the
hexahedral cell, $M=8$ and the sub-candidate stencils are selected
as
\begin{align*}
S_{i_1}^{WENO}=\{\Omega_{i},\Omega_{i_1},\Omega_{i_2},\Omega_{i_3}\},~S_{i_5}^{WENO}=\{\Omega_{i},\Omega_{i_6},\Omega_{i_2},\Omega_{i_3}\},\\
S_{i_2}^{WENO}=\{\Omega_{i},\Omega_{i_1},\Omega_{i_3},\Omega_{i_4}\},~S_{i_6}^{WENO}=\{\Omega_{i},\Omega_{i_6},\Omega_{i_3},\Omega_{i_4}\},\\
S_{i_3}^{WENO}=\{\Omega_{i},\Omega_{i_1},\Omega_{i_4},\Omega_{i_5}\},~S_{i_7}^{WENO}=\{\Omega_{i},\Omega_{i_6},\Omega_{i_4},\Omega_{i_5}\},\\
S_{i_4}^{WENO}=\{\Omega_{i},\Omega_{i_1},\Omega_{i_5},\Omega_{i_2}\},~S_{i_8}^{WENO}=\{\Omega_{i},\Omega_{i_6},\Omega_{i_5},\Omega_{i_2}\}.
\end{align*}
The linear polynomials can be determined based on above stencils,
which contain the target cell $\Omega_{i}$ and three neighboring
cells. For the tetrahedral cells, in order to avoid the centroids of
$\Omega_{i}$ and three of neighboring cells becoming coplanar,
additional cells are needed for the sub-candidate stencils. For the
tetrahedral cell,  four sub-candidate stencils are selected as
\begin{align*}
S_{i_1}^{WENO}=\{\Omega_{i},\Omega_{i_1},\Omega_{i_2},\Omega_{i_3},\Omega_{i_{11}},\Omega_{i_{12}},\Omega_{i_{13}}\},\\
S_{i_2}^{WENO}=\{\Omega_{i},\Omega_{i_1},\Omega_{i_2},\Omega_{i_4},\Omega_{i_{21}},\Omega_{i_{22}},\Omega_{i_{23}}\},\\
S_{i_3}^{WENO}=\{\Omega_{i},\Omega_{i_2},\Omega_{i_3},\Omega_{i_4},\Omega_{i_{31}},\Omega_{i_{32}},\Omega_{i_{33}}\},\\
S_{i_4}^{WENO}=\{\Omega_{i},\Omega_{i_3},\Omega_{i_1},\Omega_{i_4},\Omega_{i_{41}},\Omega_{i_{42}},\Omega_{i_{43}}\}.
\end{align*}
where $\Omega_{i_{pn}} \neq \Omega_{i}, p=1,\dots,4$ and $n=1,2,3$.
The cells of sub-candidate stencils are consist of the three
neighboring cells and three neighboring cells of one neighboring
cell. With such an enlarged sub-stencils, the linear polynomials can
be determined.

Meanwhile, the sub-stencils $S_{i_m}^{HWENO}$ in compact HWENO
reconstruction for cell $\Omega_{i}$ can be selected more simply.
The sub-candidate stencils are selected as
\begin{align*}
S_{i_m}^{HWENO}=\{\Omega_{i},\Omega_{i_f}\},
\end{align*}
where $\Omega_{i_f}$ is one of the neighboring cell of the target
cell $\Omega_{i}$, $m=1,\dots,M$ and $M$ equal to the number of the
faces of cell $\Omega_{i}$. The linear polynomials can be determined
on such small stencils with additional degree of freedoms. Noticed
that the sub-stencils $S_{i_m}^{WENO}$ for hexahedral cells also
belong to compact sub-candidate stencils.

\subsection{Non-compact WENO reconstruction}
For the target cell $\Omega_{i}$, the big stencil $S_i^{WENO}$ is rearranged as
$\{\Omega_0,\Omega_1,...,\Omega_{K}\}$ and the sub-stencil $S_{i_m}^{WENO}$ is rearranged as $\{\Omega_0,\Omega_{m_1},…,\Omega_{m_K}\}$, where $\Omega_0$ is the target cell. A quadratic polynomial and several linear polynomials can be constructed based on the big stencil $S_i^{WENO}$ and the sub-stencils $S_{i_m}^{WENO}$ as follows
\begin{equation}\label{polys}
\begin{split}
P_0(\boldsymbol{x})&=Q_{0}+\sum_{|\boldsymbol d|=1}^2a_{\boldsymbol d}p_{\boldsymbol d}(\boldsymbol{x}),\\
P_m(\boldsymbol{x})&=Q_{0}+\sum_{|\boldsymbol d|=1}b_{\boldsymbol d}^mp_{\boldsymbol d}(\boldsymbol{x}),
\end{split}
\end{equation}
where $m=1,\dots,M$, $Q_{0}$ is the cell averaged variables over $\Omega_{0}$ with newly rearranged index, the multi-index $\boldsymbol d=(d_1, d_2,
d_3)$ and $|\boldsymbol d|=d_1+d_2+d_3$. The base function
$p_{\boldsymbol d}(\boldsymbol{x})$ is defined as
\begin{align*}
\displaystyle
p_{\boldsymbol d}(\boldsymbol{x})=x^{d_1}y^{d_2}z^{d_3}-\frac{1}{\left|\Omega_{0}\right|}\iiint_{\Omega_{0}}x^{d_1}y^{d_2}z^{d_3}\text{d}V.
\end{align*}
To determine these polynomials, the following constrains need to be
satisfied
\begin{align*}
\frac{1}{\left|\Omega_{k}\right|}\iiint_{\Omega_{k}}P_0(\boldsymbol{x})\text{d}V&=Q_{k},~\Omega_{k}\in S_i^{WENO},\\
\frac{1}{\left|\Omega_{m_k}\right|}\iiint_{\Omega_{m_k}}P_m(\boldsymbol{x})\text{d}V&=Q_{m_k},~\Omega_{m_k}\in S_{i_m}^{WENO},
\end{align*}
where $k=0, \dots, K$, $m_k=0,\dots,m_K$, $Q_{k}$ and
$Q_{m_k}$ are the conservative variable with newly rearranged index.
The over-determined linear systems can be generated and the least
square method is used to obtain the coefficients $a_{\boldsymbol d}$
and $b_{\boldsymbol d}^m$.

\subsection{Compact HWENO reconstruction}
For the target cell $\Omega_{i}$, the big stencil $S_i^{HWENO}$ is rearranged as $\{\Omega_0,\Omega_1,...,\Omega_{K}\}$ and the sub-stencil $S_{i_m}^{HWENO}$ is rearranged as $\{\Omega_0,\Omega_{m_1}\}$, where $\Omega_0$ is the target cell. The quadratic polynomial and linear polynomials in Eq.\eqref{polys} can be also reconstructed
based on the compact big stencil $S_i^{HWENO}$ and sub-stencils
$S_{i_m}^{HWENO}$ respectively. To determine these polynomials with
a smaller stencil, the additional constrains need to be added for
all cells as follows
\begin{equation}\label{compact-big-1}
\begin{split}
\frac{1}{\left|\Omega_{k}\right|}\iiint_{\Omega_{k}}P_0(\boldsymbol{x})\text{d}V&=Q_{k},~\Omega_{k}\in S_i^{HWENO},\\
\frac{1}{\left|\Omega_{k}\right|}\iiint_{\Omega_{k}}\frac{\partial}{\partial {\boldsymbol \tau}}P_0(\boldsymbol{x})\text{d}V&=(Q_{\boldsymbol{\tau}})_{k},~\Omega_{k}\in S_i^{HWENO},
\end{split}
\end{equation}
and
\begin{equation}\label{compact-big-2}
\begin{split}
\frac{1}{\left|\Omega_{m_k}\right|}\iiint_{\Omega_{m_k}}P_m(\boldsymbol{x})\text{d}V&=Q_{m_k},~\Omega_{m_k}\in S_{i_m}^{HWENO},\\
\frac{1}{\left|\Omega_{m_k}\right|}\iiint_{\Omega_{m_k}}\frac{\partial}{\partial {\boldsymbol \tau}}P_m(\boldsymbol{z})\text{d}V&=(Q_{\boldsymbol{\tau}})_{m_k},~\Omega_{m_k}\in S_{i_m}^{HWENO},
\end{split}
\end{equation}
where $\boldsymbol{\tau}$ is the unit directions of local coordinate, i.e., $\boldsymbol{n}_x, \boldsymbol{n}_y, \boldsymbol{n}_z$, $k=0, \dots, K$, $m_k=0,m_1$, $Q_{k}$, $Q_{k}$ and $Q_{m_k}$ are the cell averaged conservative variables and $(Q_{\boldsymbol{n}_x})_{k}$, $(Q_{\boldsymbol{n}_y})_{k}$, $(Q_{\boldsymbol{n}_z})_{k}$ and $(Q_{\boldsymbol{n}_x})_{m_k}$,
$(Q_{\boldsymbol{n}_y})_{m_k}$, $(Q_{\boldsymbol{n}_z})_{m_k}$ are
the cell averaged directional derivatives
over $\Omega_{k}\in S_i^{HWENO}$ and $\Omega_{m_k}\in
S_{i_m}^{HWENO}$ respectively with newly rearranged index.

In order to solve the systems in Eq.\eqref{compact-big-1} and Eq.\eqref{compact-big-2}, the cell
averaged directional derivatives on the right-hand side need to be
calculated. Different from the traditional Riemann solvers, the
gas-kinetic scheme provides a time-dependent gas distribution
function by Eq.\eqref{flux}. Meanwhile,  the macroscopic conservative variables can be obtained by taking
moments of the distribution function as well 
\begin{align}\label{point}
Q(\boldsymbol{x}_G,t)=\int\boldsymbol\psi f(\boldsymbol{x}_G,t,\boldsymbol{u},\xi)\text{d}\Xi.
\end{align}
According to the Gauss's theorem, the cell averaged gradient of the
flow variable $Q$ can be calculated as follows
\begin{equation}\label{derivative}
\begin{split}
|\Omega_k|(\nabla Q)_k(t)&=\iiint_{\Omega_k}\nabla Q(t)\text{d}V=\iint_{\partial\Omega_k}Q(t)\boldsymbol{\tau}\text{d}S \\
&=\sum_{i_p\in N(k)} \Big( \sum_{G}\omega_G \Big(\int\boldsymbol\psi f(\boldsymbol{x}_G,t,\boldsymbol{u},\xi)\text{d}\Xi\Big) \boldsymbol{\tau} S_{i_p} \Big),
\end{split}
\end{equation}
where $\Omega_k$ is arbitrary cell in computational mesh, $(\nabla Q)_k$ is the cell averaged gradient of the flow 
variable $Q_k$ over cell $\Omega_k$, $\nabla Q$ is the distribution of flow gradient, 
$\boldsymbol{\tau}=( \boldsymbol{n_x}, \boldsymbol{n_y}, \boldsymbol{n_z})$ are the unit directions of local coordinate on the cell interface.
With Gaussian quadrature rule, the cell-averaged gradient $(\nabla Q)_k$ at $t=t^n$ can be
obtained. It means that the computation
of cell-averaged gradient $(\nabla Q)_k$ can be turned to calculate
the point value of $Q$ of cell interface in a local orthogonal
coordinate. In practice, the derivative parts in
Eq.\eqref{compact-big-1} and  Eq.\eqref{compact-big-2}  are scaled by $h=|\Omega_{k}|^{1/3}$ and
$|\Omega_{m_k}|^{1/3}$ for a smaller condition number of the matrix
in the linear system of $a_{\boldsymbol d}$ and $b_{\boldsymbol d}$.
If the coefficients $a_{\boldsymbol d}$ and $b_{\boldsymbol d}$ are solved by the
least square method, the linear instability will be introduced for
the compact reconstruction. In order to overcome this drawback, the
constrained least-square method is used for solving the above linear
systems, where the conservative variable equations are set as
strictly satisfied and others are satisfied in the sense of least
square \cite{CLS}.

\subsection{Non-linear combination}
With the reconstructed polynomial $P_m(\boldsymbol{x}), m=0,...,M$,
the point-value $Q(\boldsymbol{x}_{G})$ and the spatial derivatives
$\partial_{x,y,z} Q(\boldsymbol{x}_{G})$  for reconstructed
variables  at Gaussian quadrature point can be given by the
non-linear combination
\begin{equation}\label{weno}
\begin{split}
Q(\boldsymbol{x}_{G})=\overline{\omega}_0(\frac{1}{\gamma_0}P_0(\boldsymbol{x}_{G})-&
\sum_{m=1}^{M}\frac{\gamma_m}{\gamma_0}P_m(\boldsymbol{x}_{G}))+\sum_{m=1}^{M}\overline{\omega}_mP_m(\boldsymbol{x}_{G}),\\
\partial_{x,y,z} Q(\boldsymbol{x}_{G})=\overline{\omega}_0(\frac{1}{\gamma_0}\partial_{x,y,z}
P_0(\boldsymbol{x}_{G})-&\sum_{m=1}^{M}\frac{\gamma_m}{\gamma_0}\partial_{x,y,z}
P_m(\boldsymbol{x}_{G}))+\sum_{m=1}^{M}\overline{\omega}_m\partial_{x,y,z}
P_m(\boldsymbol{x}_{G}),
\end{split}
\end{equation}
where $\gamma_0, \gamma_1,\dots,\gamma_M$ are the linear weights. The
non-linear weights $\omega_m$ and normalized non-linear weights
$\overline{\omega}_m$ are defined as
\begin{align*}
\overline{\omega}_{m}=\frac{\omega_{m}}{\sum_{m=0}^{M} \omega_{m}},~
\omega_{m}=\gamma_{m}\Big[1+\Big(\frac{\tau_Z}{\beta_{m}+\epsilon}\Big)\Big],~
\tau_Z=\sum_{m=1}^{M}\Big(\frac{|\beta_0-\beta_m|}{M}\Big),
\end{align*}
where $\epsilon$ is a small positive number.
The smooth indicator $\beta_{m}$ is defined as
\begin{align*}
\beta_m=\sum_{|l|=1}^{r_m}|\Omega_{i}|^{\frac{2|l|}{3}-1}\int_{\Omega_{i}}
\Big(\frac{\partial^lP_m}{\partial_x^{l_1}\partial_y^{l_2}\partial_z^{l_3}}(x,y,z)\Big)^2\text{d}V,
\end{align*}
where $r_0=2$ and $r_m=1$ for $m=1,\dots,M$. It can be proved that
Eq.\eqref{weno} ensures third-order accuracy and more details can be
found in \cite{GKS-WENO-2, GKS-HWENO-2}. In the computation, the
linear weights are set as $\gamma_i=0.025$, $\gamma_0=1-\gamma_i M$
for both non-compact and compact scheme without special statement.

\section{Implicit method}
\subsection{Implicit method for non-compact scheme}
The backward Euler method for the semi discretized scheme
Eq.\eqref{semi} at $t^{n+1}$ is given by
\begin{equation}\label{impEuler}
\frac{|\Omega_i|}{\Delta t}\Delta Q_i^n=\mathcal{L}(Q_{i}^{n+1}),
\end{equation}
where $\Delta Q_i^n=Q_{i}^{n+1}-Q_{i}^{n}$, $\Delta t$ is the time
step and $\mathcal{L}(Q_{i}^{n+1})$ can be linearized as
\begin{equation*}
\mathcal{L}(Q_{i}^{n+1})\approx \mathcal{L}(Q_{i}^{n})+(\frac{\mathrm{d} \mathcal{L}}{\mathrm{d} Q})\Delta Q_i^n.
\end{equation*}
where $\displaystyle\frac{\mathrm{d}\mathcal{L}}{\mathrm{d} Q}$ is the Jacobian matrix. 
Eq.\eqref{impEuler} can be rewritten as
\begin{equation}\label{impsystem}
\Big(\frac{|\Omega_i|}{\Delta t}I-(\frac{\mathrm{d} \mathcal{L}}{\mathrm{d} Q})\Big)\Delta Q_i^n=\mathcal{L}(Q_{i}^{n}),
\end{equation}
where $I$ is an identity matrix. In order to preserve the temporal evolution advantage brought by the
GKS solver, the residual term at time step $t=t^n$ for both
non-compact and compact schemes is given by taking integration of Eq.\eqref{finite} over the time interval
\begin{equation*}
\displaystyle \mathcal{L}(Q_{i}^{n})=-\sum_{i_p\in N(i)}\Big(\sum_{G}\omega_G\Big(\frac{1}{\Delta t} \int_{0}^{\Delta t}\int\psi \boldsymbol{u}f(\boldsymbol{x}_G,t,\boldsymbol{u},\xi)\boldsymbol{n}_G\text{d}\Xi\text{d}t\Big)S_{i_p}\Big).
\end{equation*}

In the previous study, LUSGS method was implemented for the
high-order GKS on the unstructured meshes \cite{GKS-WENO-2, GKS-WENO-3}. For the steady state
problem, LUSGS method converges more efficiently than the explicit methods. However, it still requires a great number of time steps to
achieve a satisfactory final residual.  In order to speed up the
convergence, the Newton-GMRES method \cite{GMRES-1, GMRES-2} will be
implemented for the high-order GKS on the three-dimensional
unstructured meshes. The linearized equation Eq.\eqref{impsystem} can be simply written into the following form
\begin{equation*}
\boldsymbol{A}\Delta Q^n=\boldsymbol{R}^n.
\end{equation*}
The purpose of the GMRES method is finding an approximate solution
over a finite dimensional orthogonal space
\begin{equation*}
\Delta Q^n\approx\Delta Q^{(m)}=\Delta Q^{(0)}+\boldsymbol{u}_m,
\end{equation*}
such that
\begin{equation*}
\| \boldsymbol{r}^{(m)}\|_2=\mathop{\min}\limits_{\boldsymbol{u}_m\in K_m}\|\boldsymbol{R}^n-\boldsymbol{A}(\Delta Q^{(0)}+\boldsymbol{u}_m)\|_2
=\mathop{\min}\limits_{\boldsymbol{u}_m\in K_m}\|\boldsymbol{r}^{(0)}-\boldsymbol{A}\boldsymbol{u}_m\|_2.
\end{equation*}
With the initial condition $Q^n$ and arbitrary initial solution
$\Delta Q^{(0)}$, which is chosen as $\Delta Q^{(0)}=0$, the initial
residual is given by
\begin{equation*}
\boldsymbol{r}^{(0)}=\boldsymbol{R}^n-\boldsymbol{A}\Delta Q^{(0)}.
\end{equation*}
The Krylov subspace of $\boldsymbol{r}^{(0)}$, which is used for
searching the approximation solution,  can be generated by
\begin{equation*}
K_m\equiv\mathrm{span}\{\boldsymbol{v}_{1}, \boldsymbol{v}_{2},\dots,\boldsymbol{v}_{m}\}=
\mathrm{span}\{\boldsymbol{r}^{(0)}, \boldsymbol{A}\boldsymbol{r}^{(0)},\dots,\boldsymbol{A}^{m-1}\boldsymbol{r}^{(0)}\}.
\end{equation*}
To save the memory storage, the matrix-free GMRES method is widely
employed  \cite{GMRES-3, GMRES-4}. In order to avoid calculating the
Jabobian matrix in $\boldsymbol{A}$ directly, the $F$-derivative is
used as follows
\begin{equation}\label{frechet}
(\frac{\partial \mathcal{L}}{\partial Q}) \Delta Q^n \approx\frac{\mathcal{L}(Q^n+ \sigma \Delta Q^n)-\mathcal{L}(Q^n)}{\sigma},
\end{equation}
where the subscript is omitted and $\sigma$ is a small scalar to be chosen carefully. Although
the above matrix-free approximation reduces the memory storage
requirement, the reconstruction processes combined with GKS solver
involved in Eq.\eqref{frechet} will cost too much time, which will
reduce the computational efficiency greatly. Taking the balance of
computational efficiency and memory storage into account, the GMRES
method with an approximation matrix is employed in this paper. The
system Eq.\eqref{impEuler} can be written into the flux form 
\begin{equation}\label{A-system-1}
\frac{|\Omega_i|}{\Delta t}\Delta Q_i^n+\sum_{i_p\in N(i)}\left(F_{i,i_p}^{n+1}-F_{i,i_p}^{n}\right) S_{i_p}=-\sum_{i_p\in N(i)}F_{i,i_p}^{n} S_{i_p}= \mathcal{L}(Q_{i}^{n}).
\end{equation}
According to the total differential formulation,
\begin{align*}
F_{i,i_p}^{n+1}-F_{i,i_p}^{n}&=
F(Q_i^n+\Delta Q_i^n,Q_{i_p}^n+\Delta Q_{i_p}^n)-F(Q_i^n,Q_{i_p}^n)\\
&=(\frac{\partial F}{\partial Q_i})^n\Delta Q_i^n+(\frac{\partial F}{\partial Q_{i_p}})^n\Delta Q_{i_p}^n.
\end{align*}
Substituting the term above into Eq.\eqref{A-system-1}, the backward
Euler method can be rewritten as
\begin{equation}\label{A-system-2}
\Big(\frac{|\Omega_i|}{\Delta t}I+\sum_{i_p\in N(i)}(\frac{\partial F}{\partial Q_i})^n S_{i_p}\Big)\Delta Q_i^n
+\sum_{i_p\in N(i)}\Big((\frac{\partial F}{\partial Q_{i_p}})^n S_{i_p}\Big)\Delta Q_{i_p}^n= \mathcal{L}(Q_{i}^{n}).
\end{equation}
According to Roe's flux function
\begin{equation}\label{deri-matrix}
\begin{cases}
\displaystyle\frac{\partial F}{\partial Q_i}=\frac{1}{2}(J(Q_i)+|\lambda_{i,i_p}|I),\\
\displaystyle\frac{\partial F}{\partial Q_{i_p}}=\frac{1}{2}(J(Q_{i_p})-|\lambda_{i,i_p}|I),
\end{cases}
\end{equation}
where
\begin{equation*}
|\lambda_{i,i_p}|\geq|\boldsymbol{u}_{i,i_p}\cdot \boldsymbol{n}_{i,i_p}|+a_{i,i_p},
\end{equation*}
$\boldsymbol{u}_{i,i_p}$ and $a_{i,i_p}$ are the velocity and the speed of sound as at
 the interface $S_{i_p}$, $\boldsymbol{n}_{i,i_p}$ is the unit normal direction on the cell 
 interface $S_{i_p}$, and $J(Q)$ represents the Jacobian of the inviscid flux. Substituting Eq.\eqref{deri-matrix} into
Eq.\eqref{A-system-2}, the coefficient matrix $\boldsymbol{A}$ of
Eq.\eqref{A-system-2} can be determined. In order to avoid storing
the whole matrix $\boldsymbol{A}$ directly and be able to calculate
the matrix multiplication in parallel, a suitable strategy for the
storage of matrix $\boldsymbol{A}$ and its parallel computation is
designed, and the details will be shown in following section.

\subsection{Implicit method for compact scheme}
In this section, the GMRES method is combined with third-order
compact HWENO reconstruction for solving the steady state problems.
In practice, the coefficient matrix $\boldsymbol{A}$ of
Eq.\eqref{A-system-2} can be fixed with the known cell averaged
conservative variable $Q_i^n$ at each time step, which has nothing
to do with the cell averaged derivatives. Meanwhile, the right side
$\mathcal{L}(Q_{i}^{n})$ is calculated by the compact HWENO
reconstruction with the known cell averaged conservative variable
$Q_i^n$ and cell averaged directional derivatives
$(Q_{\boldsymbol{n}_x})_i^n$, $(Q_{\boldsymbol{n}_y})_i^n$,
$(Q_{\boldsymbol{n}_z})_i^n$. After solving the system
Eq.(\ref{A-system-2}) at each time step, the cell averaged
conservative variable $Q_i^{n+1}$ and cell averaged derivatives
$(Q_{\boldsymbol{n}_x})_i^{n+1}$, $(Q_{\boldsymbol{n}_y})_i^{n+1}$,
$(Q_{\boldsymbol{n}_z})_i^{n+1}$ should be updated simultaneously.
So that $\mathcal{L}(Q_{i}^{n+1})$ can be obtained as the right side
of Eq.\eqref{A-system-2} in the next time step.

For the implicit compact scheme, the GMRES method need to associate
with a suitable strategy for updating cell averaged directional
derivatives at each time step. Due to the small change of flow
variables in the steady state problems, the update of cell averaged
directional derivatives can be driven by time evolution directly. 
Similar technique has been used in other equations \cite{HWENO-FSM}.
According to Eq.\eqref{point} and
Eq.\eqref{derivative}, the cell averaged directional derivatives at
$t=t^{n+1}$ can be updated directly by
\begin{align*}
(Q_{\boldsymbol{\tau}})_i^{n+1}&=(Q_{\boldsymbol{\tau}})_i^{n}=
\frac{1}{|\Omega_i|}\sum_{i_p\in N(i)}\Big(\sum_{G}\omega_G\Big(\int\psi f(\boldsymbol{x}_G,t^n,\boldsymbol{u},\xi)\boldsymbol{\tau}\text{d}\Xi\Big)S_{i_p}\Big).
\end{align*}
where $\boldsymbol{\tau}=\boldsymbol{n}_x,\boldsymbol{n}_y,\boldsymbol{n}_z$.

The numerical tests show that the above update strategy is suitable
for steady flows, and the high-order spatial accuracy can be
maintained in the computation. In the future, the optimization of
above implicit update strategy will be studied for unsteady flows,
for example, designing nonlinear weights in time direction for
preserving temporal and spatial accuracy.

Algorithm.\ref{PGMRES-algorithm} shows the whole
process of  preconditioned GMRES method combined with WENO and HWENO
reconstruction, where the red lines are special steps of the
non-compact scheme, the blue lines are special steps of the compact
scheme and the black lines are common steps for both two schemes.

To ensure the convergence rate of the GMRES method, the condition
number of the corresponding matrix $\boldsymbol{A}$ has to be as
small as possible and the initial solution of flow evolution is
preferably within the convergence domain of Newton iteration. In the
computation, the GMRES method combined with preconditioning
technique is considered for non-compact and compact schemes. It
means that instead of solving the system
\begin{equation*}
\boldsymbol{A} \boldsymbol{x}=\boldsymbol{b},
\end{equation*}
an equivalent preconditioned linear system \cite{GMRES-3} is considered
\begin{equation}\label{precond}
\boldsymbol{P}^{-1}\boldsymbol{A} \boldsymbol{x}=\boldsymbol{P}^{-1}\boldsymbol{b},
\end{equation}
where $\boldsymbol{P}$ is the preconditioning matrix which is an
approximation of matrix $\boldsymbol{A}$.  Usually, the LUSGS method
is adopted as preconditioner of GMRES method, but the serial sweep
part in LUSGS method is not easy to be implemented for parallel
computation. In this paper, Jacobi iteration is adopted as a
preconditioner, which is easy to be implemented in parallel. According to Eq.\eqref{A-system-1}, the matrix
$\boldsymbol{A}$ of Eq.\eqref{A-system-1} can be written as
\begin{equation*}
\boldsymbol{A}=\boldsymbol{D}+\boldsymbol{L}+\boldsymbol{U},
\end{equation*}
where $\boldsymbol{D}$ represents the diagonal area,
$\boldsymbol{L}$ represents the lower triangular area and
$\boldsymbol{U}$ represents the upper triangular area of matrix
$\boldsymbol{A}$.  The Jacobi iteration is used as preconditioner to provide $\boldsymbol{P}^{-1}\boldsymbol{b}$ as follows
\begin{align*}
&\boldsymbol{b}^{0}=\boldsymbol{D}^{-1}\boldsymbol{b}^{\text{ini}},\\
&\boldsymbol{b}^{k}=\boldsymbol{D}^{-1}(\boldsymbol{b}^{\text{ini}}-(\boldsymbol{L}+\boldsymbol{U})\boldsymbol{b}^{k-1}),\\
&\boldsymbol{P}^{-1}\boldsymbol{b}=\boldsymbol{b}^{k_{\text{max}}}.
\end{align*}
where $1\le k\le k_{\max}$ and $\boldsymbol{b}^{\text{ini}}$ represents  $\boldsymbol{r}^{0}$ and $\boldsymbol{A}\boldsymbol{v}_{j}$
 in Line 10 and Line 15 of Algorithm.\ref{PGMRES-algorithm}. With the process above, 
$\boldsymbol{P}^{-1}\boldsymbol{r}^{0}$ and $\boldsymbol{P}^{-1}\boldsymbol{A}\boldsymbol{v}_{j}$ can be fully given.

\begin{algorithm}[!h]
    \SetAlgoLined
    \textcolor{red}{Initial condition $Q^n$ for non-compact scheme}\;
    \textcolor{blue}{Initial condition $Q^n$, $(Q_{\boldsymbol{n}_x})^n$, $(Q_{\boldsymbol{n}_y})^n$, $(Q_{\boldsymbol{n}_z})^n$ for compact scheme}\;
    \While{residual $ \leq $ tolerance}{
        Calculation of time step\;
        \textcolor{red}{WENO reconstruction to calculate $\mathcal{L}(Q^n)$}\;
        \textcolor{blue}{HWENO reconstruction to calculate $\mathcal{L}(Q^n)$ and $f^{n}$}\;
        Calculation of $\boldsymbol{A}$\;
        \For{$i=1, restart\,\,times$}{
        Calculation of $\boldsymbol{r}^{0}$ and residual\;
        Jacobi preconditioning: $\hat{\boldsymbol{r}^0}=\boldsymbol{P}^{-1}\boldsymbol{r}^{0}$\;
        Arnoldi process:\\
        $\boldsymbol{v}_{1}=\hat{\boldsymbol{r}^0}/||\hat{\boldsymbol{r}^0}||_2$\;
        \For{$j=1,  dim K_m$}{
        $\boldsymbol{y}_j=\boldsymbol{A}\boldsymbol{v}_{j}$\;
        Jacobi preconditioning: $\boldsymbol{w}_j=\boldsymbol{P}^{-1}\boldsymbol{y}_j$\;
        \For{$i=1, j$}{
        $h_{ij}=(\boldsymbol{w}_j,\boldsymbol{v}_{i})$\;
        $\boldsymbol{w}_j=\boldsymbol{w}_j-h_{ij}\boldsymbol{v}_{i}$\;
        }
        $h_{j+1,j}=||\boldsymbol{w}_j||_2$\;
        $\boldsymbol{v}_{j+1}=\boldsymbol{w}_j/h_{j+1,j}$\;
        }
        Minimization process: \\
        QR decomposition and solve an upper triangular matrix\;
        Get the solution vector $y_m$\;
        $\Delta Q^n=\Delta Q^n+K_m y_m$\;
        }
        Update:\\
        $Q^{n+1}=Q^{n}+\Delta Q^n$\;
        \textcolor{blue}{Calculate $(Q_{\boldsymbol{n}_x})^{n+1}$, $(Q_{\boldsymbol{n}_y})^{n+1}$, $(Q_{\boldsymbol{n}_z})^{n+1}$ by $f^{n}$ for compact scheme (HWENO)}\;}
    Output of flow field and residual convergence history.
\caption{\label{PGMRES-algorithm} Program for preconditioned GMRES method}
\end{algorithm}

\subsection{Parallel computation for implicit method}
In the previous work, the computation of HGKS is mainly based on the central processing unit (CPU) code. 
To improve the efficiency, the OpenMP directives and message passing interface (MPI) are used for parallel computation  \cite{GKS-WENO-1}. 
However, the CPU computation is usually limited in the number of threads
which are handled in parallel. Graphics processing unit (GPU) is a
form of hardware acceleration, which is originally developed for
graphics manipulation and execute highly-parallel computing tasks.
Recently,  GPU has applied to HGKS scheme for large-scale scientific
computation \cite{GKS-WENO-5}. In this work, the implicit HGKS
scheme with WENO and HWENO reconstruction on unstructured meshes is
implemented with both CPU using OpenMP and GPU using compute unified
device architecture (CUDA). 

\begin{figure}[!h]
\centering
\includegraphics[width=0.45\textwidth]{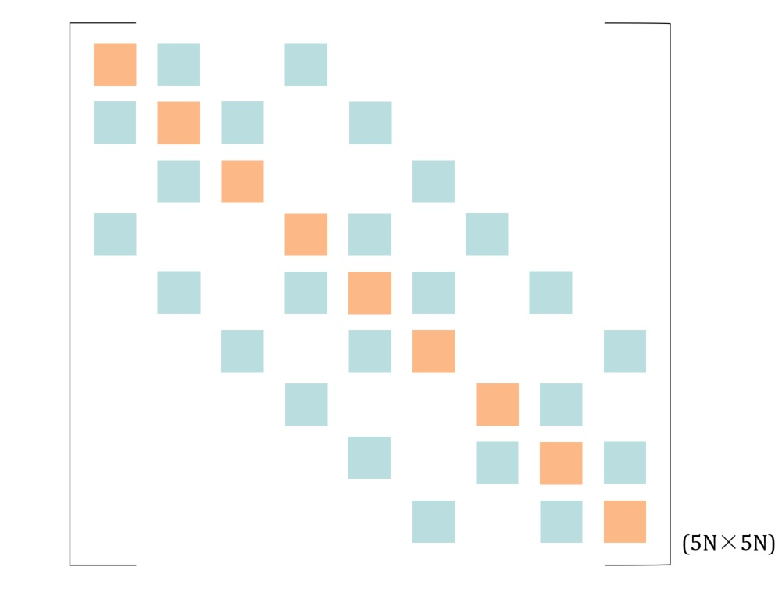}
\caption{\label{jacobian-A} The non-zero block distribution of $\boldsymbol{A}$.}
\end{figure}

The matrix $\boldsymbol{A}$ in Eq.\eqref{A-system-2} is a large asymmetric sparse block matrix of
$(5N, 5N)$ dimension, where $N$ is the total number of computational
cells. The amount of memory storage would be unbearable if the whole
matrix were to be stored. In this section, a simple strategy of storing the matrix $\boldsymbol{A}$
in Eq.\eqref{A-system-2} is introduced to balance the parallel
efficiency and memory storage. It is noticed that the position of all non-zero
blocks in each row (by block) in $\boldsymbol{A}$ only corresponds
to its own cell and all neighboring cells. The non-zero block
distribution of $\boldsymbol{A}$ is shown in Figure.\ref{jacobian-A},
where the orange and blue blocks at the $i$-th row  represent the non-zero
blocks corresponding to cell $\Omega_i$ and its neighboring cells
respectively. These blocks are stored in the structure array for the $i$-th cell. With the storage of 
these non-zero blocks, the computation is easy to be implemented in parallel with single GPU.  
For example, in the computation of $\boldsymbol{r}=\boldsymbol{R}-\boldsymbol{A}\Delta Q$, 
the residual term $\boldsymbol{R}$, i.e., $\mathcal{L}(Q_{i}^{n})$,  can be calculated in parallel naturally, because WENO-type reconstruction and GKS flux solver have high parallelism. 
Meanwhile, denoting the $i$-th row of the matrix $\boldsymbol{A}$ as $\boldsymbol{A}_i=(\boldsymbol{a}_{ij})$, where $a_{ij}\neq 0$ with $j\in \{i\}\cup\{N(i)\}$. 
The vector multiplication of $\boldsymbol{A}_i$ and $\Delta Q$ is equal to $\sum_{j\in \{i\}\cup\{N(i)\}}\boldsymbol{a}_{ij}\Delta Q_j$. 
Such multiplication can be completed through cell structure array, and each row of the matrix term $\boldsymbol{A}\Delta Q$ can be also calculated in parallel. 
In the future, the above strategy is easy to develop to multiple-GPU parallel because matrix partitioning is exactly the same as cell partitioning.

\begin{figure}[!h]
\centering
\includegraphics[width=0.6\textwidth]{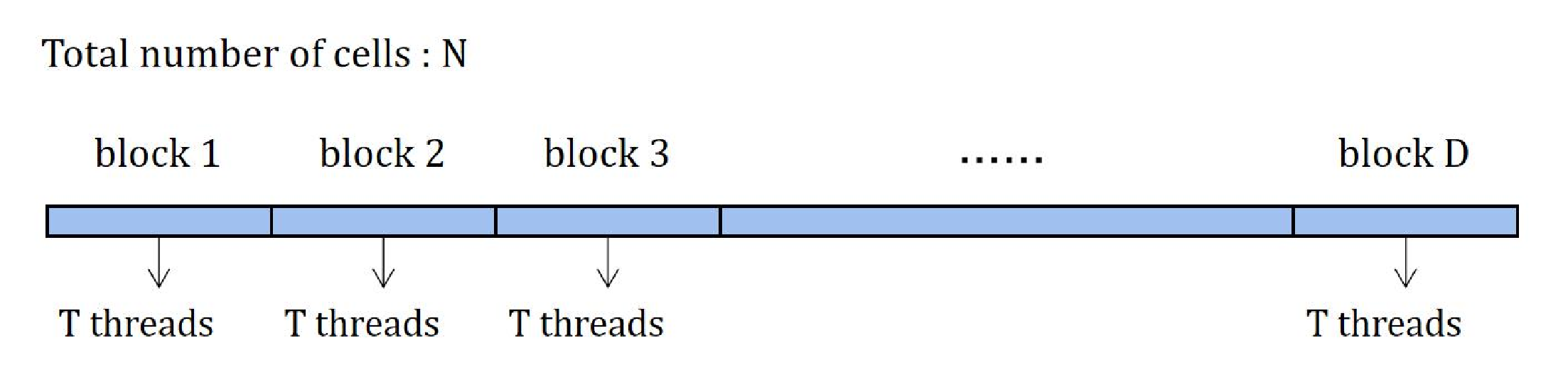}
\caption{\label{thread-block-unstructured} Connection between threads and unstructured cells.}
\end{figure}

For the unstructured meshes, the data for cells, interfaces and
nodes can be stored one-dimensionally. As shown in in
Figure.\ref{thread-block-unstructured}, the cells can be divided
into $D$ blocks, where $N$ is the total number of cells. Meanwhile,
The CUDA threads are organized into thread blocks, and thread blocks
constitute a thread grid. The one-to-one correspondence can be set
for computational cell and thread ID in parallel computation.
In practice, the whole Algorithm.\ref{PGMRES-algorithm} can be
implemented in parallel based on single GPU directly with the above
storage strategy, except Line 23, 24 and 25, i.e., the minimization process in GMRES method. But the matrix
multiplications involved in this minimization process part are all
under 20 dimensions. It does not affect the computational efficiency
in the large-scale computation even if these multiplications are
completed serially on the host CPU.

\section{Numerical tests}
In this section, the numerical tests for both inviscid and viscous flows
will be presented to validate the current scheme. For the inviscid
flows, the collision time $\tau$ takes
\begin{align*}
\tau=\epsilon \Delta t+C\displaystyle|\frac{p_l-p_r}{p_l+p_r}|\Delta t,
\end{align*}
where $\epsilon=0.1$ and $C=1$. For the viscous flows, we have
\begin{align*}
\tau=\frac{\mu}{p}+C \displaystyle|\frac{p_l-p_r}{p_l+p_r}|\Delta t,
\end{align*}
where $p_l$ and $p_r$ denote the pressure on the left and right
sides of the cell interface,  $\mu$ is the dynamic viscous
coefficient and $p$ is the pressure at the cell interface at cell
interface determined by $Q_0$.. In smooth
flow regions, it will reduce to
\begin{equation*}
\tau=\frac{\mu}{p}
\end{equation*}
In the computation, the poly gas is used and the specific heat ratio
takes $\gamma=1.4$.

In order to eliminate the spurious oscillation and improve the stability, the non-compact WENO and compact HWENO reconstructions are performed for the characteristic variables.
For a specific cell interface, denote $Q^*$ and $Q_{\boldsymbol{\tau}}^*$ are the averaged values from both side of cell interface in the local coordinate. $R_*$ is the right eigenmatrix of Jacobian matrix $\big(\displaystyle\partial \textbf{F}/\partial Q\big)$ at $Q=Q^*$, where $\textbf{F}$ is the Euler flux.
The cell averaged conservative variable and cell averaged directional derivatives can be projected into the characteristic field by $q^*=R_*^{-1}Q^*$ and 
$q_{\boldsymbol{\tau}}^*=R_*^{-1}Q_{\boldsymbol{\tau}}^*$.
With the reconstructed values, the conservative variable and directional derivatives can be obtained by the inverse projection.
In the following cases, the GMRES method with non-compact WENO
reconstruction and the GMRES method with compact HWENO
reconstruction are tested. As comparison,  the numerical results of
LUSGS method with non-compact WENO reconstruction are also provided.
For simplicity, the LUSGS, GMRES method with non-compact WENO reconstruction and the GMRES method with compact HWENO reconstruction are denoted as LUSGS-WENO, GMRES-WENO and GMRES-HWENO respectively.
In the computation, the dimension of
Krylov subspace is 10, the restart times of GMRES method is 3 and
the inner iteration times of Jacobi preconditioner is 2. The $L_1$
norm of the density residual is taken as test residual.
The CFL number is set as 10 if there is no special statement.

\begin{figure}[!h]
\centering
\includegraphics[width=0.48\textwidth]{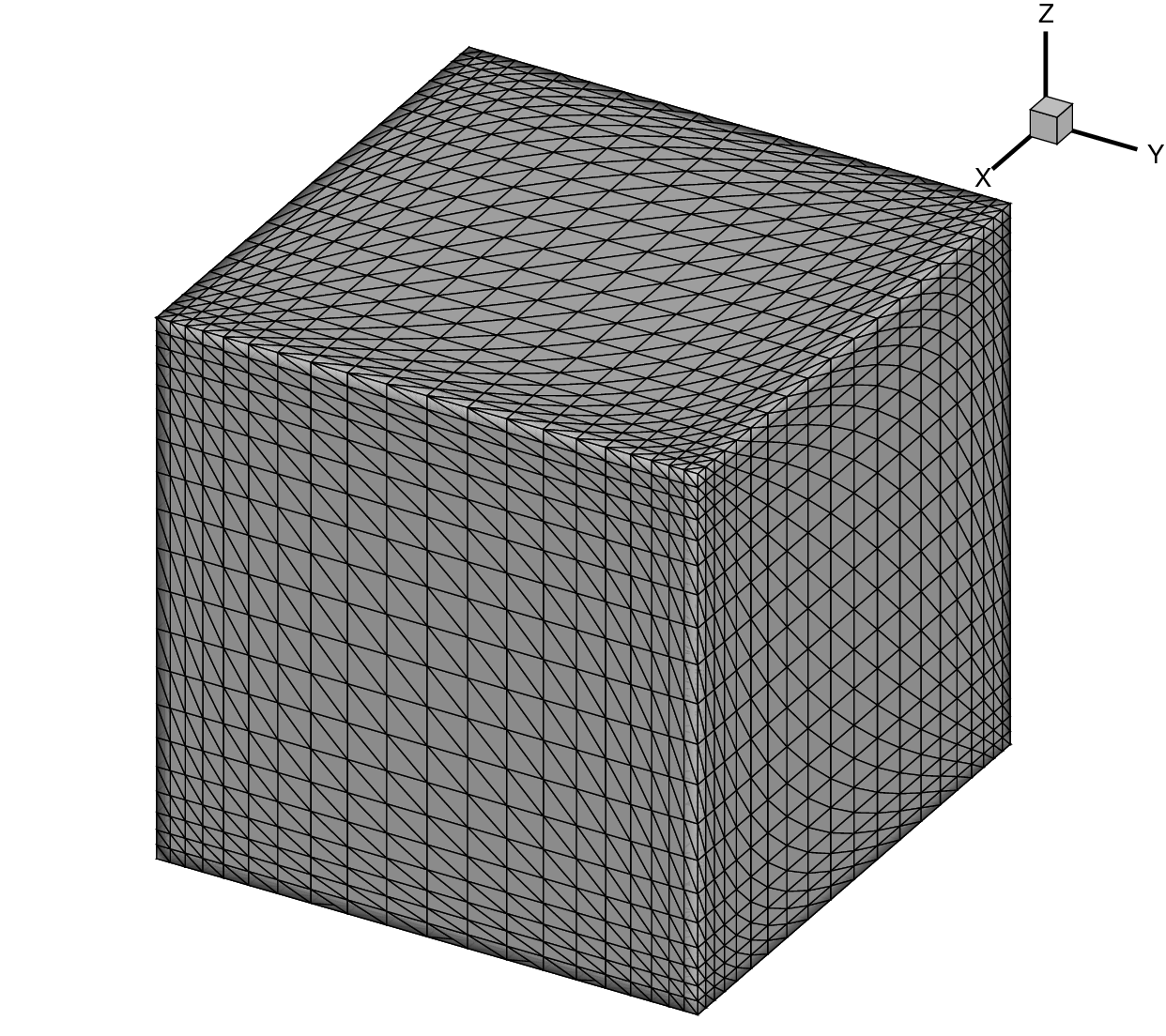}
\caption{\label{cavity-mesh} Lid-driven cavity flow: the local computational
mesh distribution with tetrahedral meshes.}
\end{figure}

\subsection{Lid-driven cavity flow}
The lid-driven cavity problem is one of the most important
benchmarks for numerical Navier-Stokes solvers. The fluid is bounded
by a unit cubic $[0, 1]\times[0, 1]\times[0, 1]$ and driven by a
uniform translation of the top boundary with $y=1$. In this case,
the flow is simulated with Mach number $Ma=0.15$ and all the
boundaries are isothermal and nonslip. Numerical simulations are
conducted with the Reynolds numbers of $Re=1000$ and $100$. This
case is performed by the tetrahedral mesh which contains
$6\times20^3$ cells, where the mesh near the wall is refined and the
size of the first layer cells is $h_{min}=2.5\times 10^{-2}$. The
mesh distribution is shown in Figure.\ref{cavity-mesh}.  
The stencils $S_{i_m}^{HWENO}$ are selected as tetrahedral compact sub-stencils, which are dominated by the gradient informations.

\begin{figure}[!h]
\centering
\includegraphics[width=0.45\textwidth]{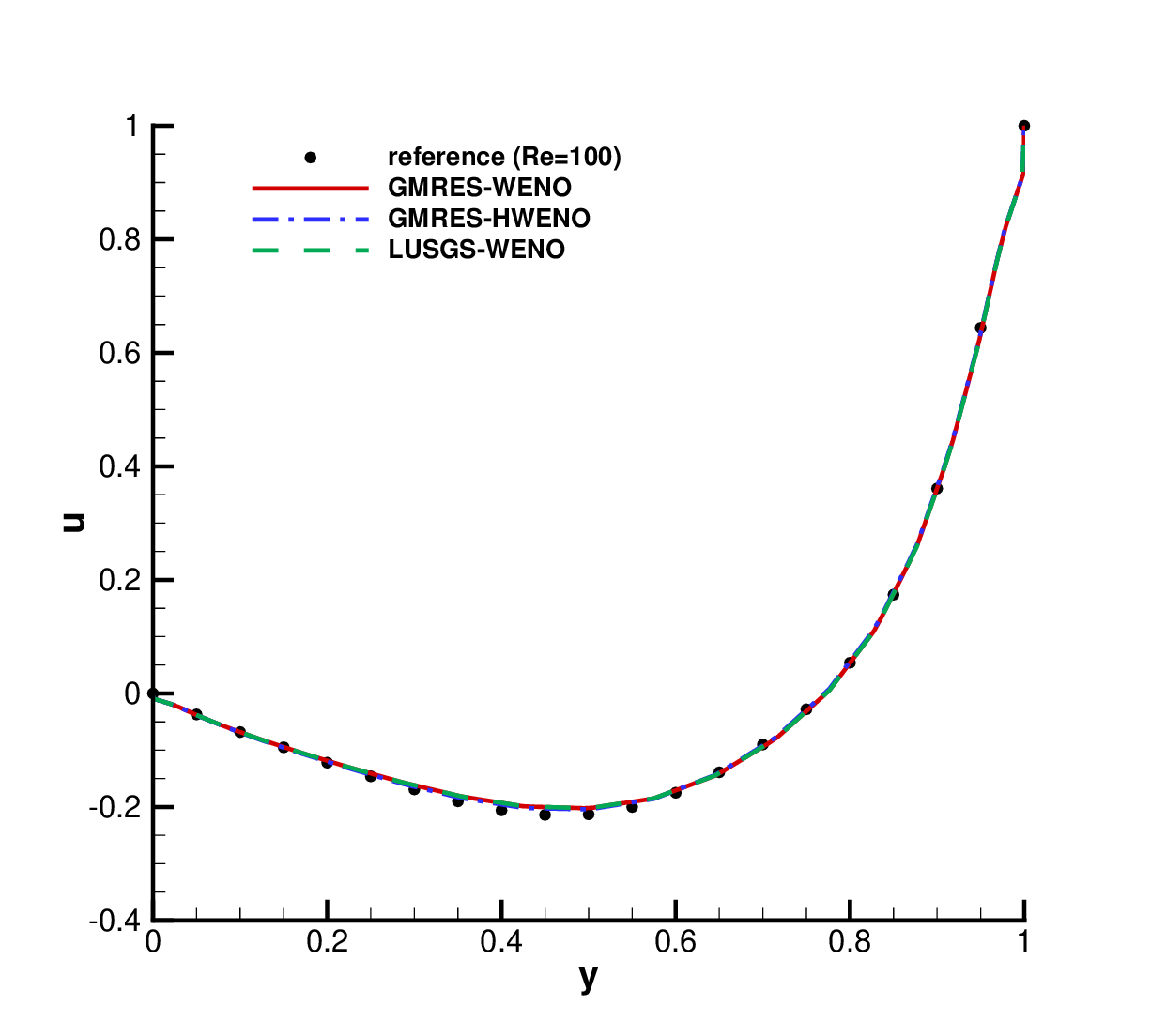}
\includegraphics[width=0.45\textwidth]{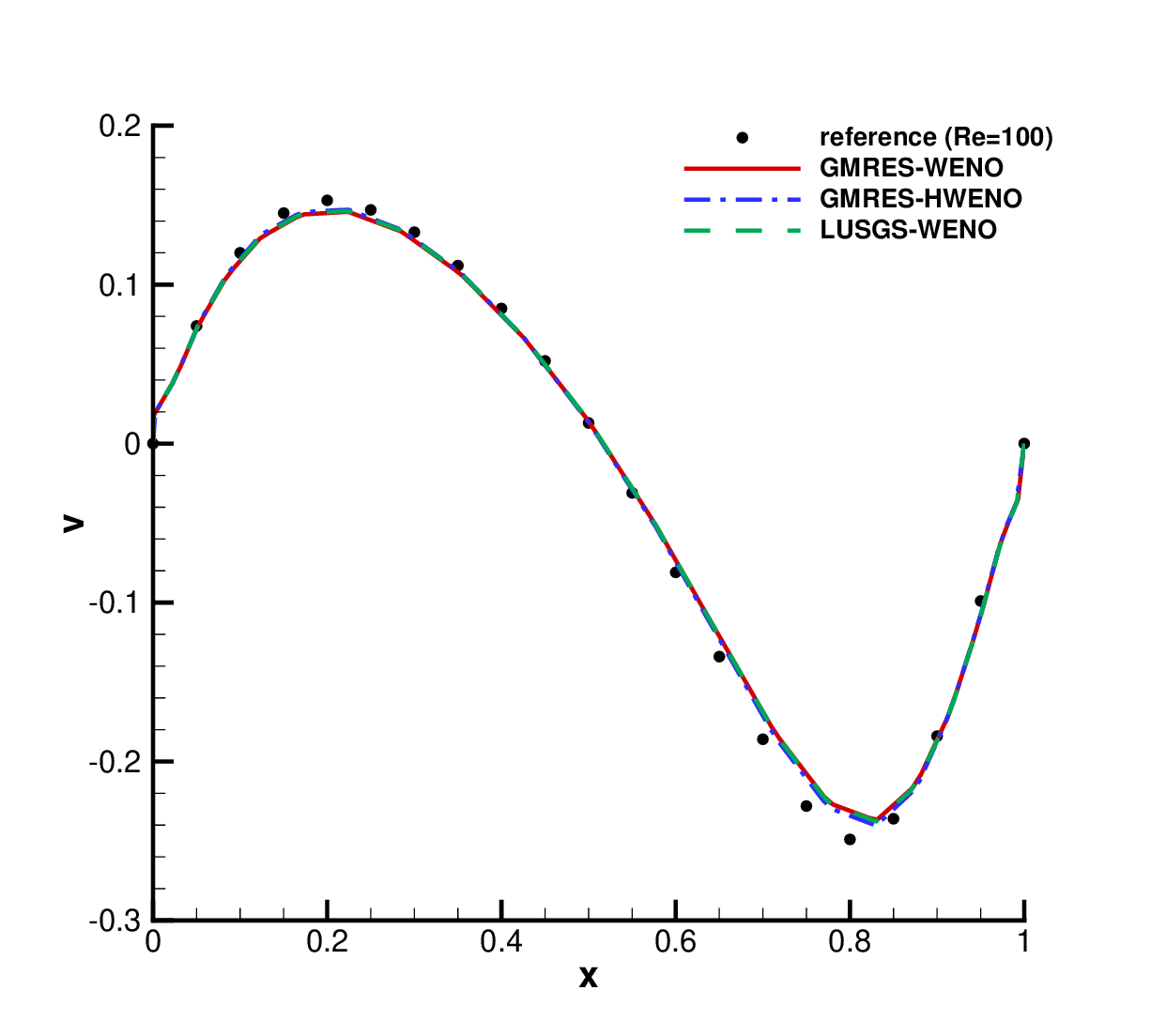}\\
\includegraphics[width=0.45\textwidth]{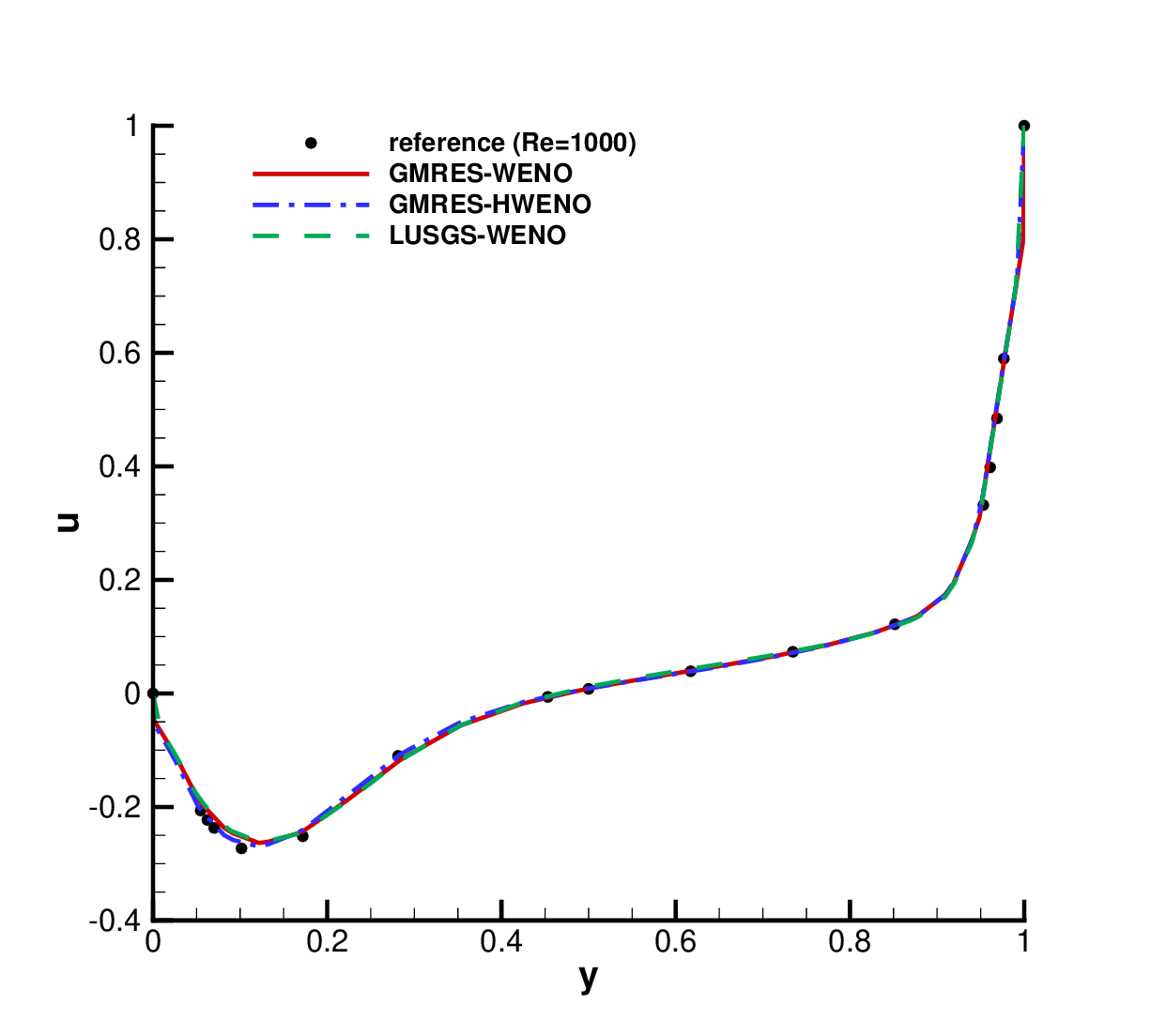}
\includegraphics[width=0.45\textwidth]{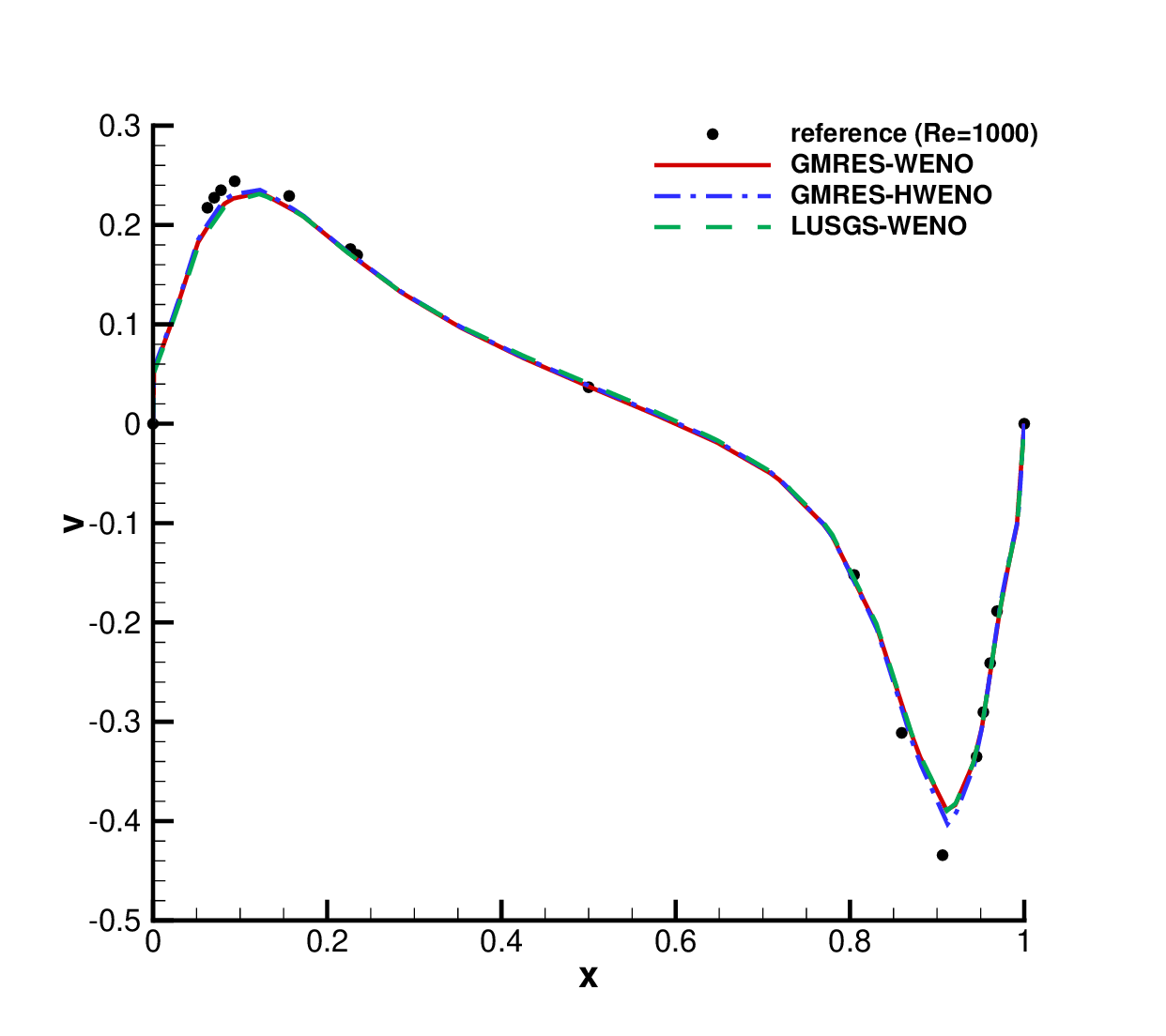}
\caption{\label{cavity-uv} Lid-driven cavity flow: the steady state $U$-velocity profiles along the vertical centerline (left), $V$-velocity
profiles along the horizontal centerline (right) for $Re=1000$ and $100$.}
\end{figure}

\begin{figure}[!h]
\centering
\includegraphics[width=0.48\textwidth]{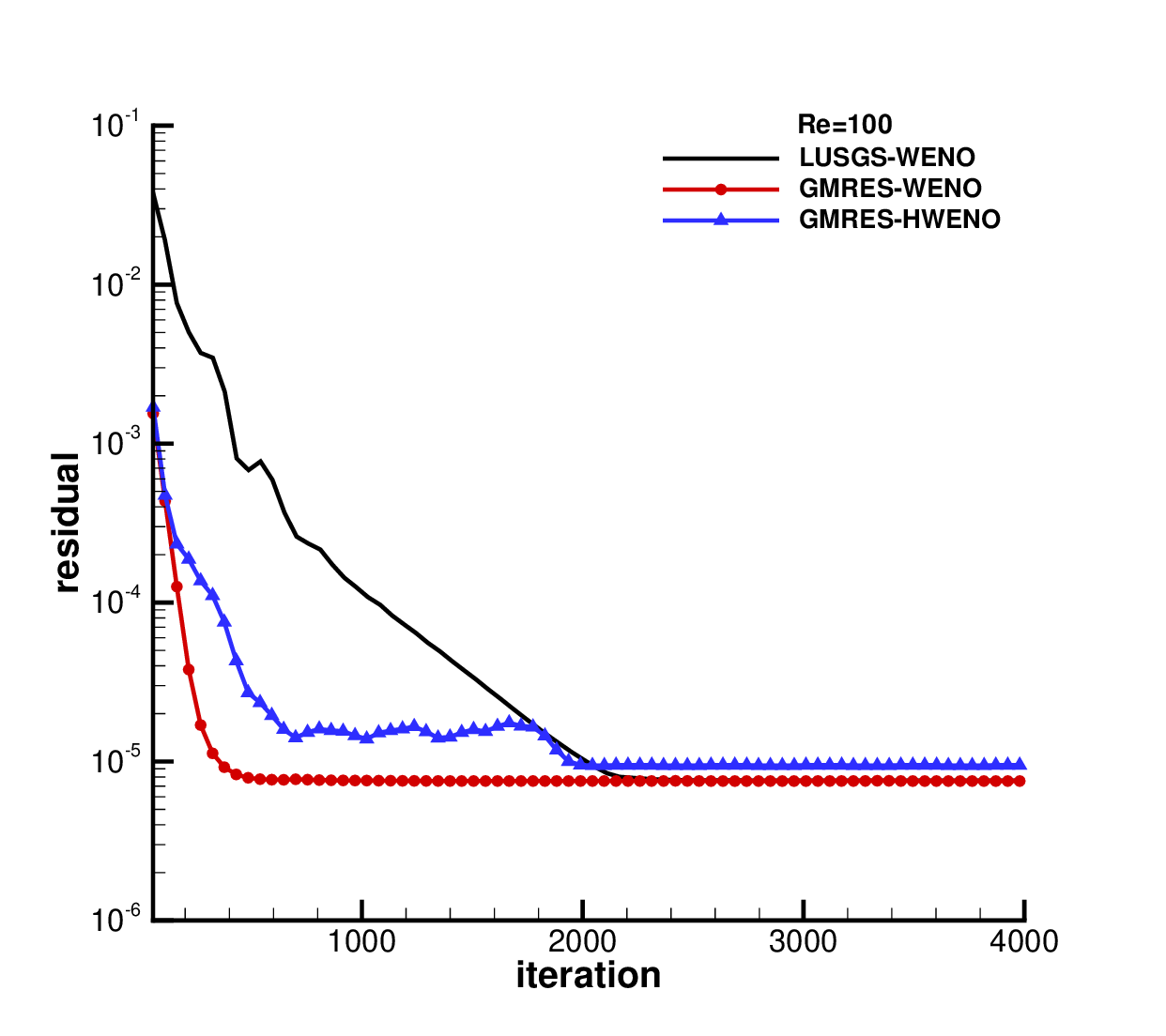}
\includegraphics[width=0.48\textwidth]{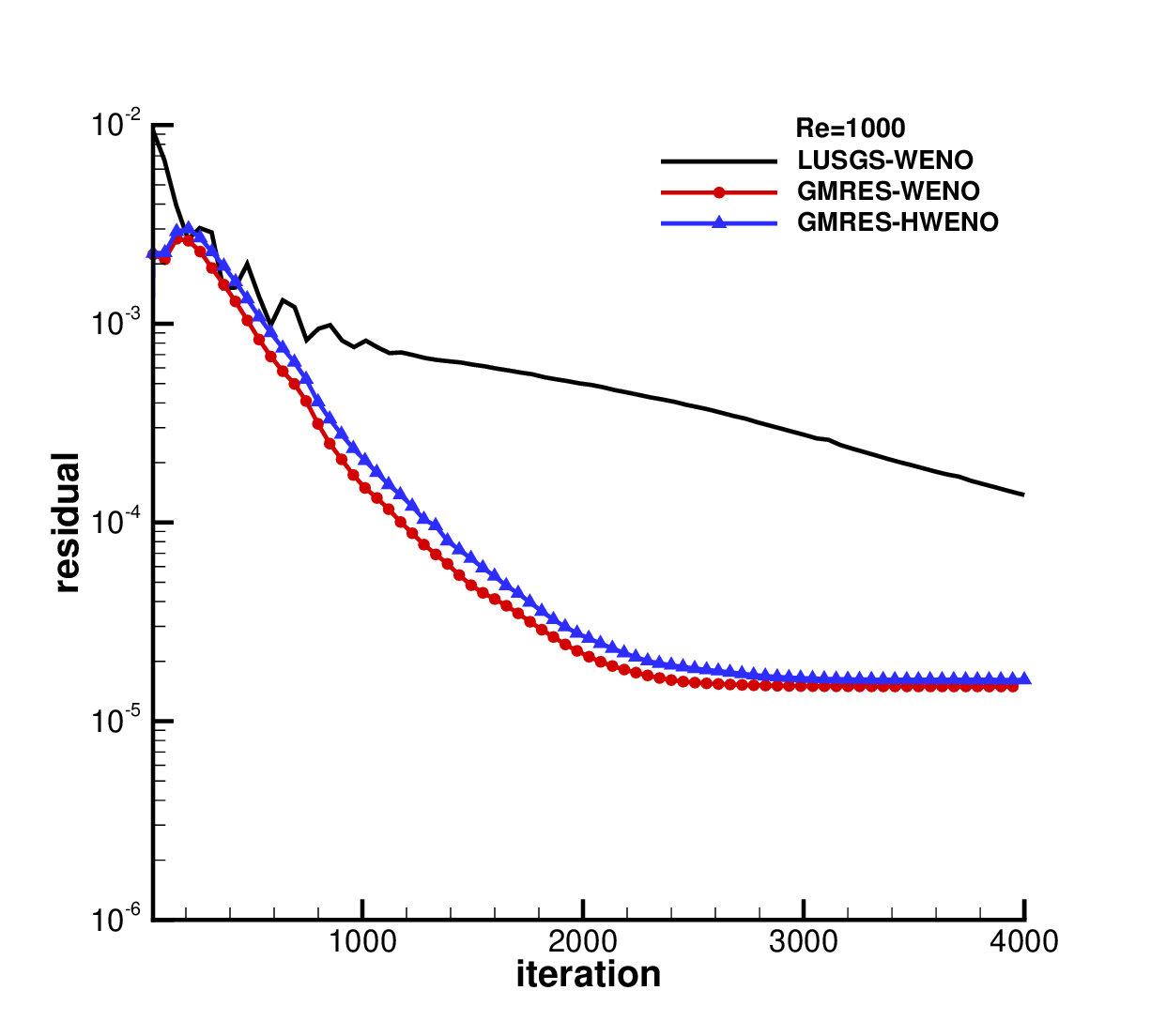}
\caption{\label{cavity-res} Lid-driven cavity flow: the residual comparison with LUSGS-WENO, GMRES-WENO and GMRES-HWENO methods for $Re=100$ and $1000$.}
\end{figure}

The steady state $U$-velocity profiles along the vertical centerline and
$V$-velocity profiles along the horizontal centerline and the
benchmark data for $Re=1000$ \cite{Case-Albensoeder} and $100$
\cite{Case-Shu} are shown in Figure.\ref{cavity-uv}. The numerical
results agree well with the benchmark data, especially the compact
scheme. The histories of residual convergence with different implicit
methods are shown in Figure.\ref{cavity-res} for $Re=1000$ and $100$.
It shows that the residual of GMRES method converges faster than the
LUSGS method. Meanwhile, the residual of compact GMRES method is
comparable with that of non-compact GMRES method.

\begin{figure}[!h]
\centering
\includegraphics[width=0.48\textwidth]{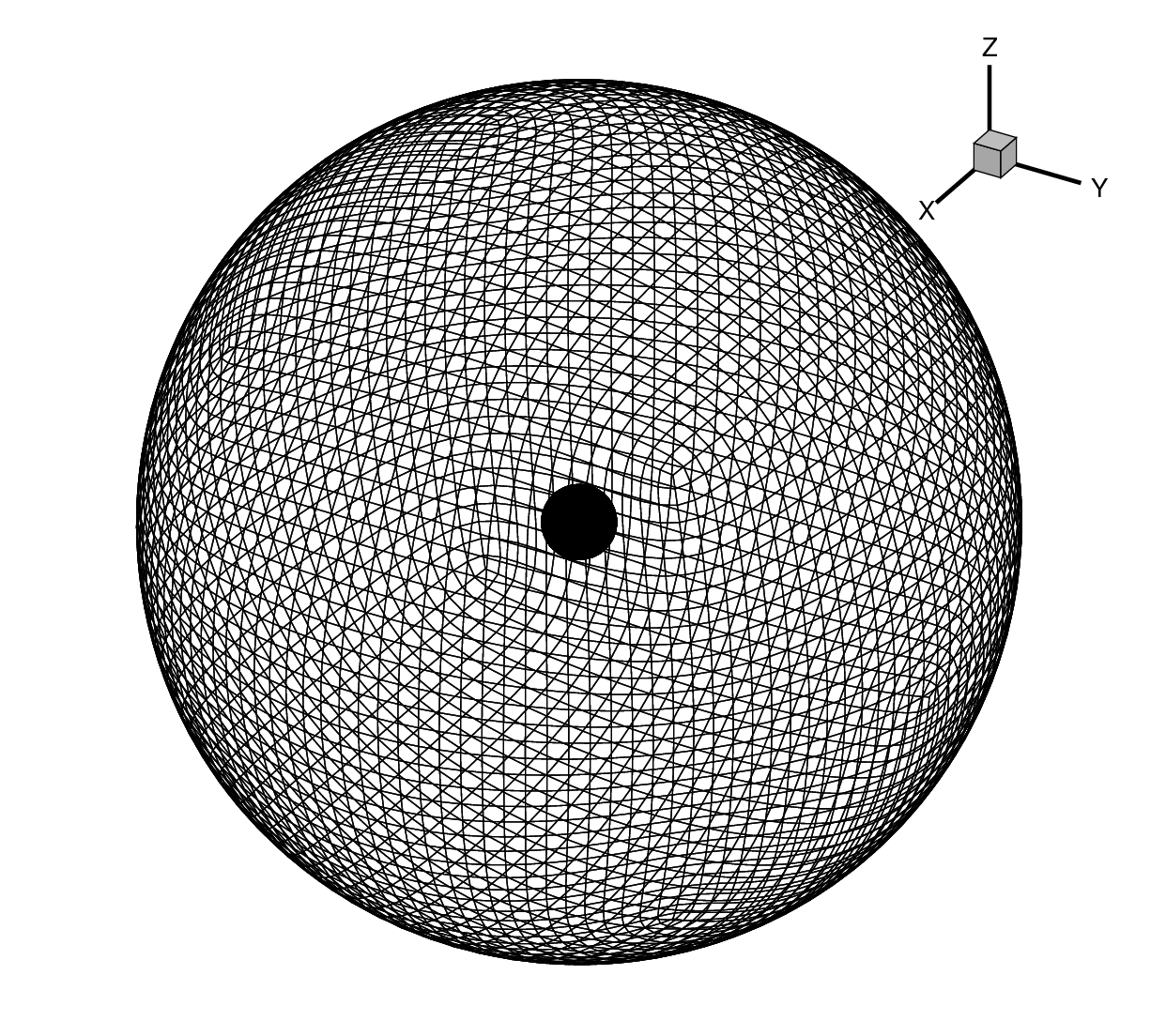}
\includegraphics[width=0.48\textwidth]{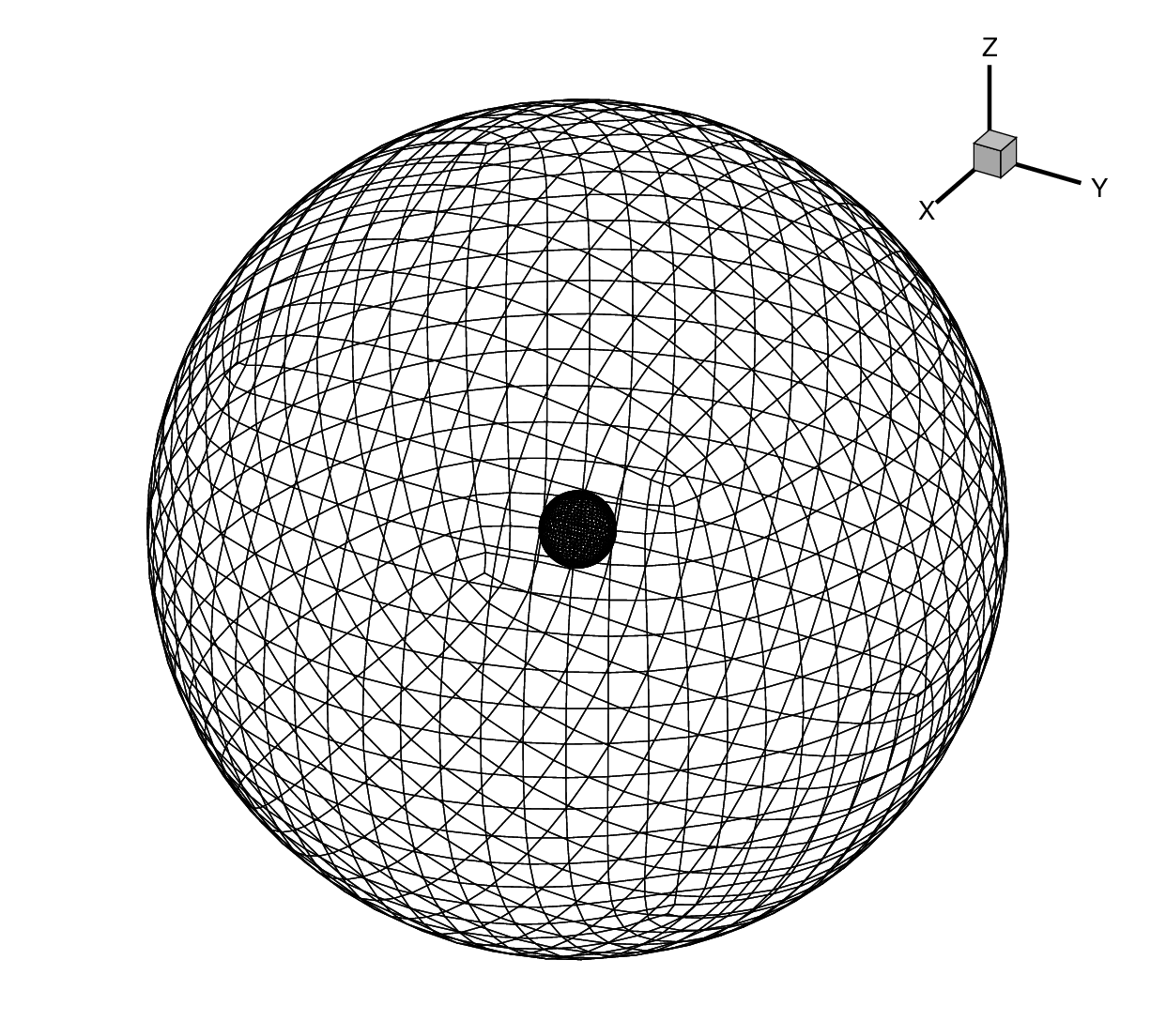}
\caption{\label{sphere-mesh} Flows passing through a sphere: the local computational mesh
distribution with hexahedral meshes for viscous flow (left) and inviscid flow (right).}
\end{figure}

\subsection{Flows passing through a sphere}
This case is used to test the capability in resolving from the low-speed
to supersonic flows, and the initial condition is given as a free
stream condition
\begin{align*}
(\rho,U,V,W,p)_{\infty} = (1, Ma_\infty,0,0, 1/\gamma),
\end{align*}
where $Ma_\infty$ is the Mach number of the free stream. As shown in Figure.\ref{sphere-mesh}, this case is
performed by two unstructured hexahedral meshes. The viscous flows
are tested on the hexahedral mesh which contains $190464$ cells,
where the size of the first layer cells on the surface of sphere is
$h_{min}=3\times 10^{-2}$. The inviscid flows are tested on the
hexahedral mesh which contains $50688$ cells and $h_{min}=2\times
10^{-2}$. The supersonic or subsonic inlet and outlet boundary
conditions are given according to far field normal velocity, the
slip adiabatic boundary condition is used for inviscid flows and the
non-slip adiabatic boundary condition is imposed for viscous flows
on the surface of sphere. 
The stencils $S_{i_m}^{WENO}$ are selected as hexahedral compact sub-stencils, which are dominated by the conservative variables informations.

\begin{figure}[!h]
\centering
\includegraphics[width=0.435\textwidth]{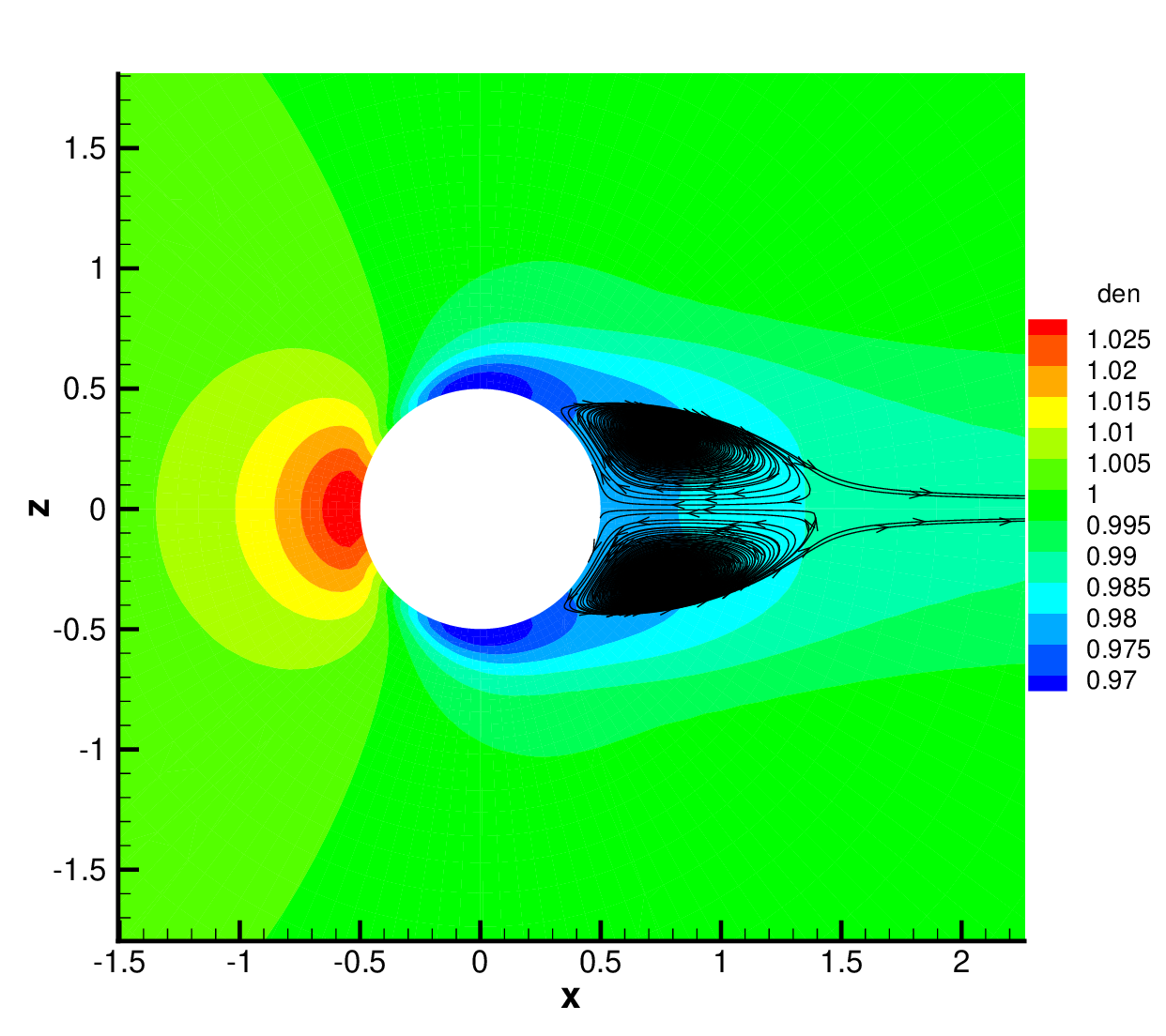}
\includegraphics[width=0.435\textwidth]{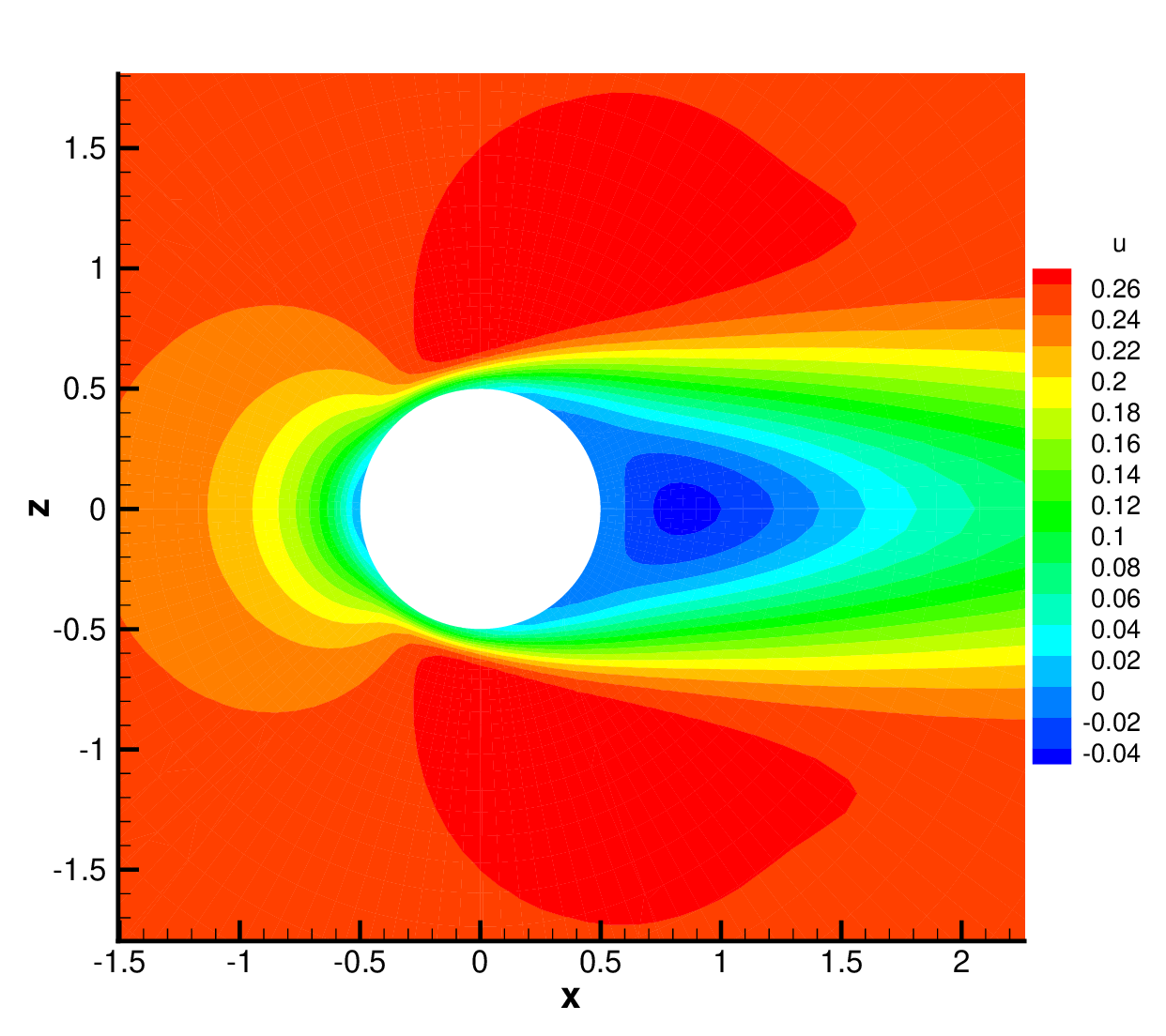}\\
\includegraphics[width=0.435\textwidth]{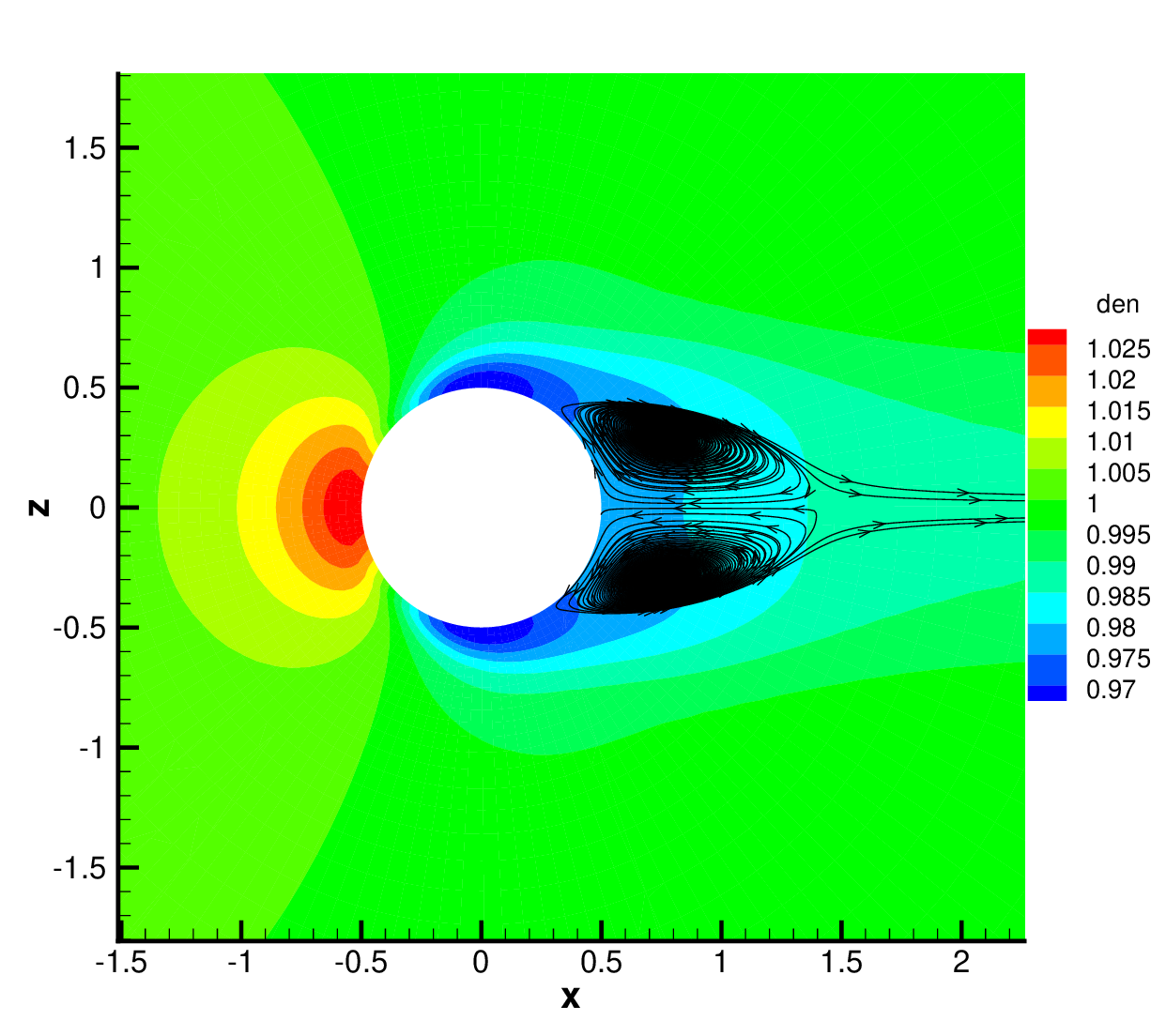}
\includegraphics[width=0.435\textwidth]{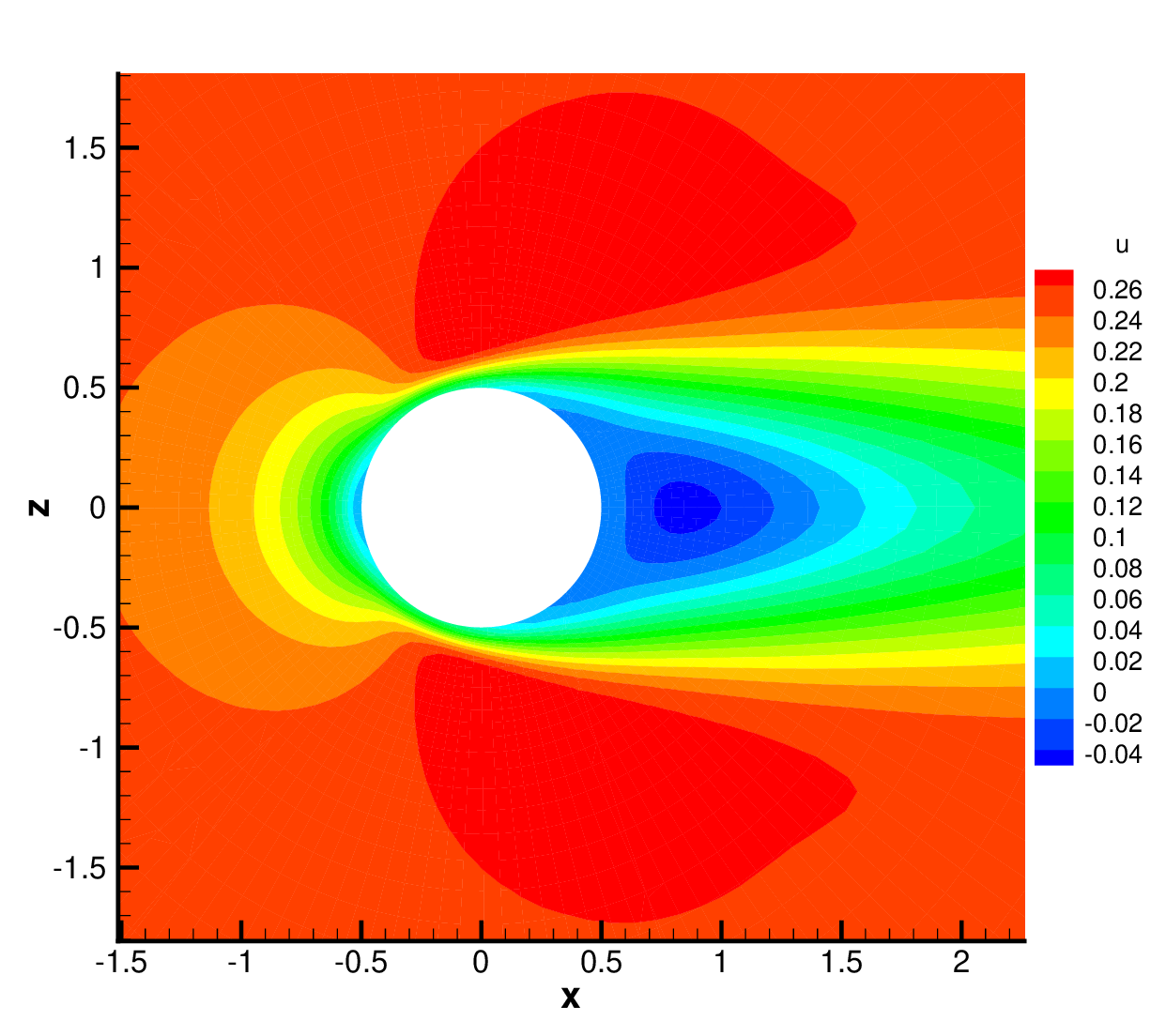}\\
\includegraphics[width=0.435\textwidth]{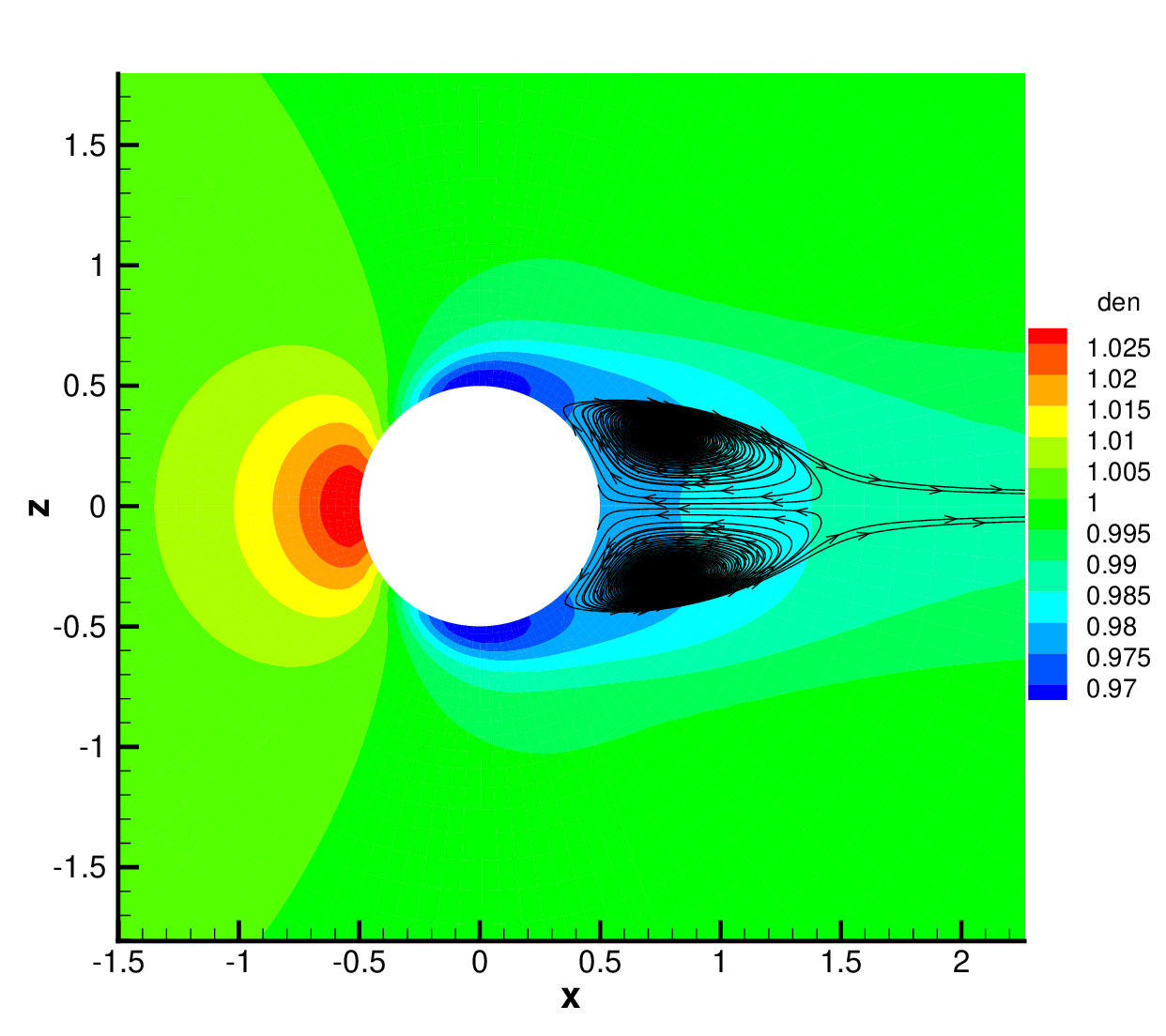}
\includegraphics[width=0.435\textwidth]{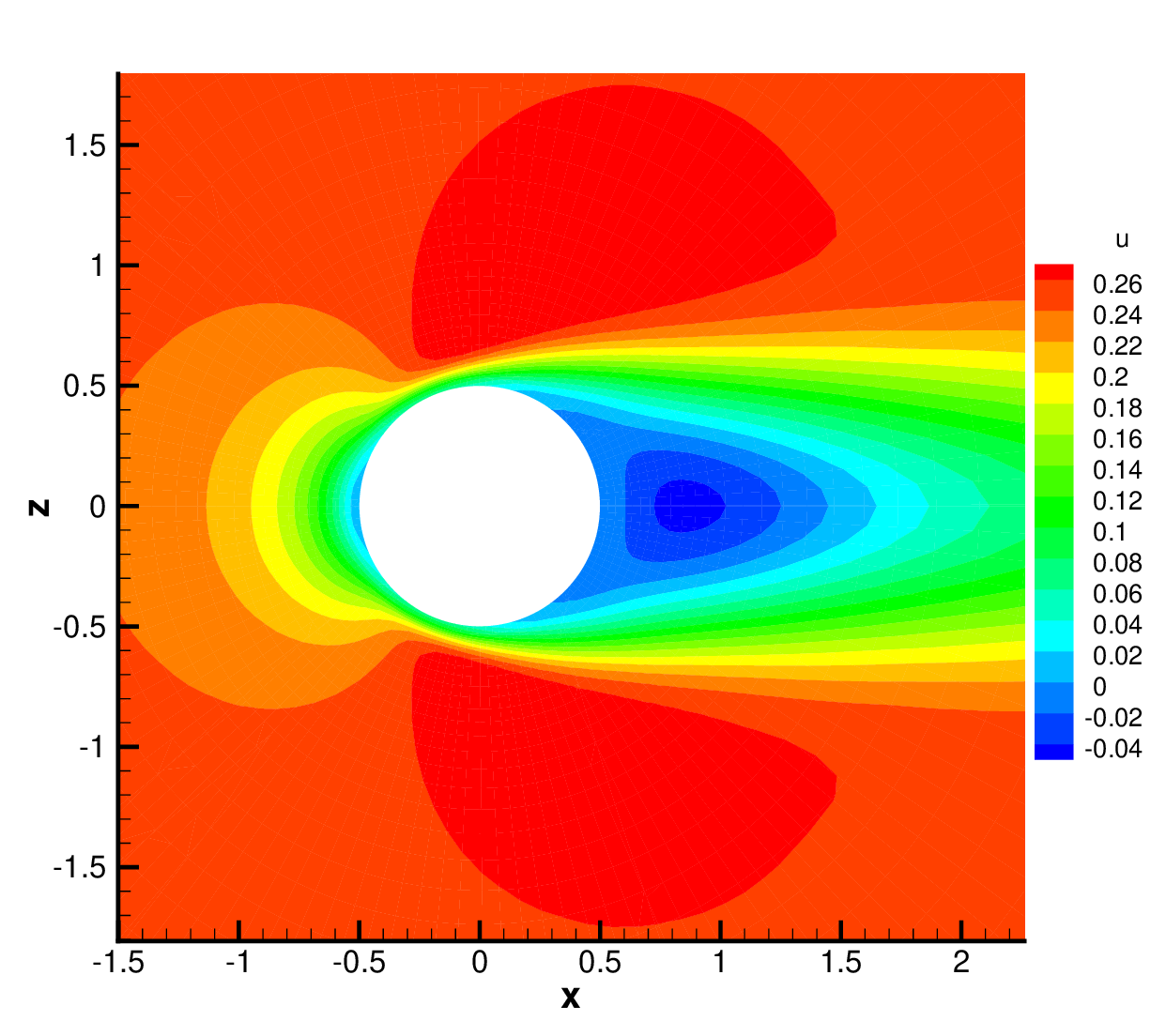}
\caption{\label{sphere-2535} Subsonic viscous flow passing through a sphere: the density and
streamline distributions (left) and velocity distribution (right) at vertical centerline planes with LUSGS-WENO method (top), GMRES-WENO method (middle) and GMRES-HWENO method (bottom) for $Re=118$ and $Ma_\infty=0.2535$. }
\end{figure}

\begin{figure}[!h]
\centering
\includegraphics[width=0.435\textwidth]{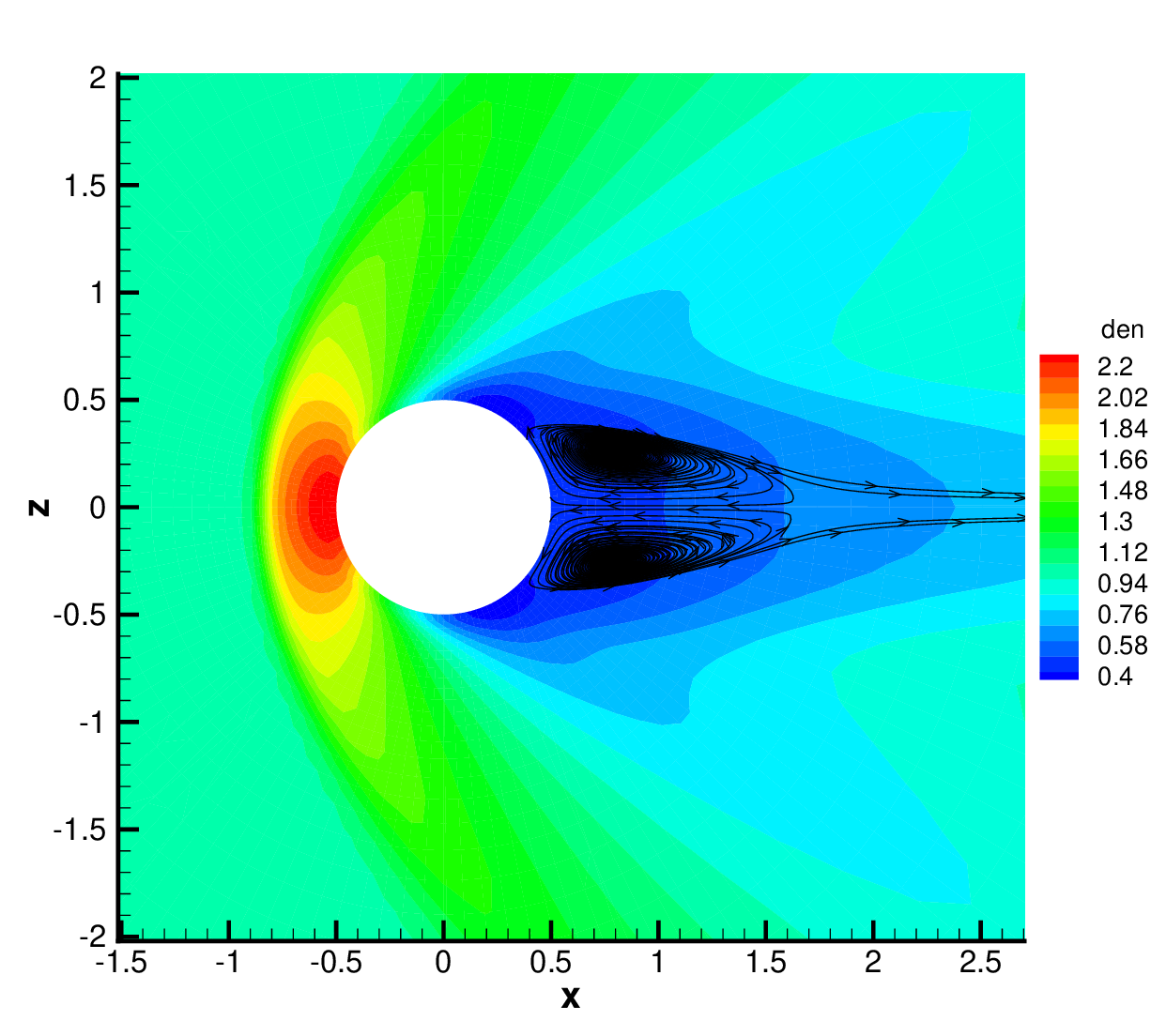}
\includegraphics[width=0.435\textwidth]{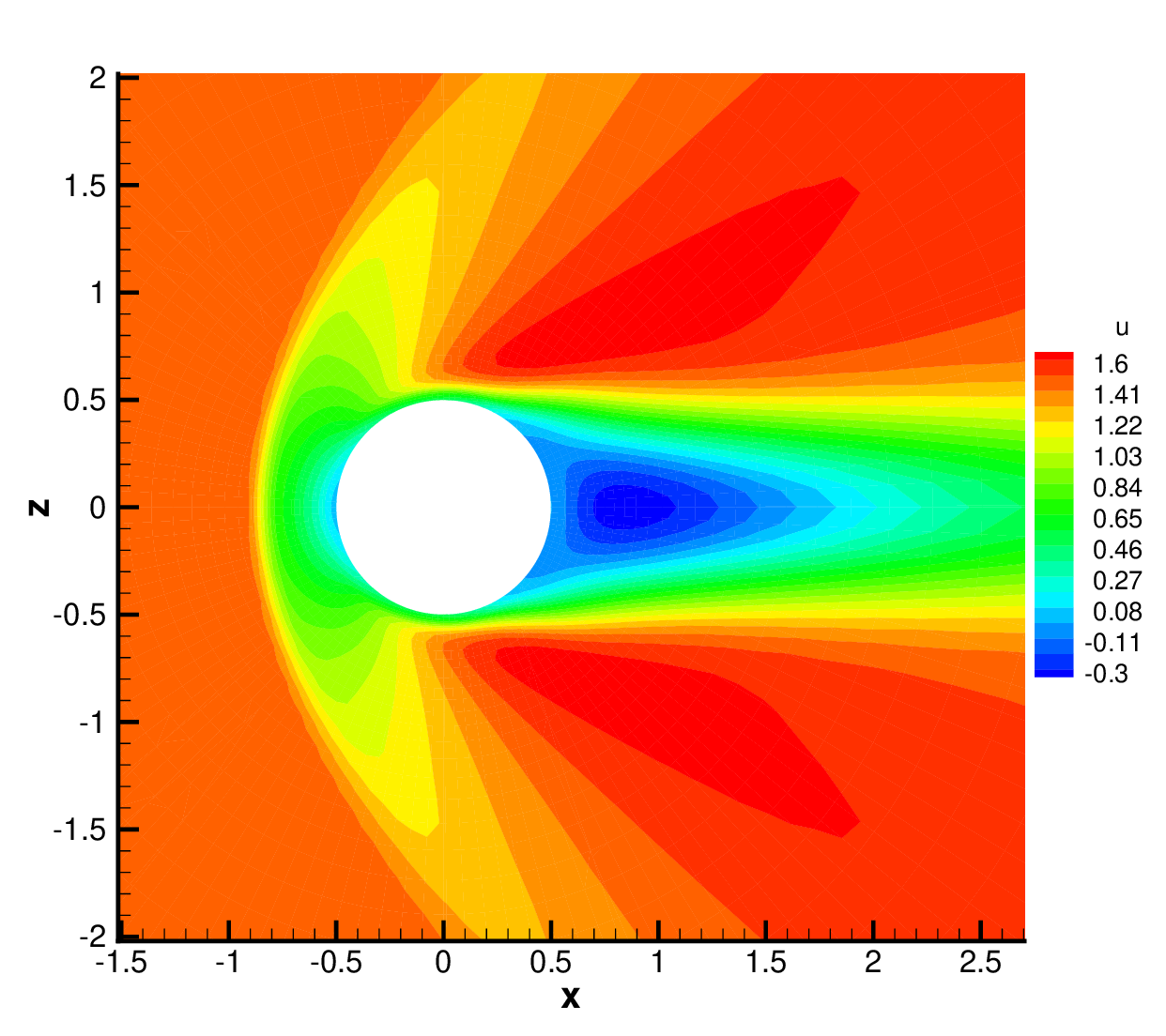}\\
\includegraphics[width=0.435\textwidth]{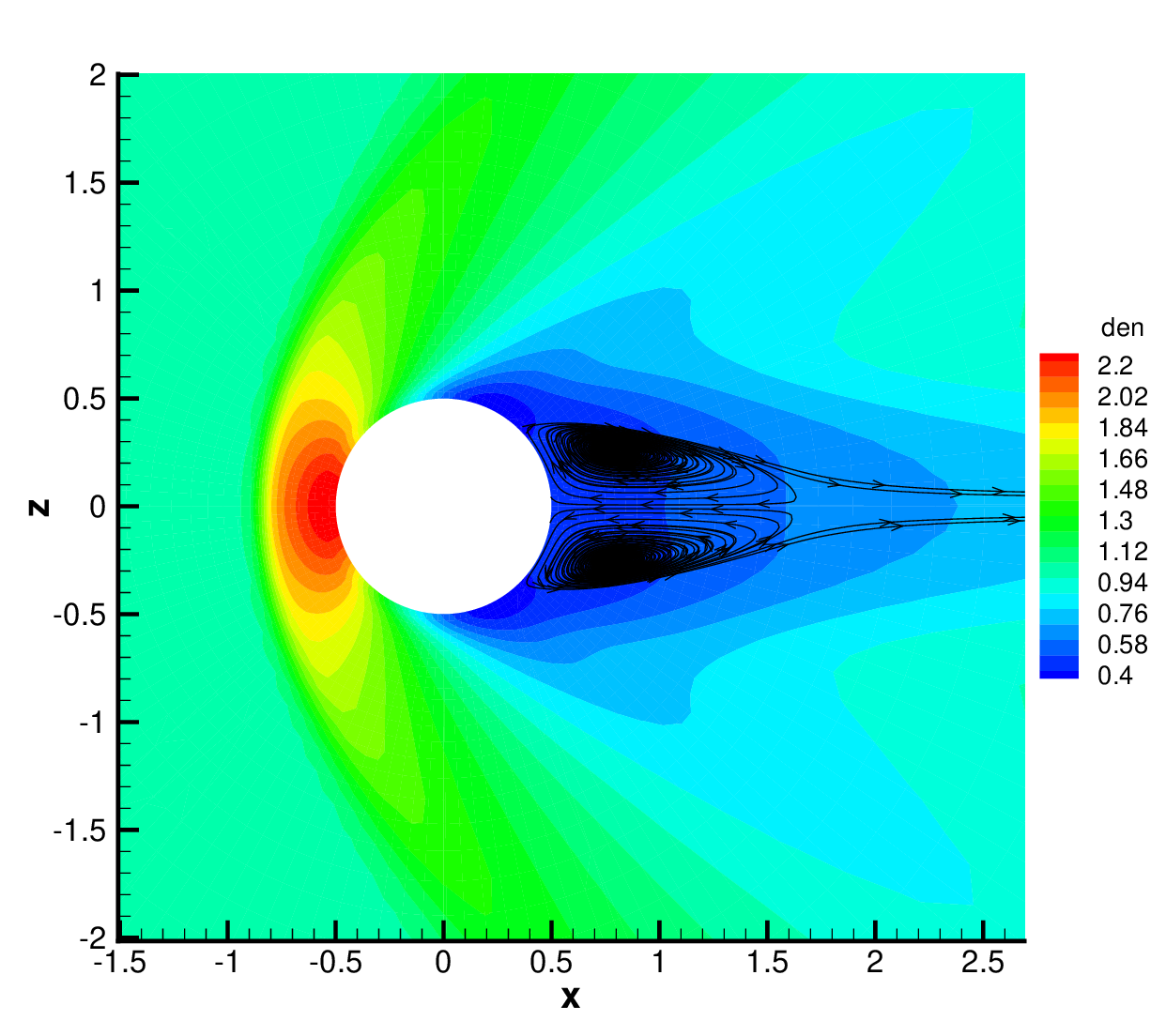}
\includegraphics[width=0.435\textwidth]{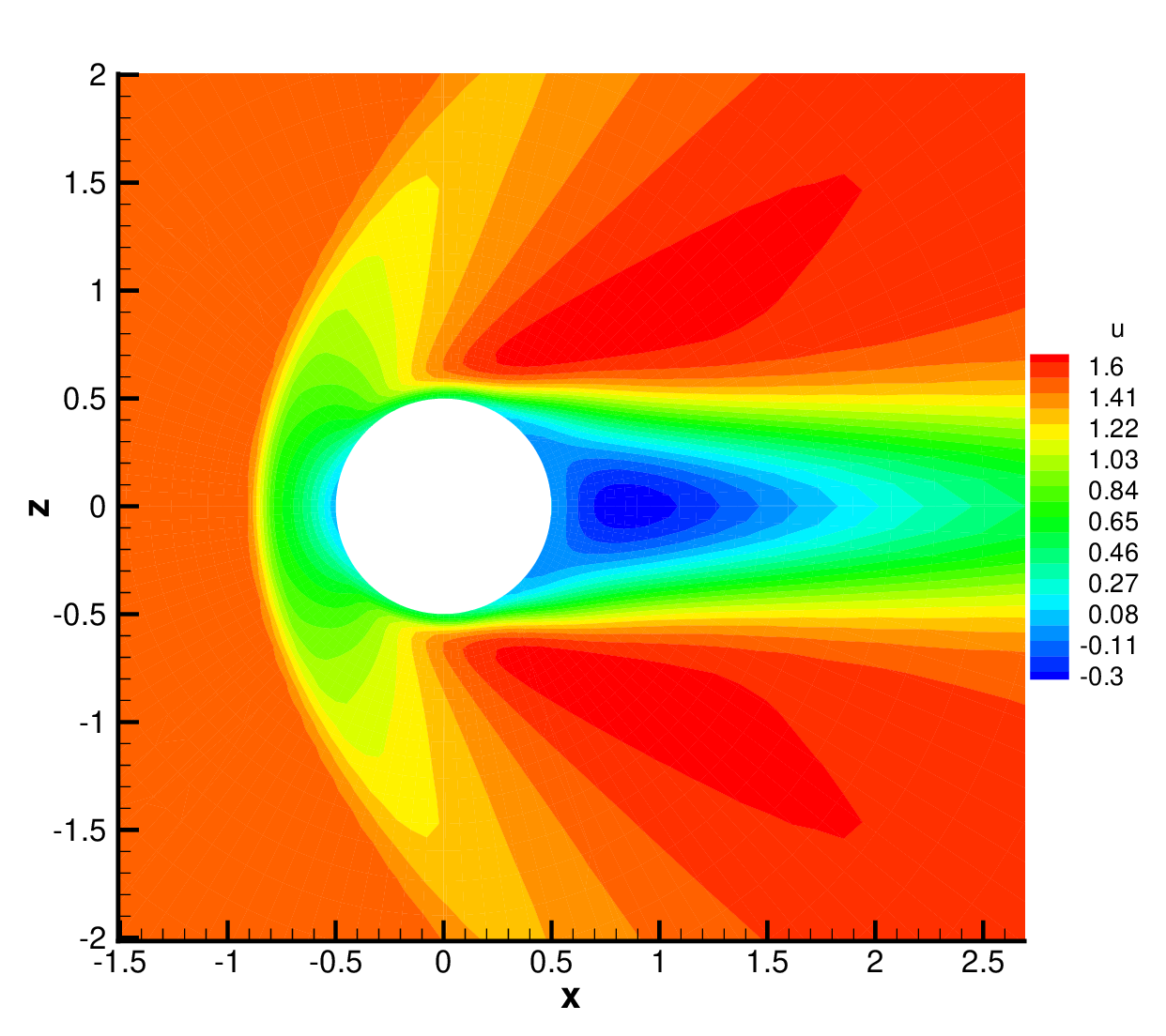}\\
\includegraphics[width=0.435\textwidth]{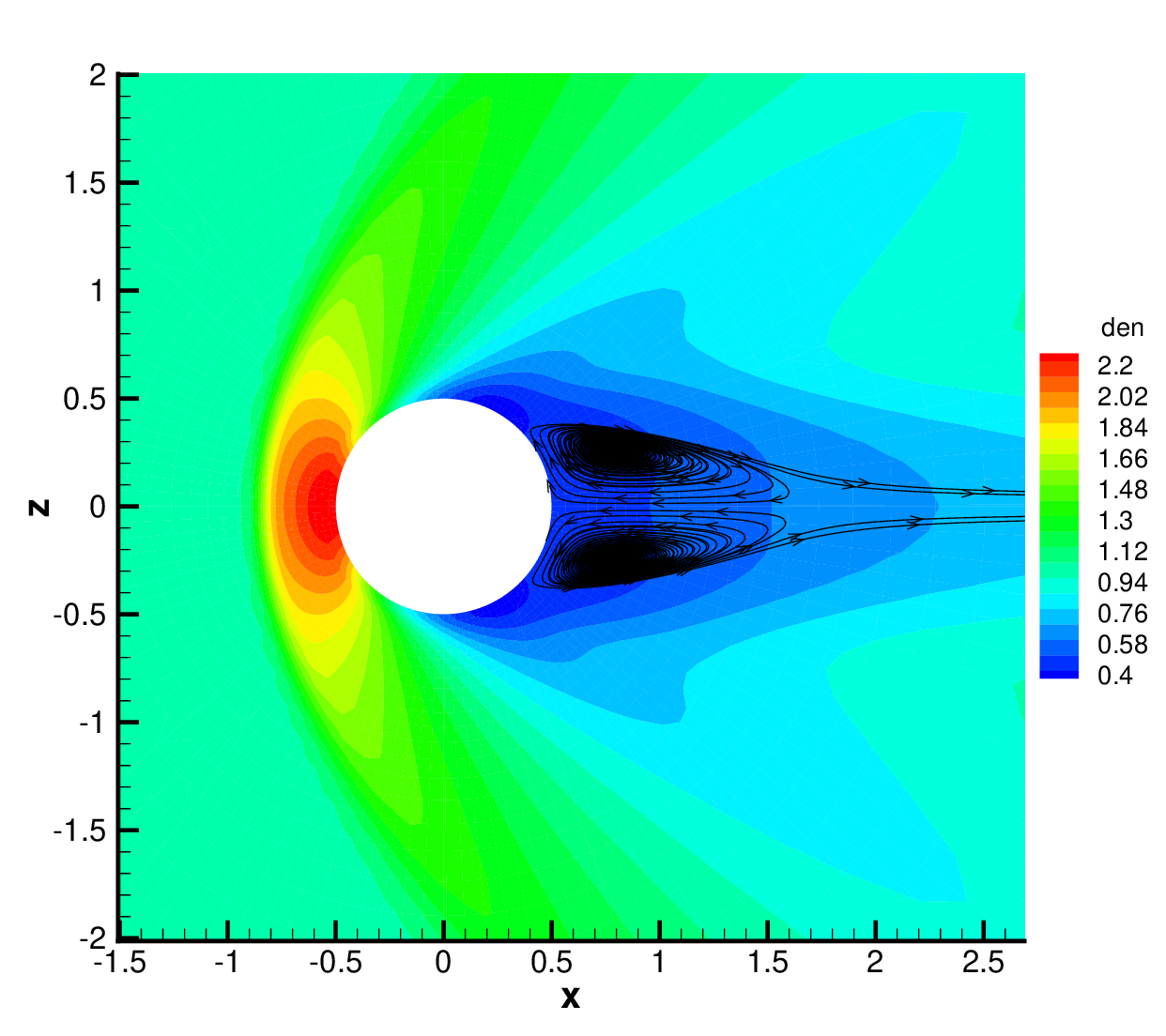}
\includegraphics[width=0.435\textwidth]{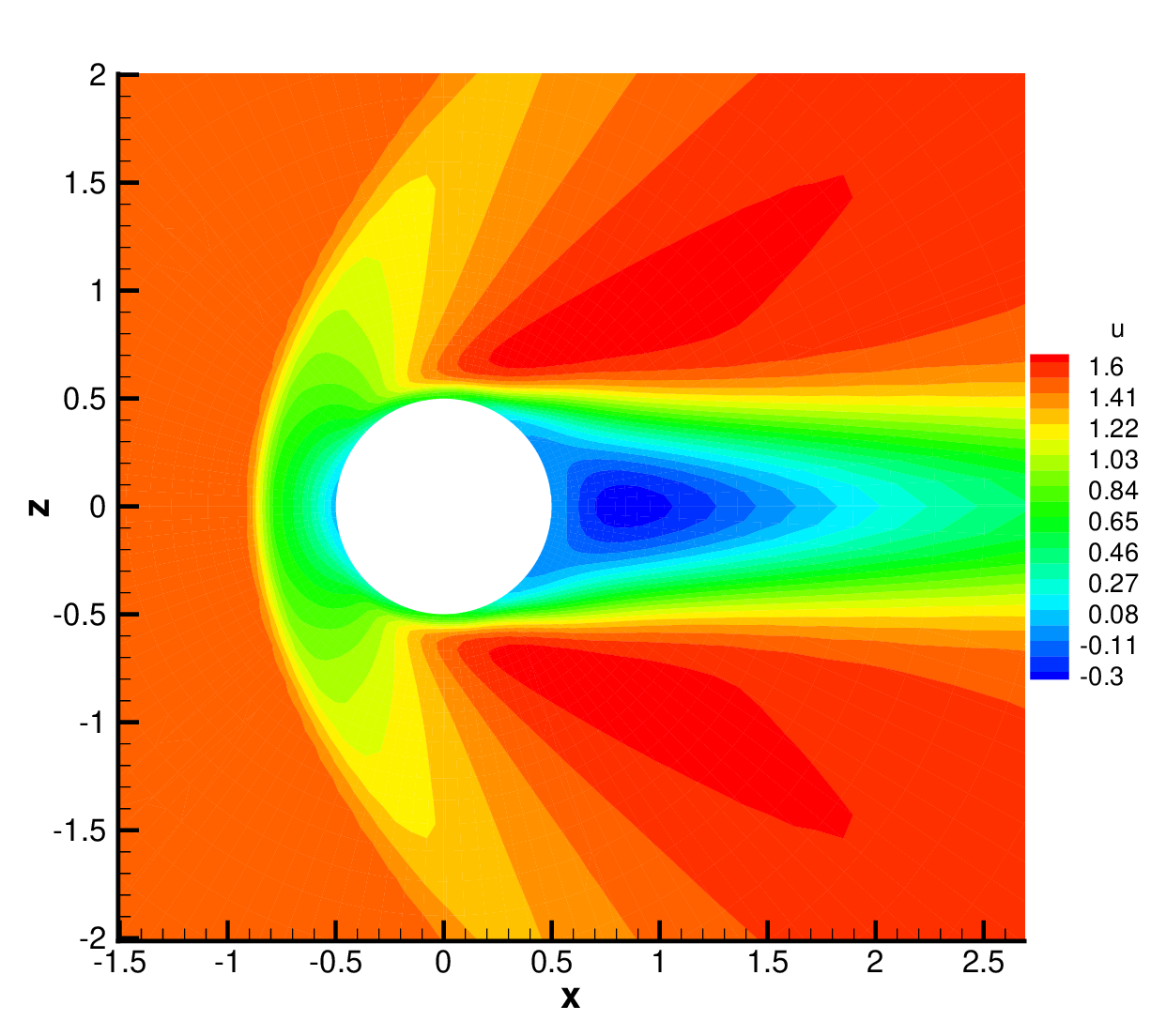}
\caption{\label{sphere-150} Supersonic viscous flow passing through a sphere: the density and streamline
distributions (left) and velocity distribution (right) at vertical centerline planes with LUSGS-WENO method (top), GMRES-WENO method (middle) and GMRES-HWENO method (bottom) for $Re=300$ and $Ma_\infty=1.5$. }
\end{figure}

\begin{table}[!h]
\begin{center}
\def\temptablewidth{0.86\textwidth}{\rule{\temptablewidth}{0.9pt}}
\begin{tabular*}{\temptablewidth}{@{\extracolsep{\fill}}c|c|c|c}
Scheme & Computational Mesh &  $L$ & $\theta$ \\
\hline
LUSGS-WENO GKS                                                 & 190464  cells, Hex        &   0.91  & 124.5\\
GMRES-WENO GKS                           & 190464  cells, Hex        &  0.91  & 124.5\\
GMRES-HWENO GKS                                 & 190464  cells, Hex        &  0.95  & 128.1\\
4th-order DDG GMRES \cite{Case-Cheng}  & 1608680 cells, Hybrid  &  0.96   & 123.7 \\
Experiment \cite{Case-Taneda}                    & -                                    & 1.07    & 151\\
\end{tabular*}
{\rule{\temptablewidth}{1.0pt}}
\end{center}
\caption{\label{sphere-2535-table} Subsonic inviscid flow passing through a sphere:
quantitative comparisons of closed wake length $L$ and separation angle $\theta$ for $Re=118$ and $Ma_\infty=0.2535$.}
\begin{center}
\def\temptablewidth{0.8\textwidth}{\rule{\temptablewidth}{1.0pt}}
\begin{tabular*}{\temptablewidth}{@{\extracolsep{\fill}}c|c|c|c}
Scheme & Computational Mesh &  $L$ & $\theta$\\
\hline
LUSGS-WENO GKS                                             & 190464 cells, Hex & 1.17 & 135.1\\
GMRES-WENO GKS                      & 190464 cells, Hex & 1.17  & 135.4\\
GMRES-HWENO GKS                             & 190464 cells, Hex & 1.16   & 135.9\\
WENO-CU6-FP RK3 \cite{Case-Nagata} & 909072 cells, Hex & 0.96 & 137.2 \\
\end{tabular*}
{\rule{\temptablewidth}{1.0pt}}
\end{center}
\caption{\label{sphere-150-table} Supersonic viscous flow passing through a sphere:
quantitative comparisons of closed wake length $L$ and separation angle $\theta$ for $Re=300$ and $Ma_\infty=1.5$.}
\end{table}

\begin{figure}[!h]
\centering
\includegraphics[width=0.48\textwidth]{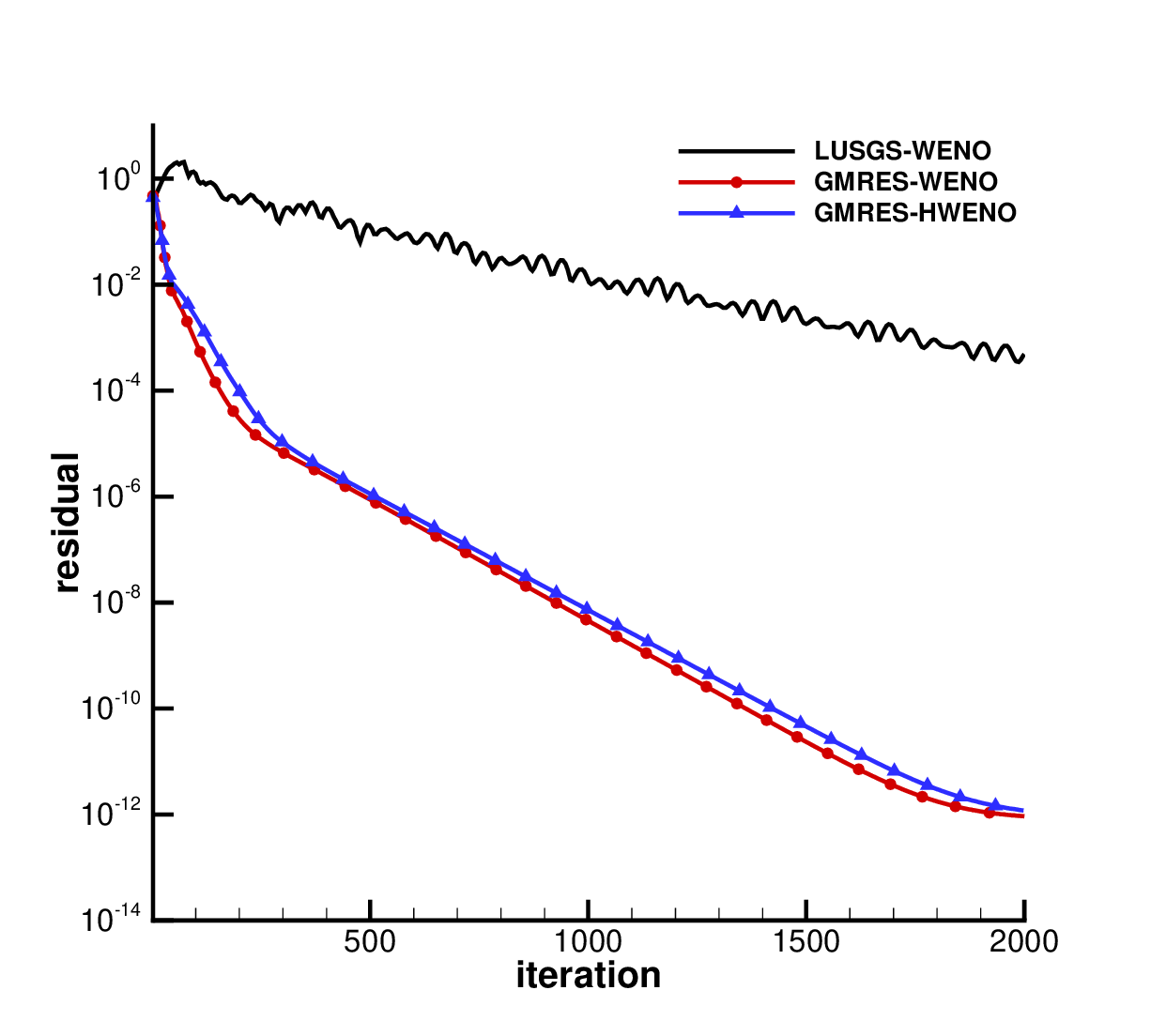}
\includegraphics[width=0.48\textwidth]{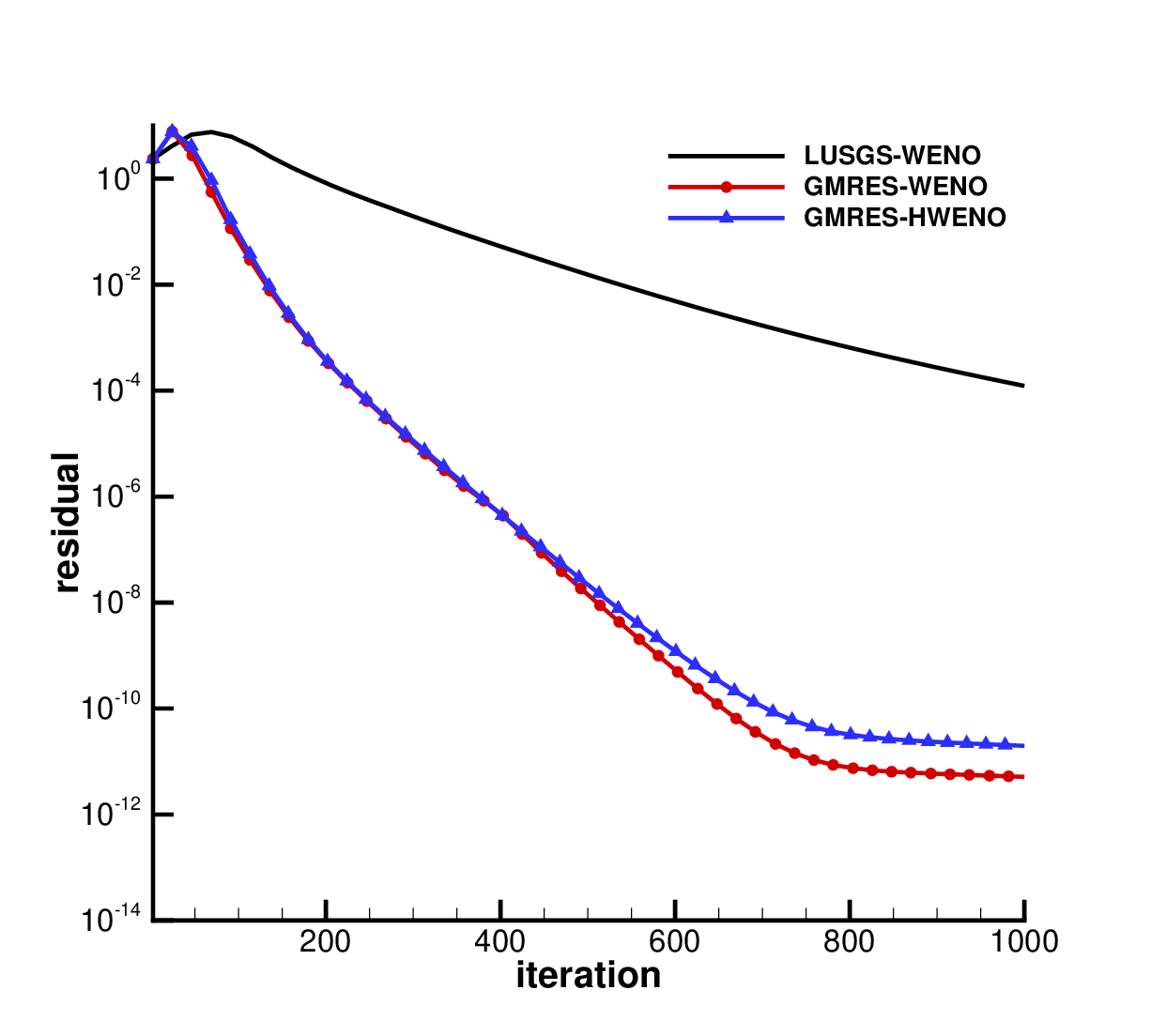}
\caption{\label{sphere-res-2535-150} Subsonic viscous flow passing through a sphere: the residual convergence comparison with LUSGS-WENO, GMRES-WENO and GMRES-HWENO methods for $Re=118$ and $Ma_\infty=0.2535$ (left) and $Re=300$ and $Ma_\infty=1.5$ (right).}
\end{figure}

To validate the linear stability of the implicit methods, the
subsonic case with $Re=118$ and $Ma_\infty=0.2535$ is provided, and
the linear weights are used in both WENO and HWENO reconstruction. In this case, the CFL number is taken as 20.
The density, velocity and streamline distributions at vertical
centerline planes are shown in Figure.\ref{sphere-2535}. The
quantitative results of separation angle $\theta$ and closed wake
length $L$ are given in Table.\ref{sphere-2535-table}. The compact
GMRES method with linear weights gives the better values of $L$ and
$\theta$, which are closer to the experiment data. It means that the
capability of implicit compact scheme to capture the structure of
viscous flow is more accurate. To verify the robustness of the
implicit methods, the supersonic viscous flow with $Re=300$ and
$Ma_\infty=1.5$ is tested as well. The non-linear weights are used
in both WENO and HWENO reconstruction. 
In this case, the CFL number is taken as 5.
The density, velocity and
streamline distributions at vertical centerline planes are shown in
Fig.\ref{sphere-150}. The quantitative results of closed wake length
$L$ and separation angle $\theta$ are given in
Table.\ref{sphere-150-table}. The current implicit compact scheme
gives a bigger value of $\theta$.  The performance
of closed wake length $L$ would be better with a refined mesh. The
residual convergence history with different implicit methods is
given in Figure.\ref{sphere-res-2535-150}. It also shows that the GMERS
method converges much faster than the LUSGS method, and the steady
state residual of GMRES method is at least five orders of magnitude
smaller than that of LUSGS method.

\begin{figure}[!h]
\centering
\includegraphics[width=0.435\textwidth]{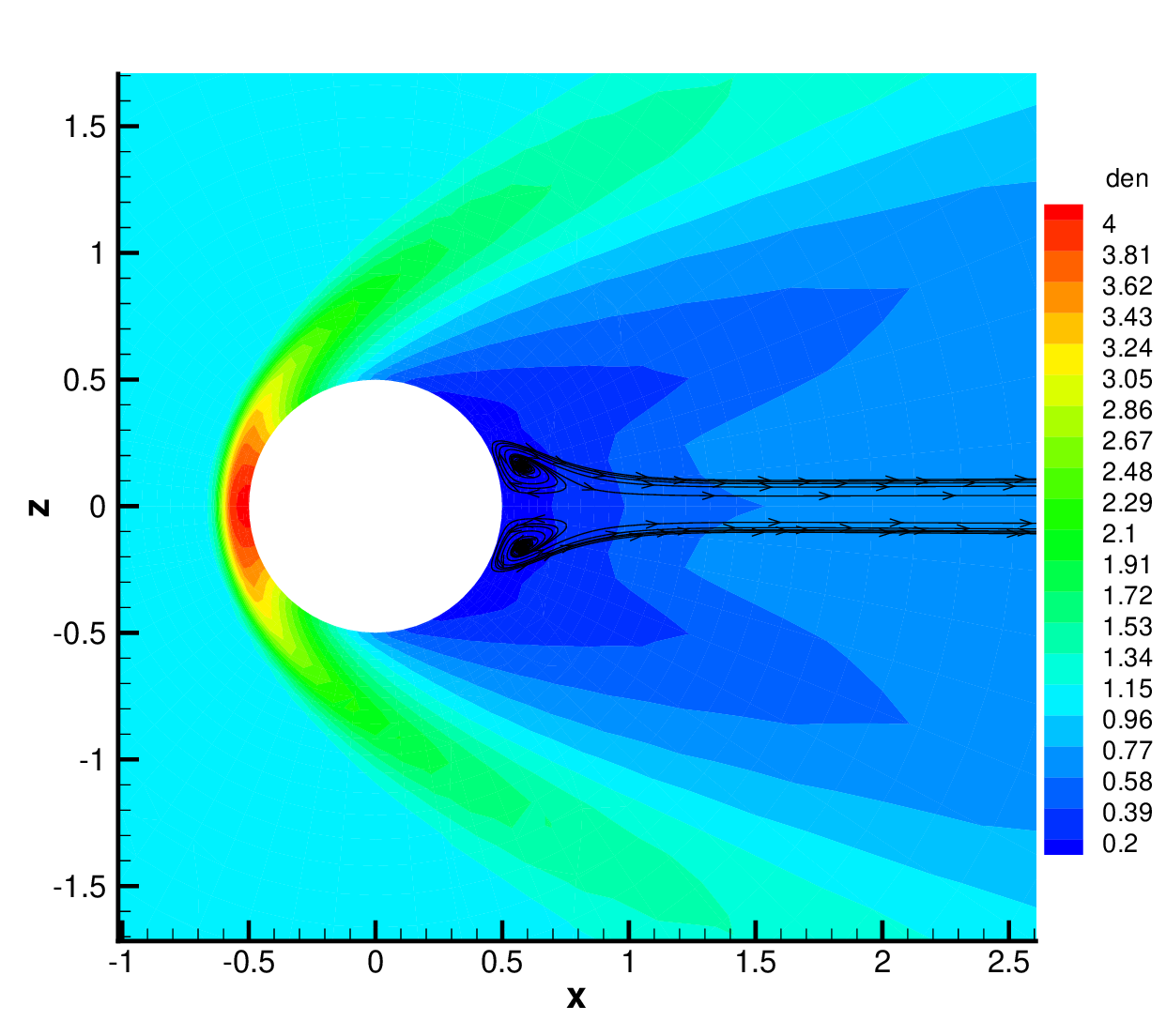}
\includegraphics[width=0.435\textwidth]{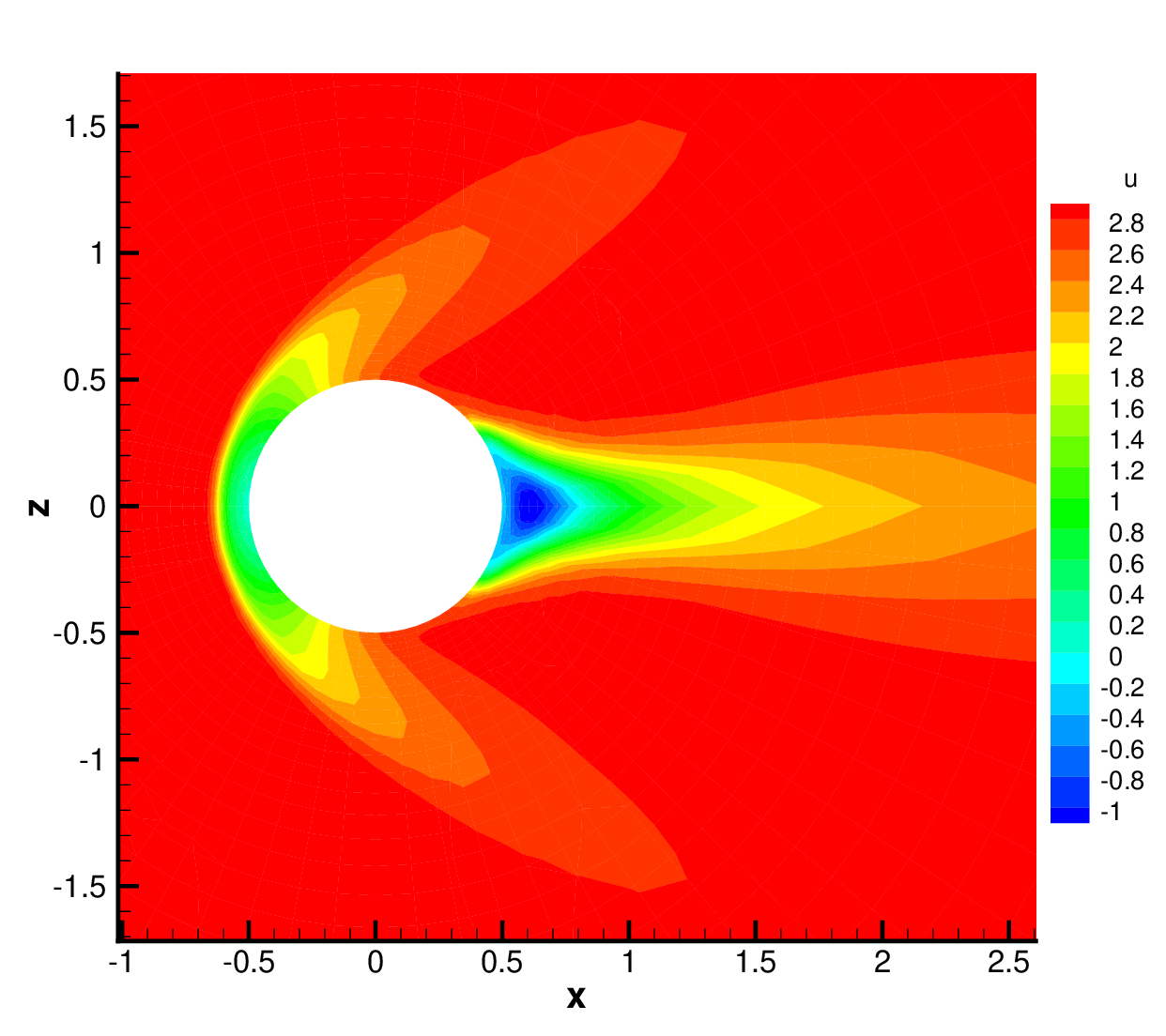}\\
\includegraphics[width=0.435\textwidth]{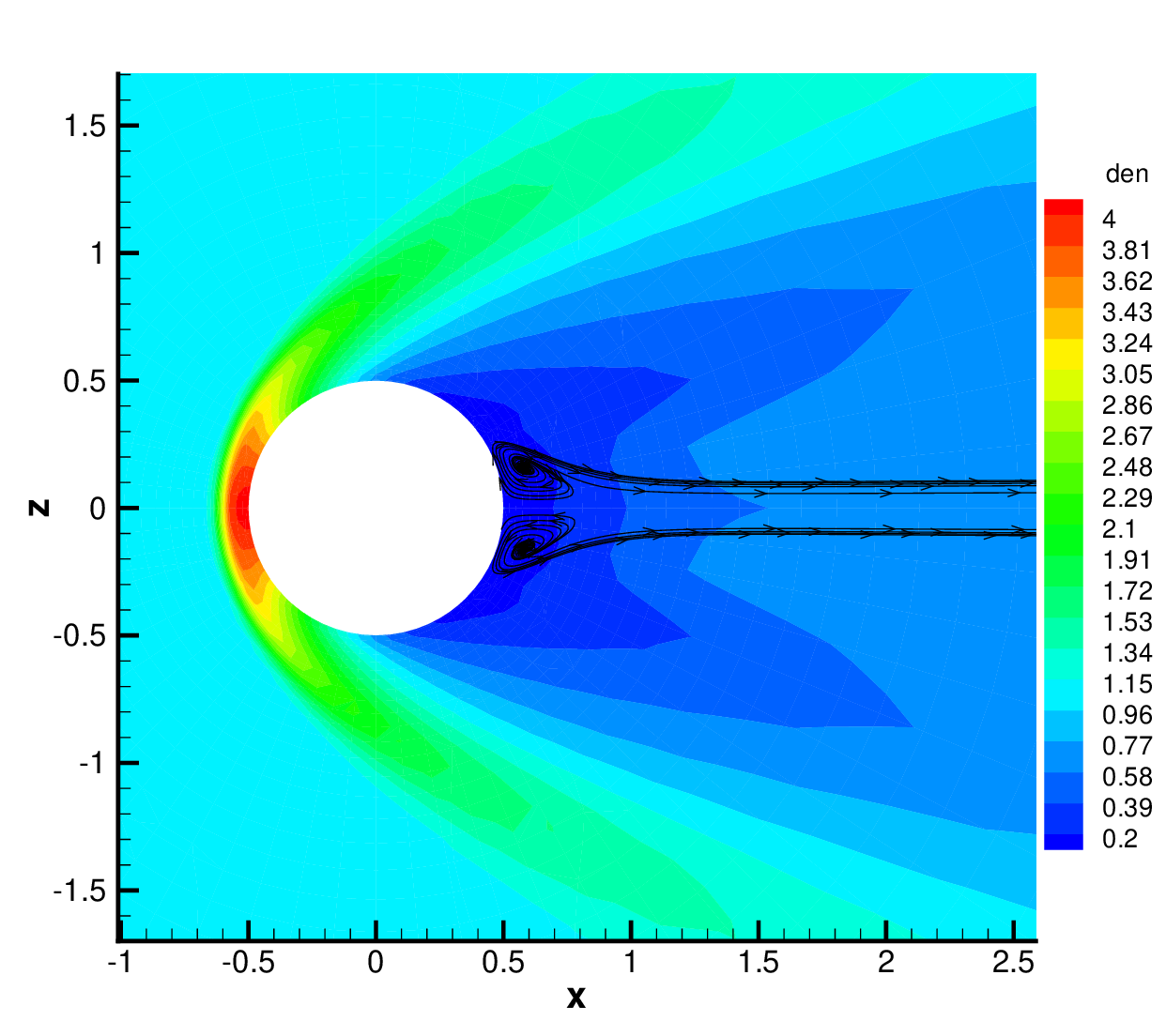}
\includegraphics[width=0.435\textwidth]{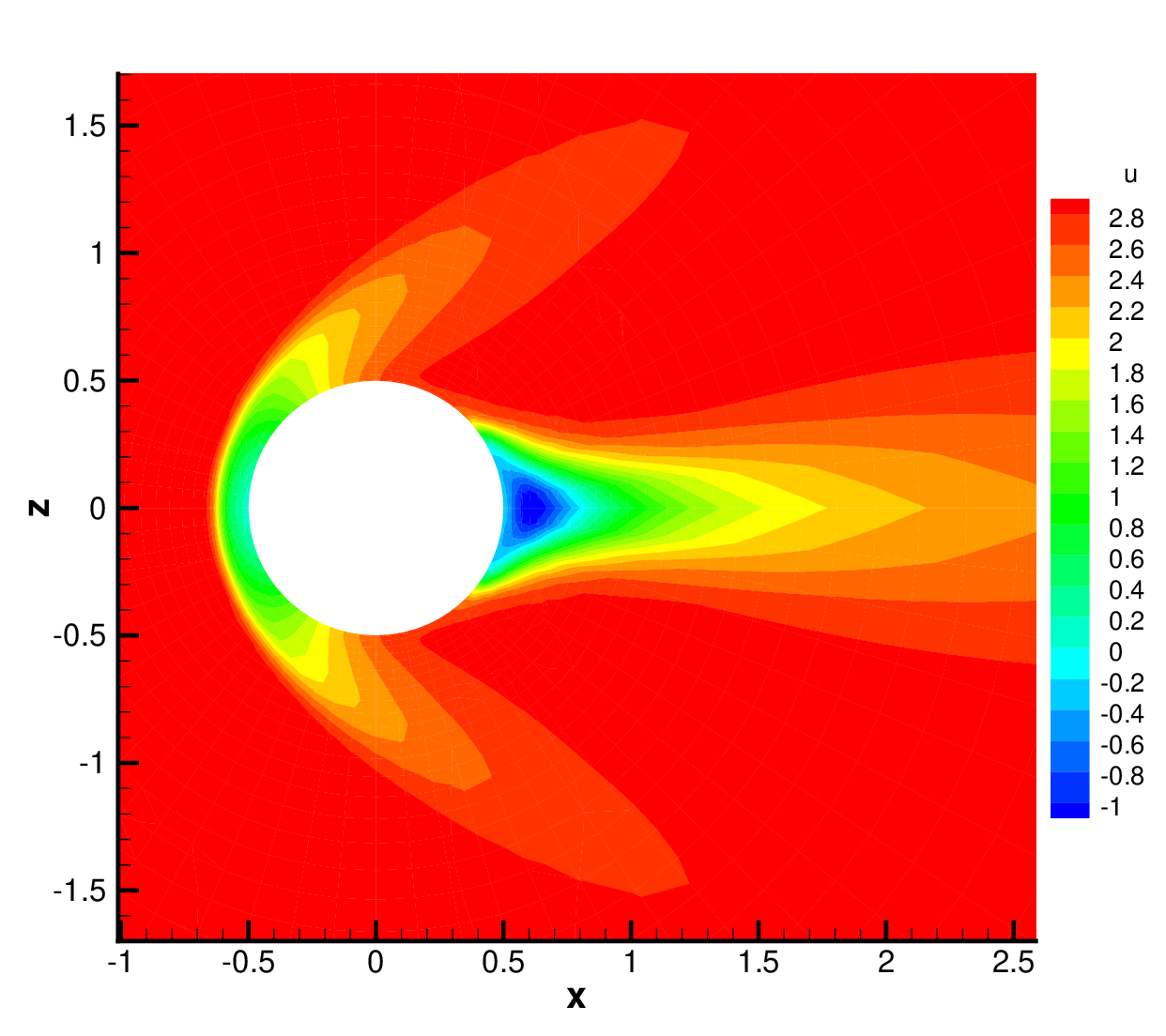}\\
\includegraphics[width=0.435\textwidth]{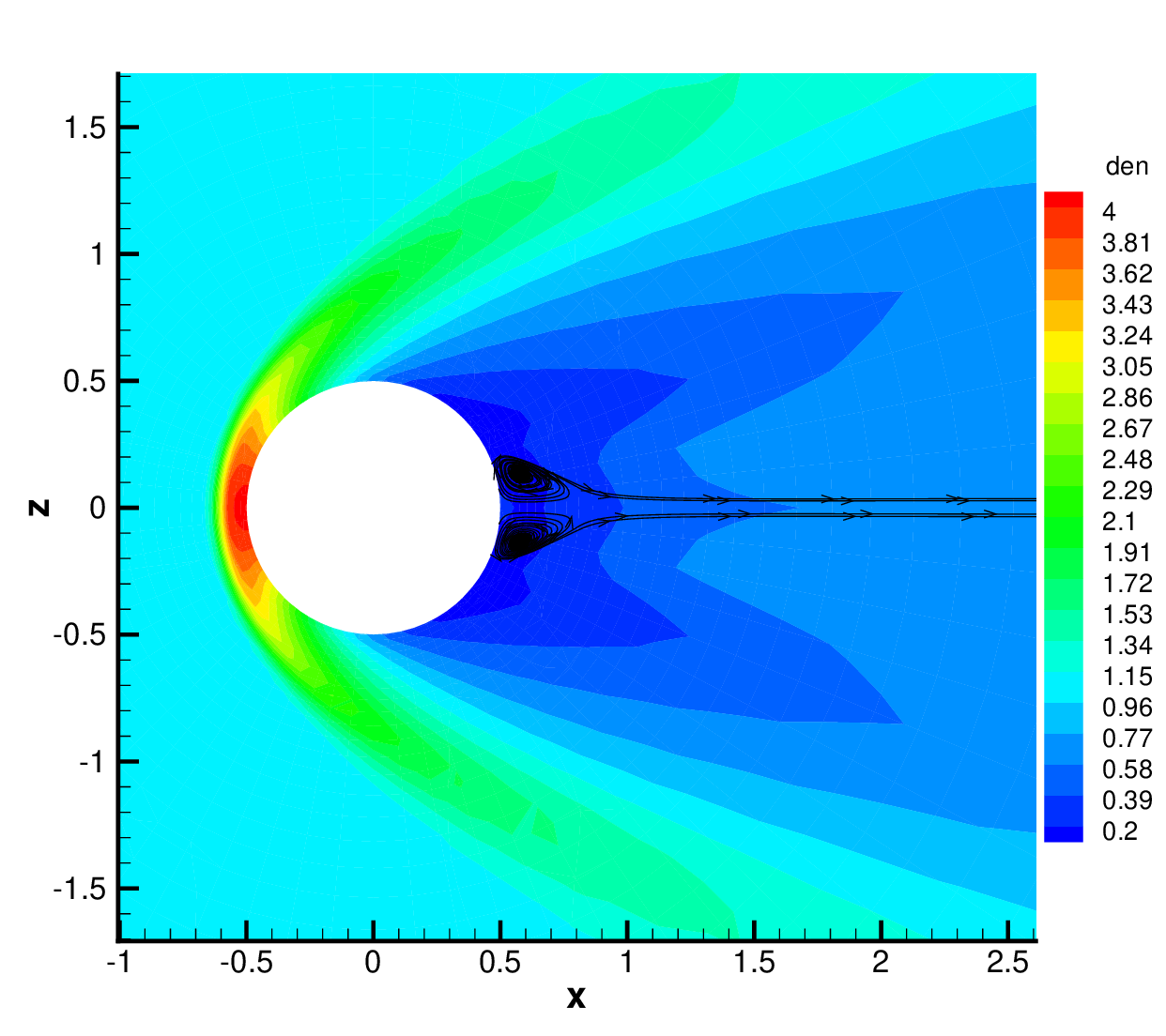}
\includegraphics[width=0.435\textwidth]{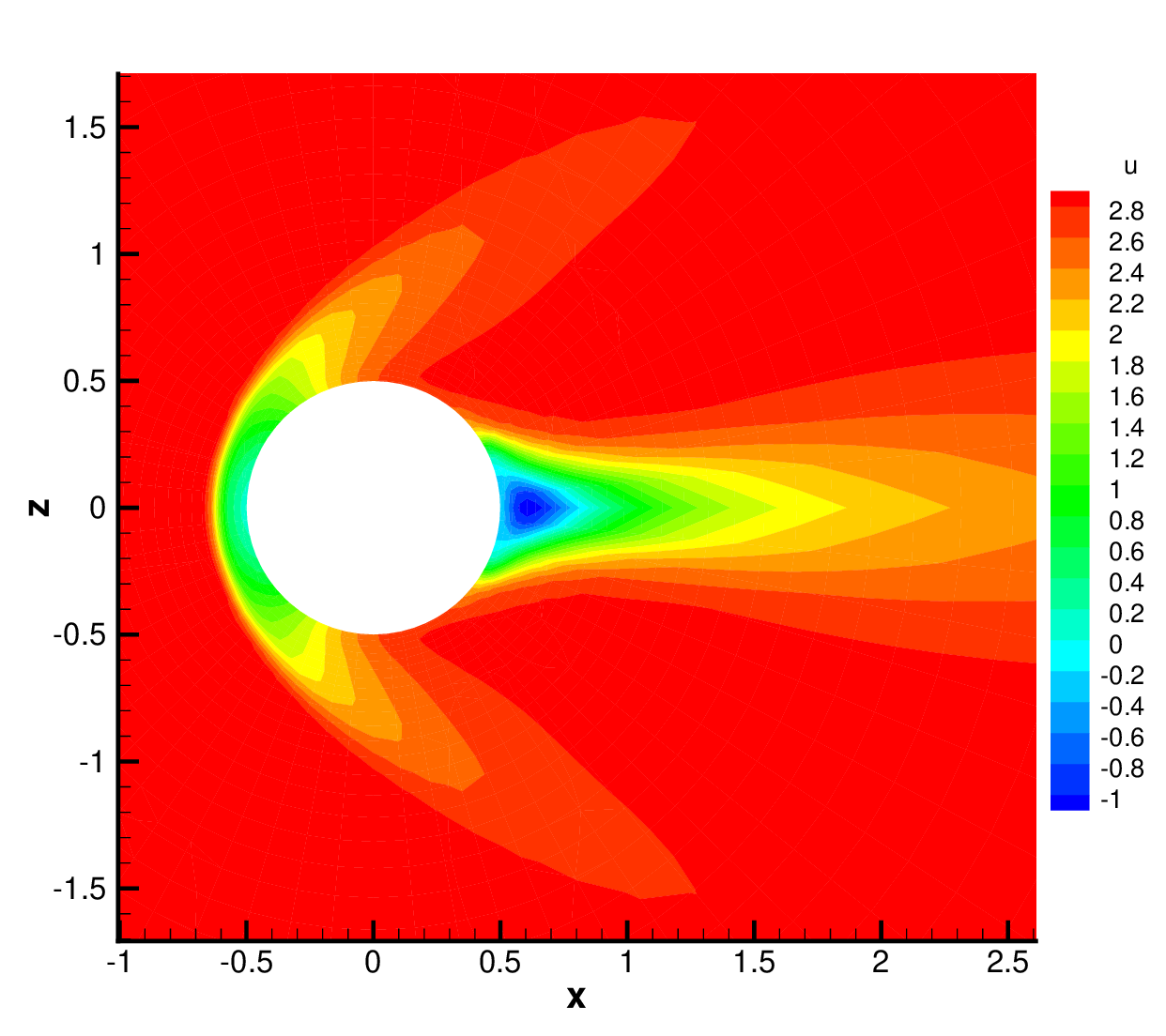}
\caption{\label{sphere-invma3-1} Supersonic inviscid flow passing through a sphere: the density and streamline distributions (left) and
velocity distribution (right) at vertical centerline planes with LUSGS-WENO method (top), GMRES-WENO method (middle) and GMRES-HWENO method (bottom) for $Ma_\infty=3.0$.}
\end{figure}

\begin{figure}[!h]
\centering
\includegraphics[width=0.48\textwidth]{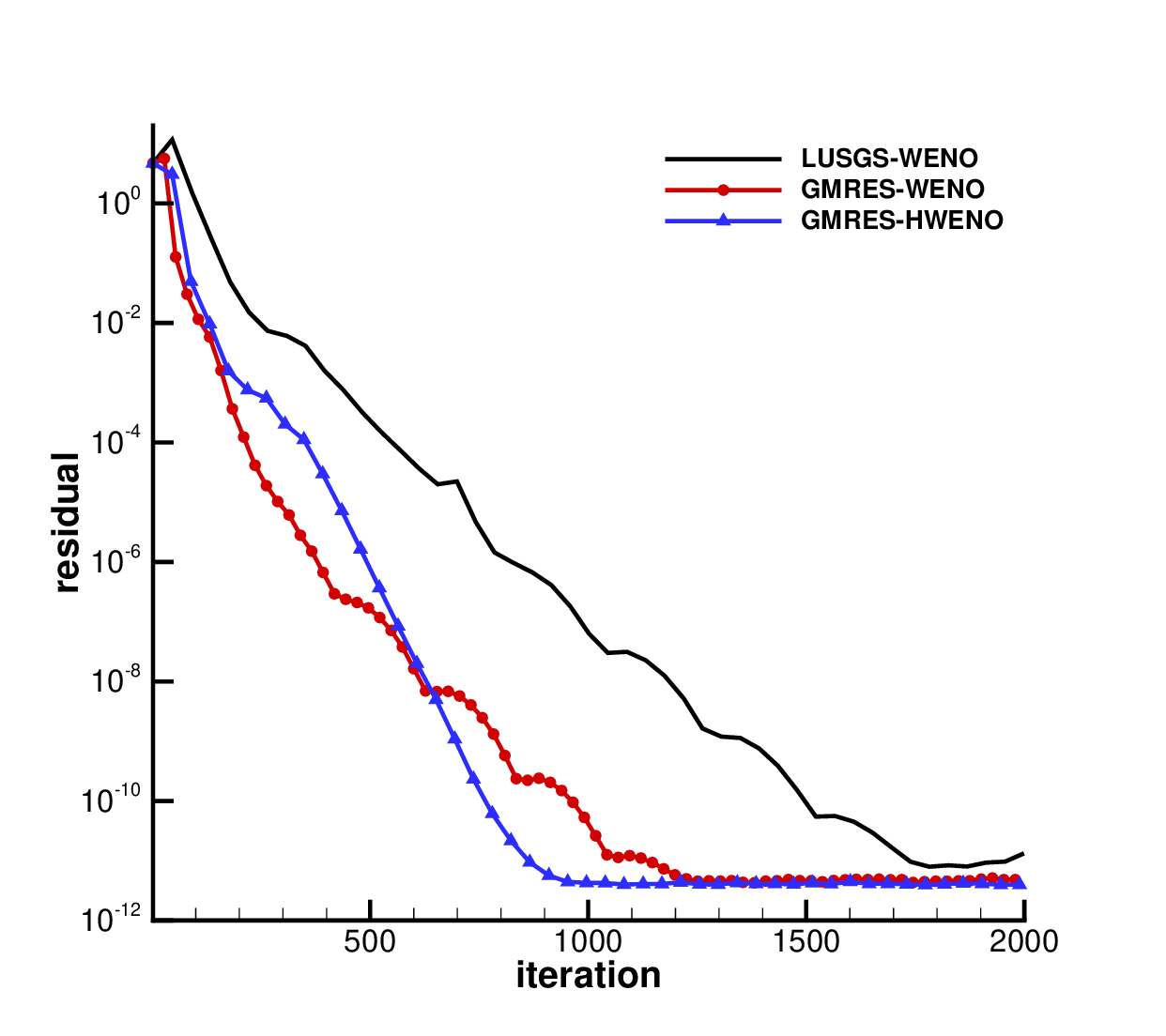}
\caption{\label{sphere-invma3-2} Supersonic inviscid flow passing through a sphere: the residual convergence comparison with LUSGS-WENO, GMRES-WENO and GMRES-HWENO methods for $Ma_\infty=3.0$. }
\end{figure}

For the inviscid flows, the case with $Ma_\infty=3.0$ is performed
to test the robustness of implicit scheme. 
In this case, the CFL number is taken as 3.
The density, velocity and streamline distributions
at vertical centerline planes are shown in
Figure.\ref{sphere-invma3-1}. The shock structure at the leeward side
of the sphere is well captured with current implicit schemes. The
residual convergence histories with different implicit methods are
also shown in Figure.\ref{sphere-invma3-2}. The GMRES method also shows
advantage on the convergence rate and the order of magnitude of
residual.
\begin{figure}[!h]
\centering
\includegraphics[width=0.45\textwidth]{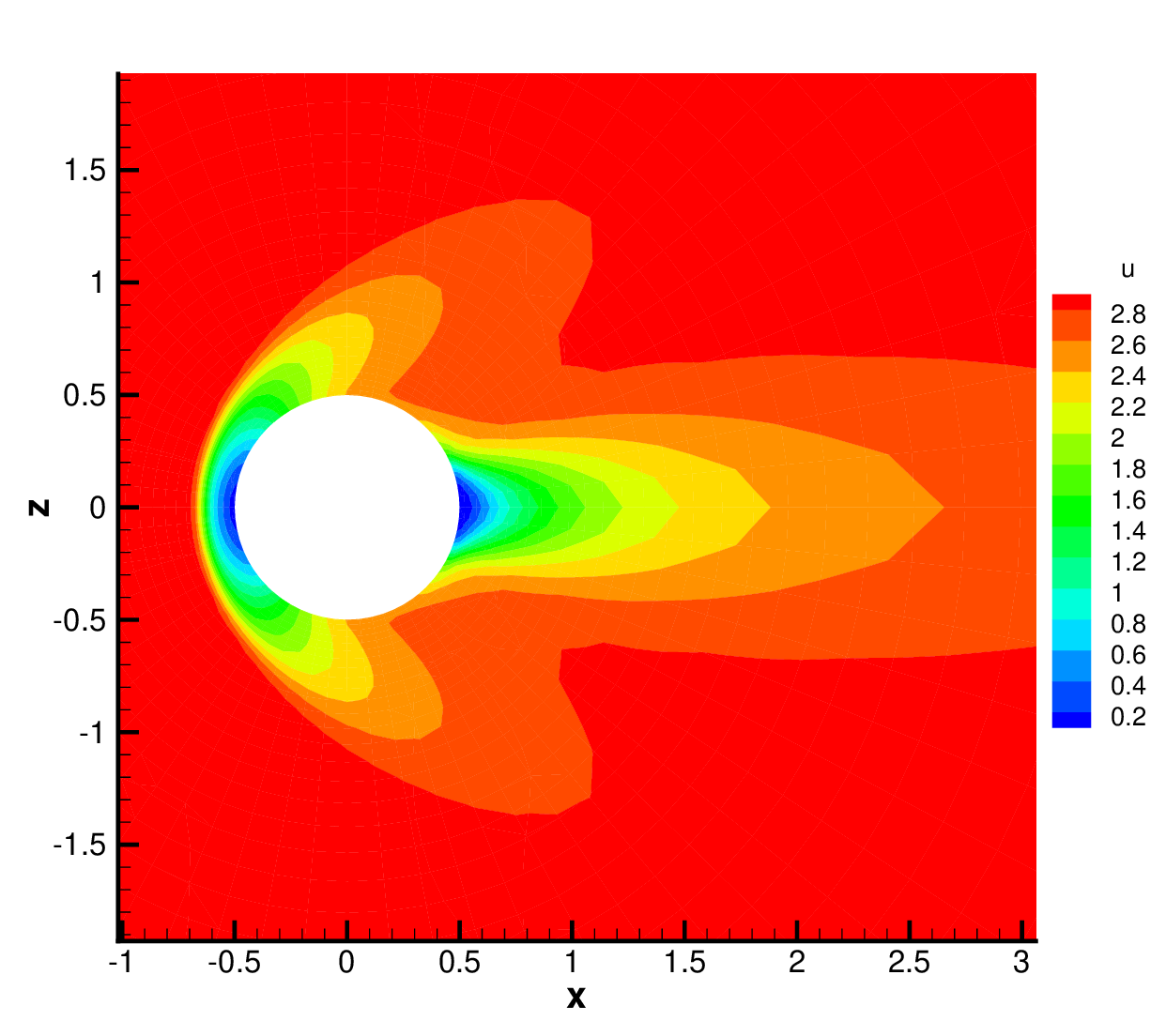}
\includegraphics[width=0.45\textwidth]{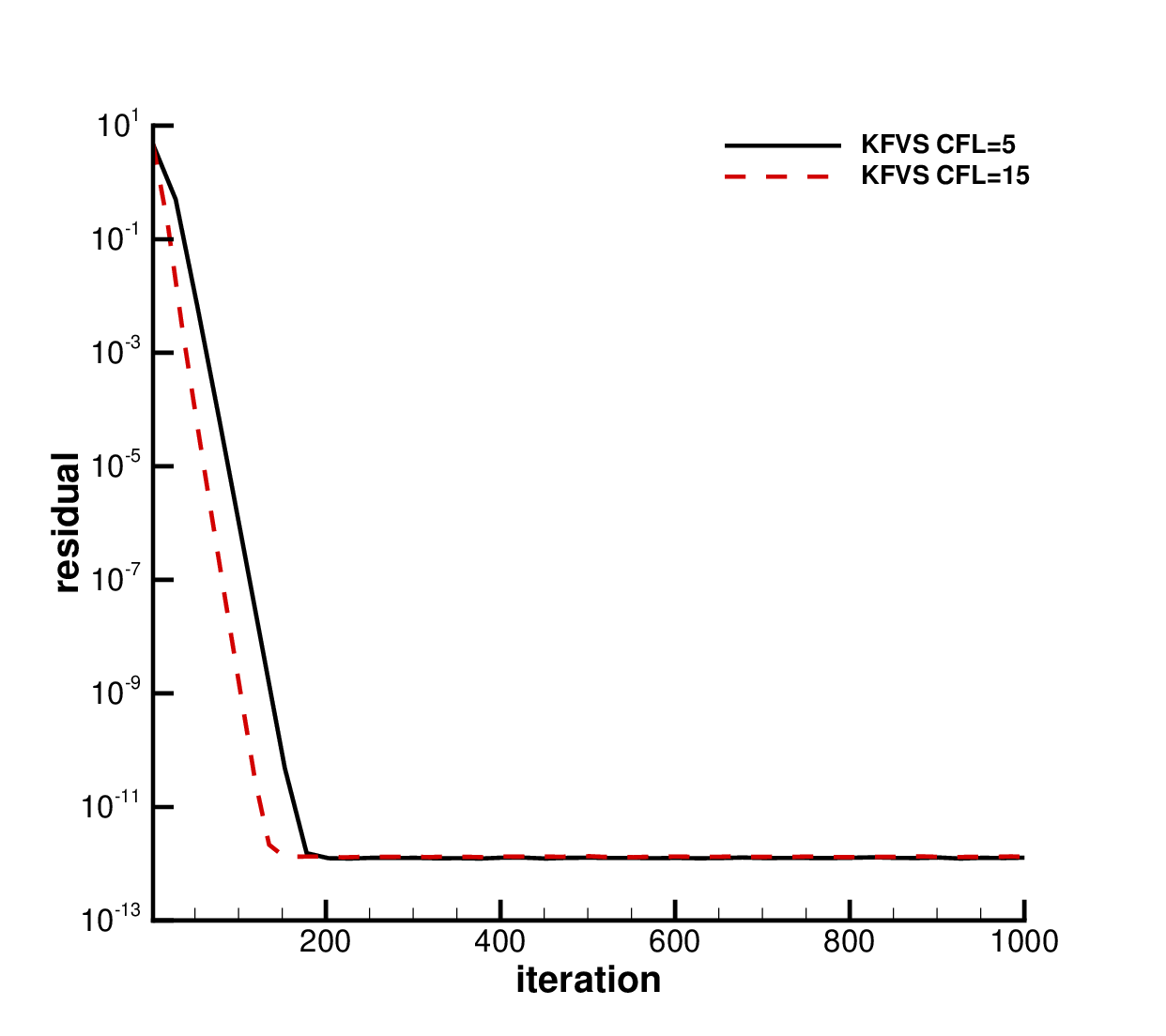}
\caption{\label{sphere-res-invma3-kfvs} Supersonic inviscid flow passing through a sphere: the velocity distribution at
vertical centerline planes with $CFL=5$ (left) and the residual convergence comparison (right) tested by
the GMRES method combined with first order KFVS reconstruction for $Ma_\infty=3.0$.}
\end{figure}
In this case, the minimum number of CFL number is not sufficiently
large. In order to explore the relationship between the
magnification of CFL number and the characteristic of high-order
scheme, the above supersonic inviscid case with $Ma_\infty=3.0$ is
retested by the implicit GMRES method with first order KFVS flux
solver. The velocity distribution at vertical centerline planes and
residual convergence histories with different CFL numbers are shown
in Figure.\ref{sphere-res-invma3-kfvs}. It shows that the magnification
of the CFL number for lower-order schemes on good quality meshes
will be more ideal. In the follow-up optimization work, a more
reasonable combination between the implicit method and high-order
scheme will be investigated.

\begin{figure}[!h]
\centering
\includegraphics[width=0.45\textwidth]{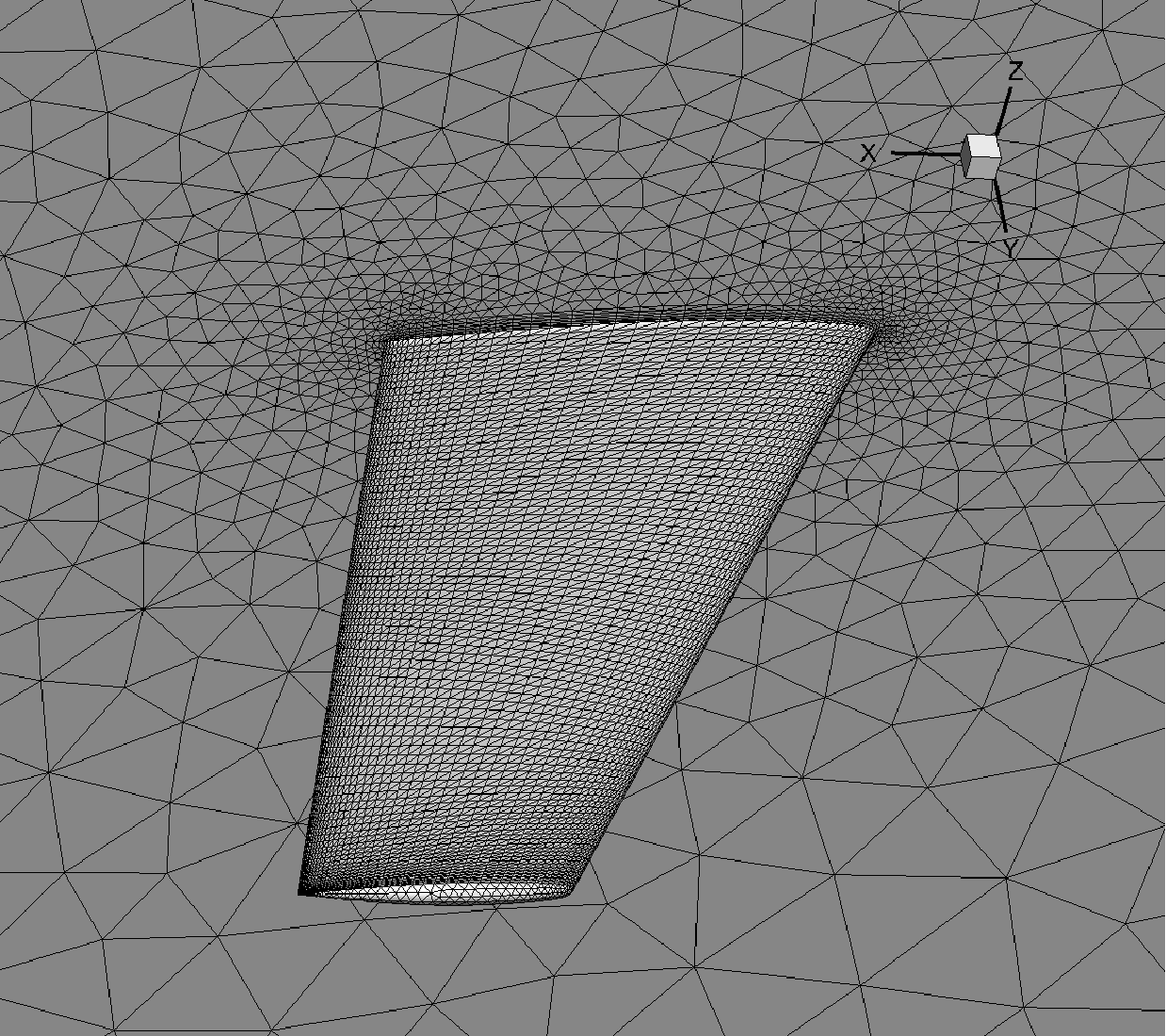}
\caption{\label{m6wing-mesh} Transonic flow around ONERA M6 wing: the local computational mesh distribution with tetrahedral meshes.}
\end{figure}

\begin{figure}[!h]
\centering
\includegraphics[width=0.42\textwidth]{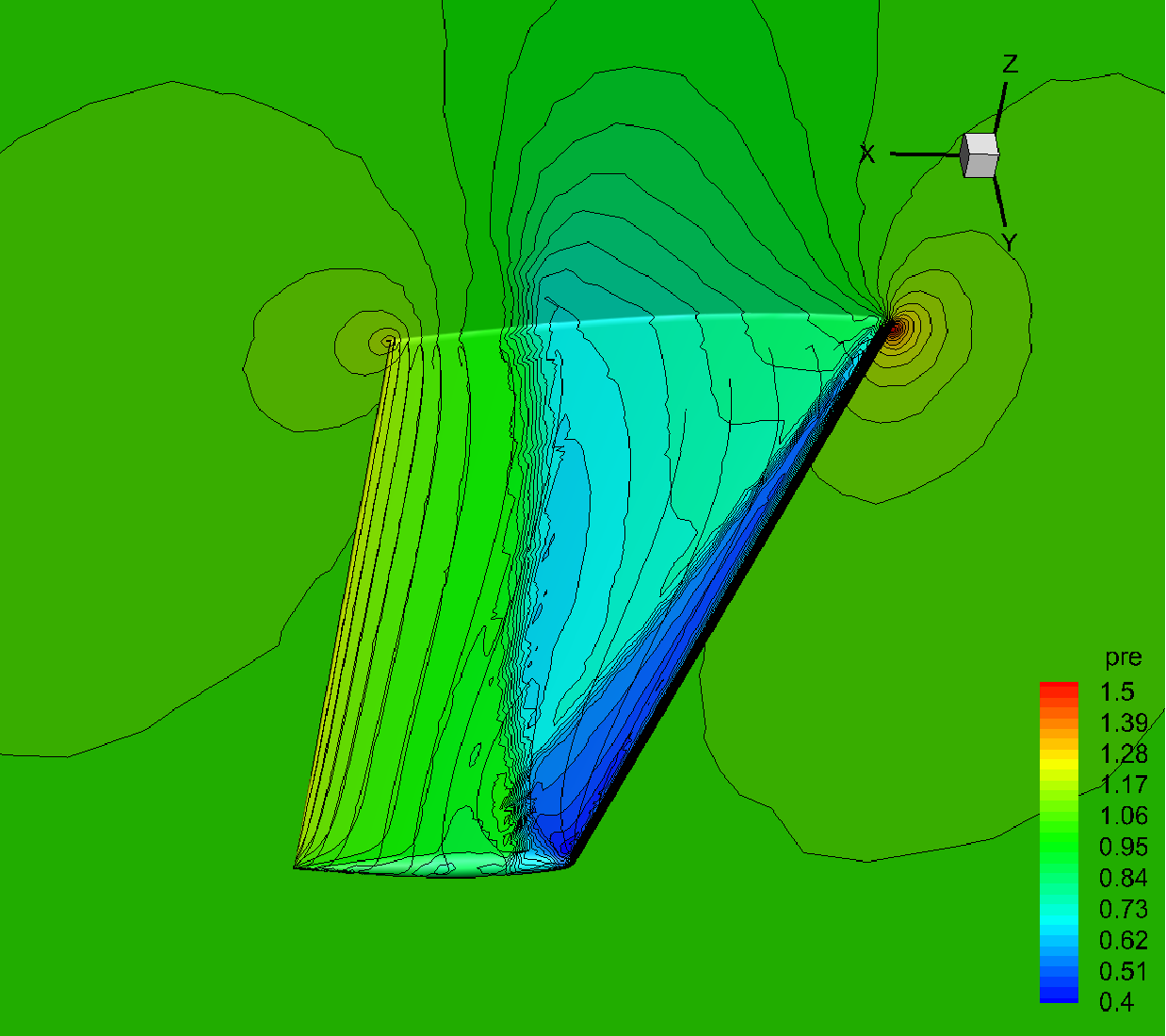}
\includegraphics[width=0.42\textwidth]{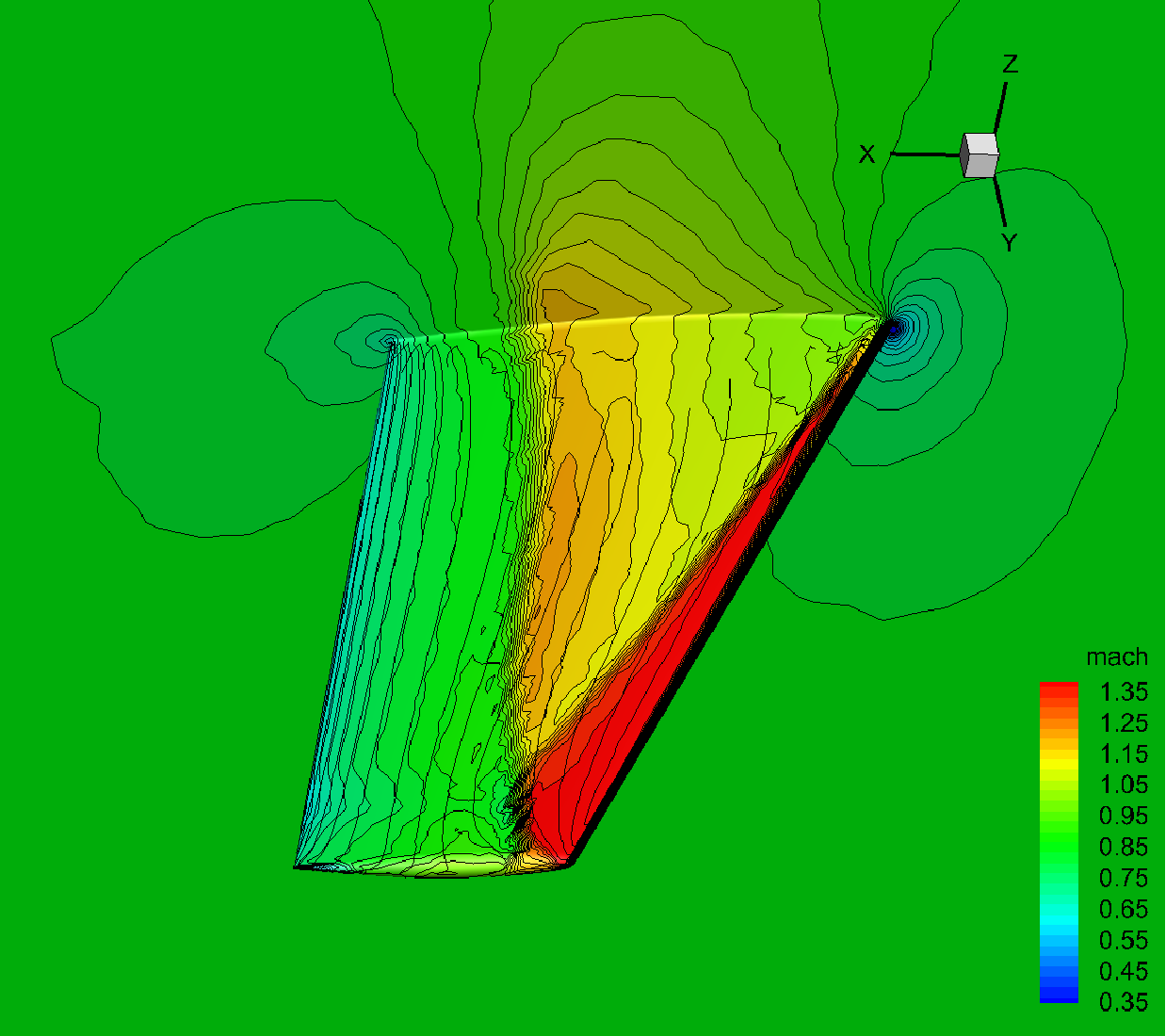}\\
\includegraphics[width=0.42\textwidth]{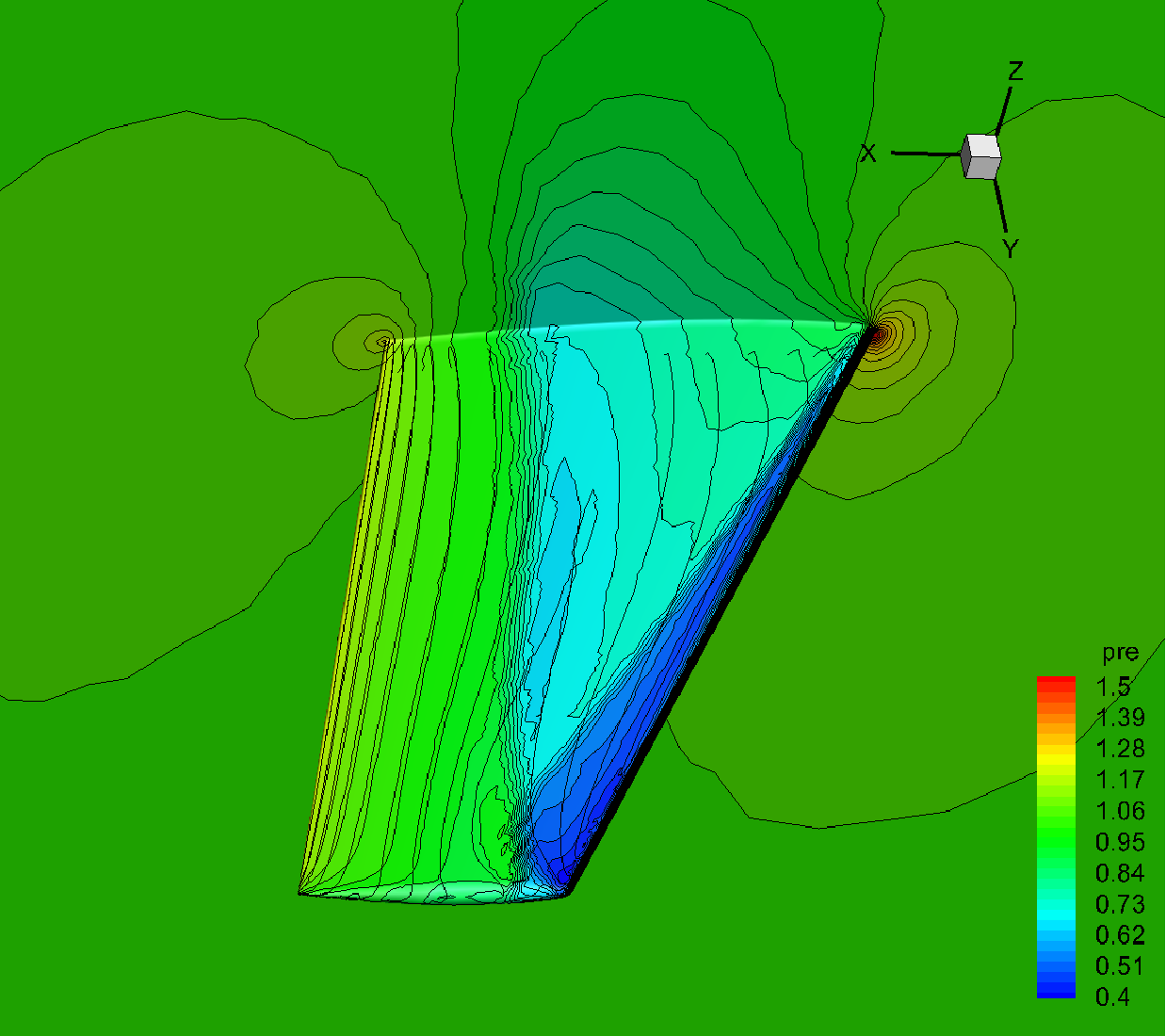}
\includegraphics[width=0.42\textwidth]{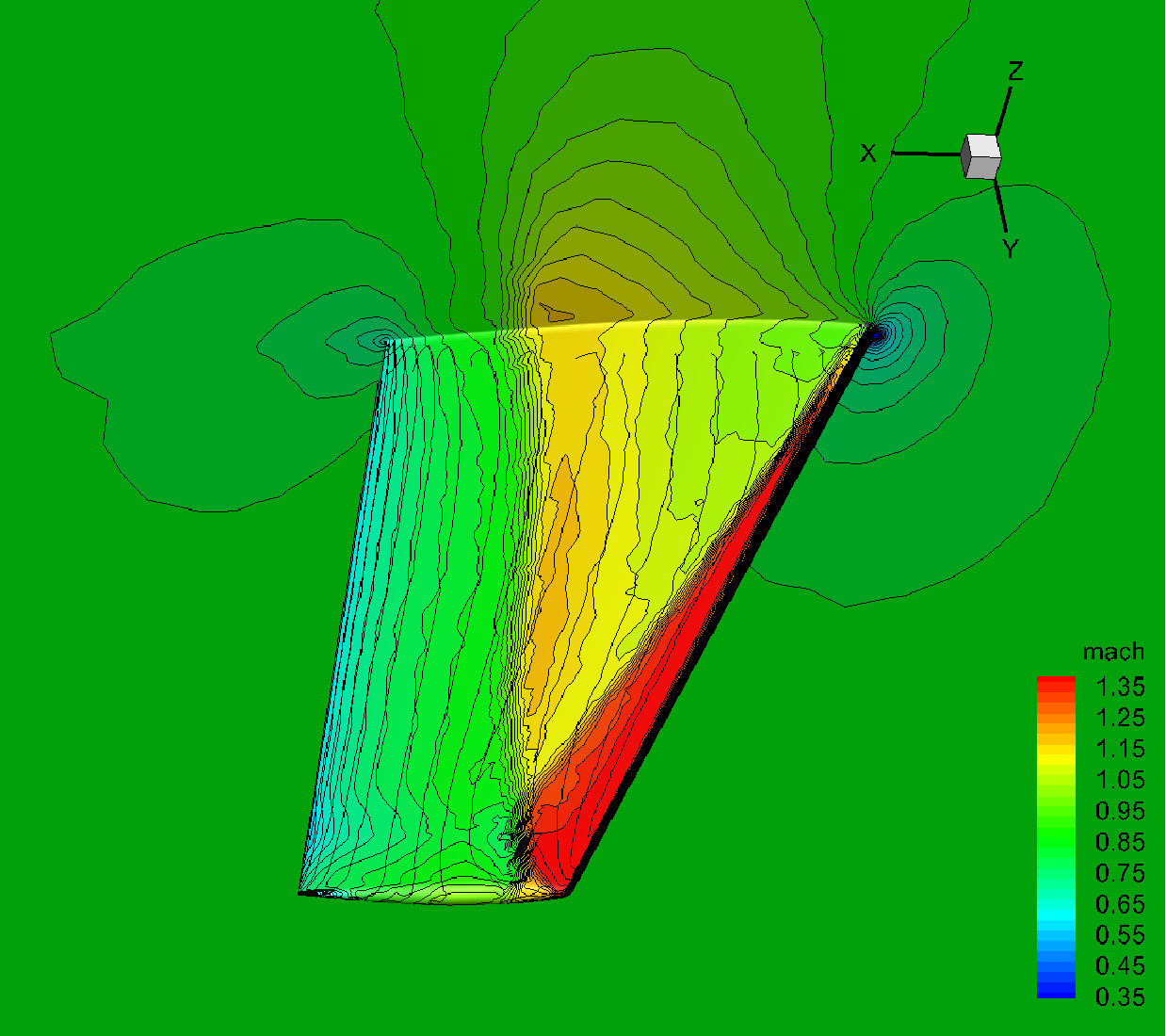}\\
\includegraphics[width=0.42\textwidth]{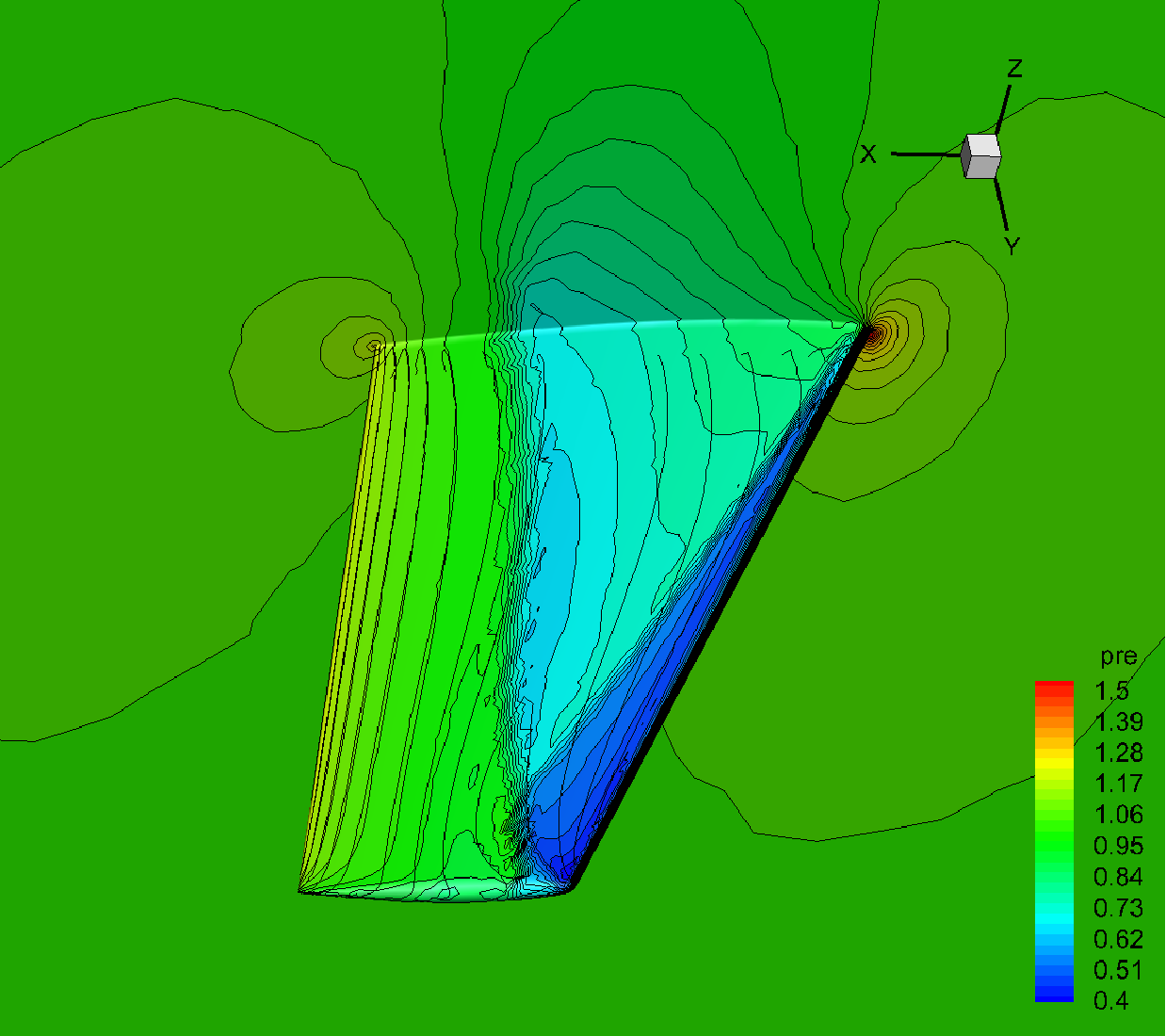}
\includegraphics[width=0.42\textwidth]{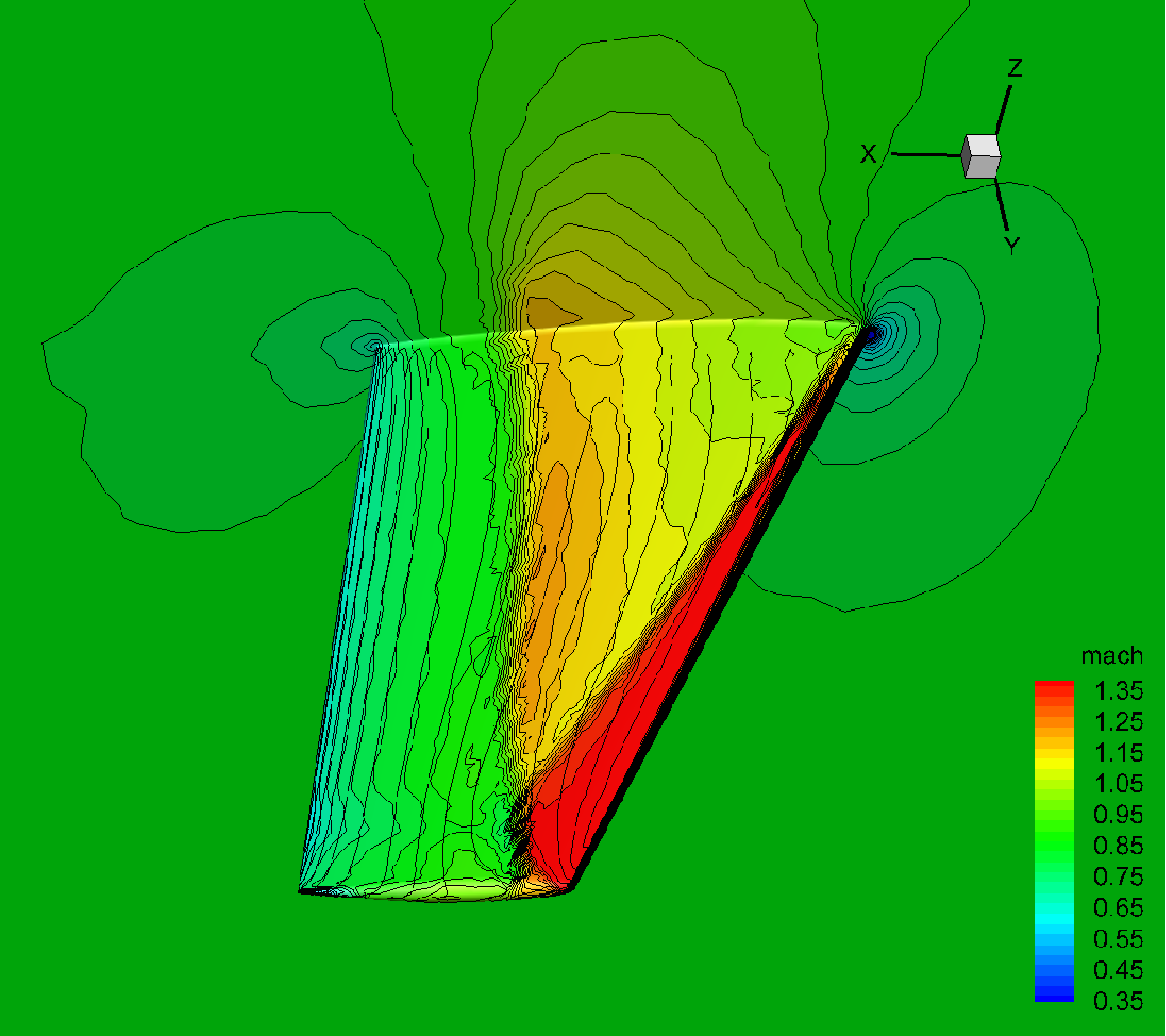}
\caption{\label{M6-wing-1} Transonic flow around ONERA M6 wing:  the local pressure distribution with LUSGS-WENO method (top), GMRES-WENO method (middle) and GMRES-HWENO method (bottom).}
\end{figure}

\begin{figure}[!h]
\centering
\includegraphics[width=0.42\textwidth]{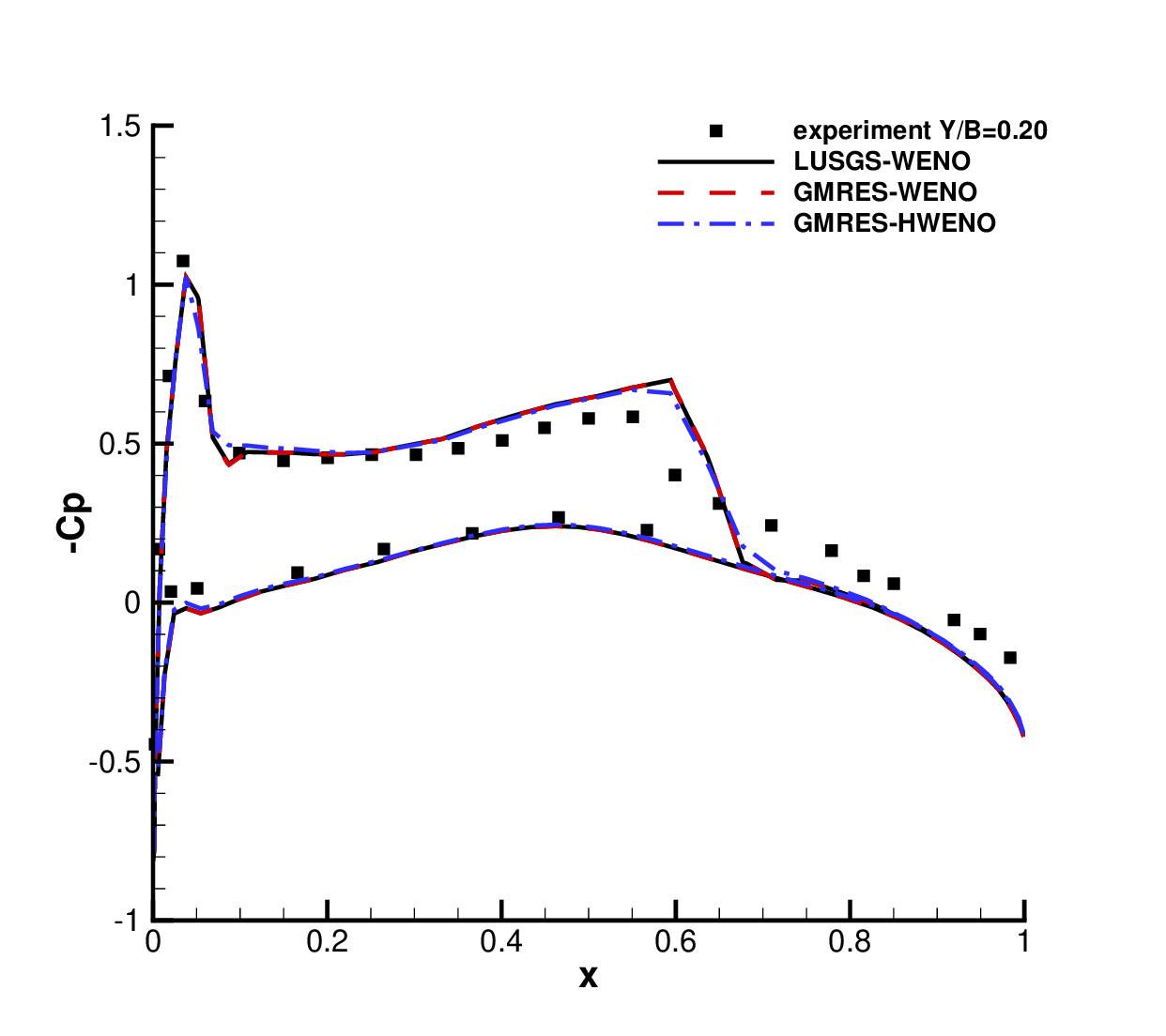}
\includegraphics[width=0.42\textwidth]{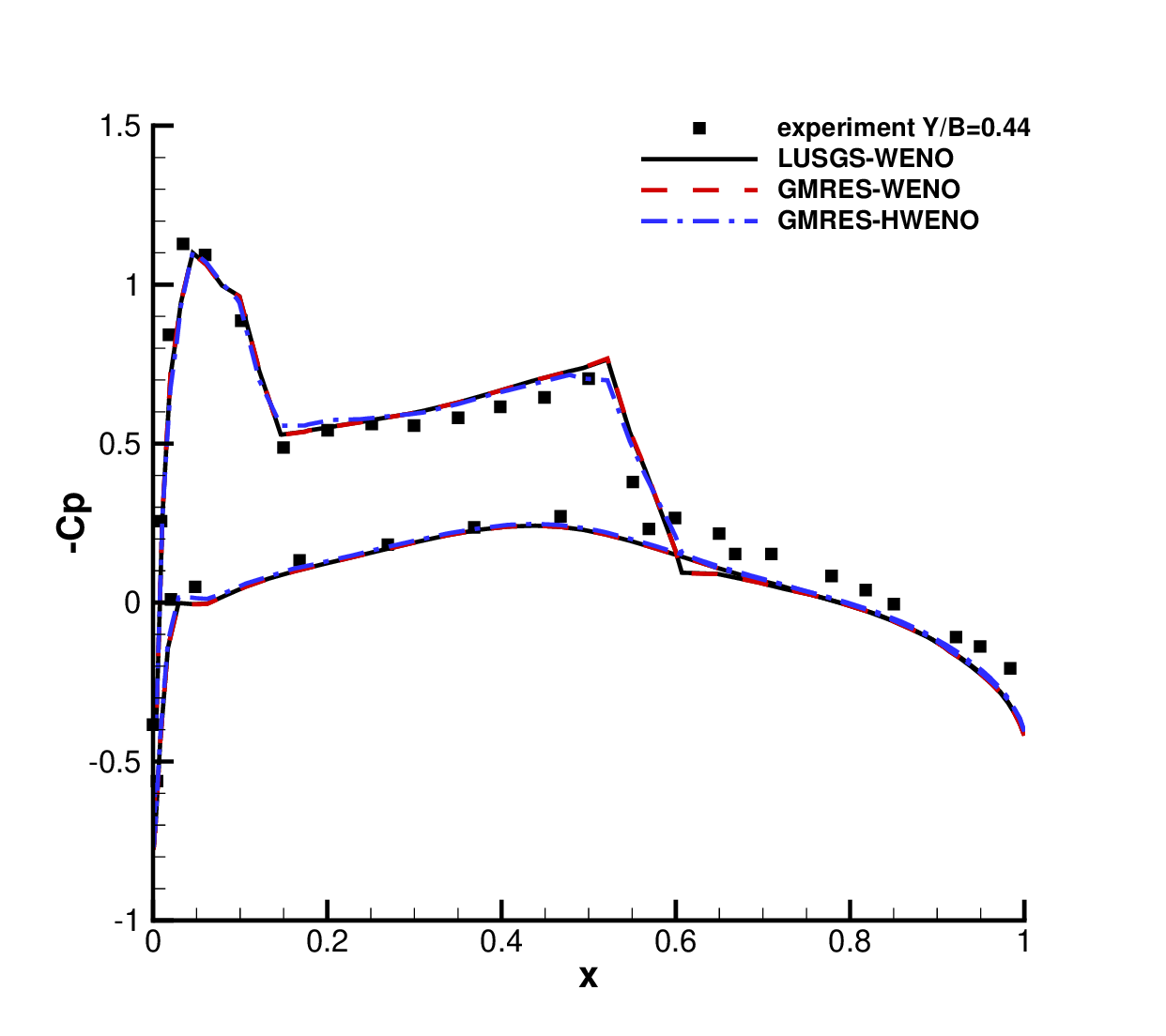}\\
\includegraphics[width=0.42\textwidth]{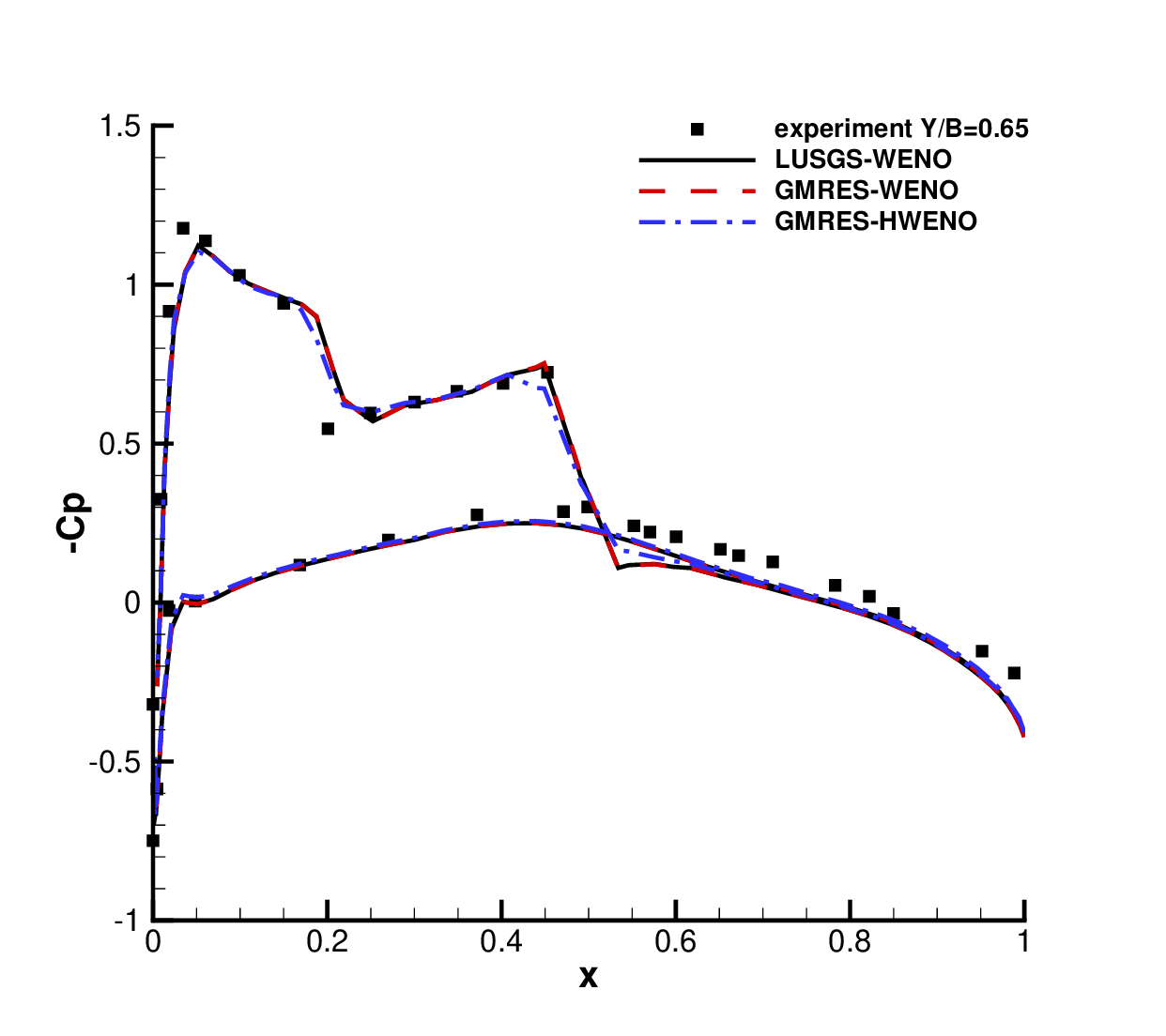}
\includegraphics[width=0.42\textwidth]{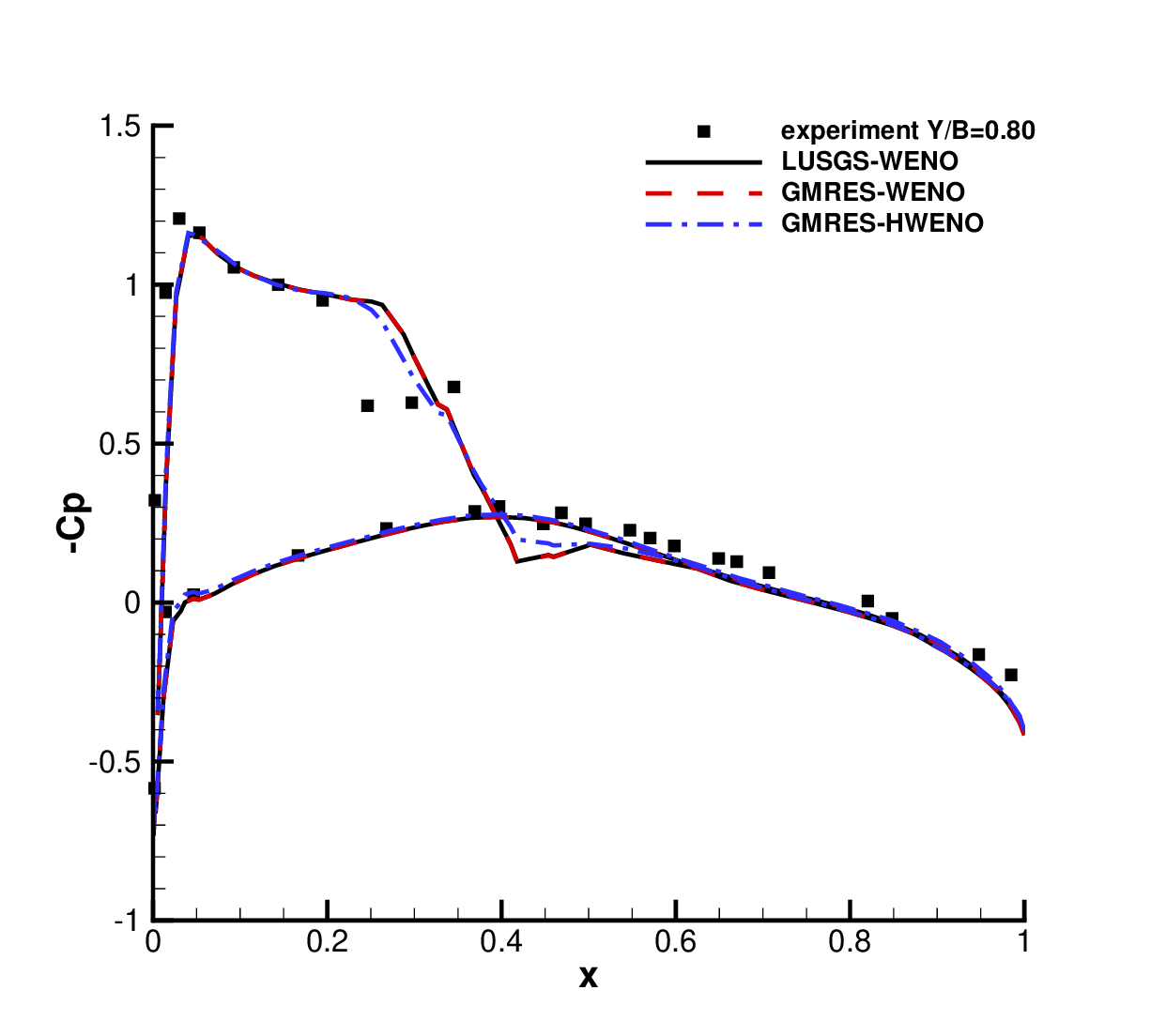}\\
\includegraphics[width=0.42\textwidth]{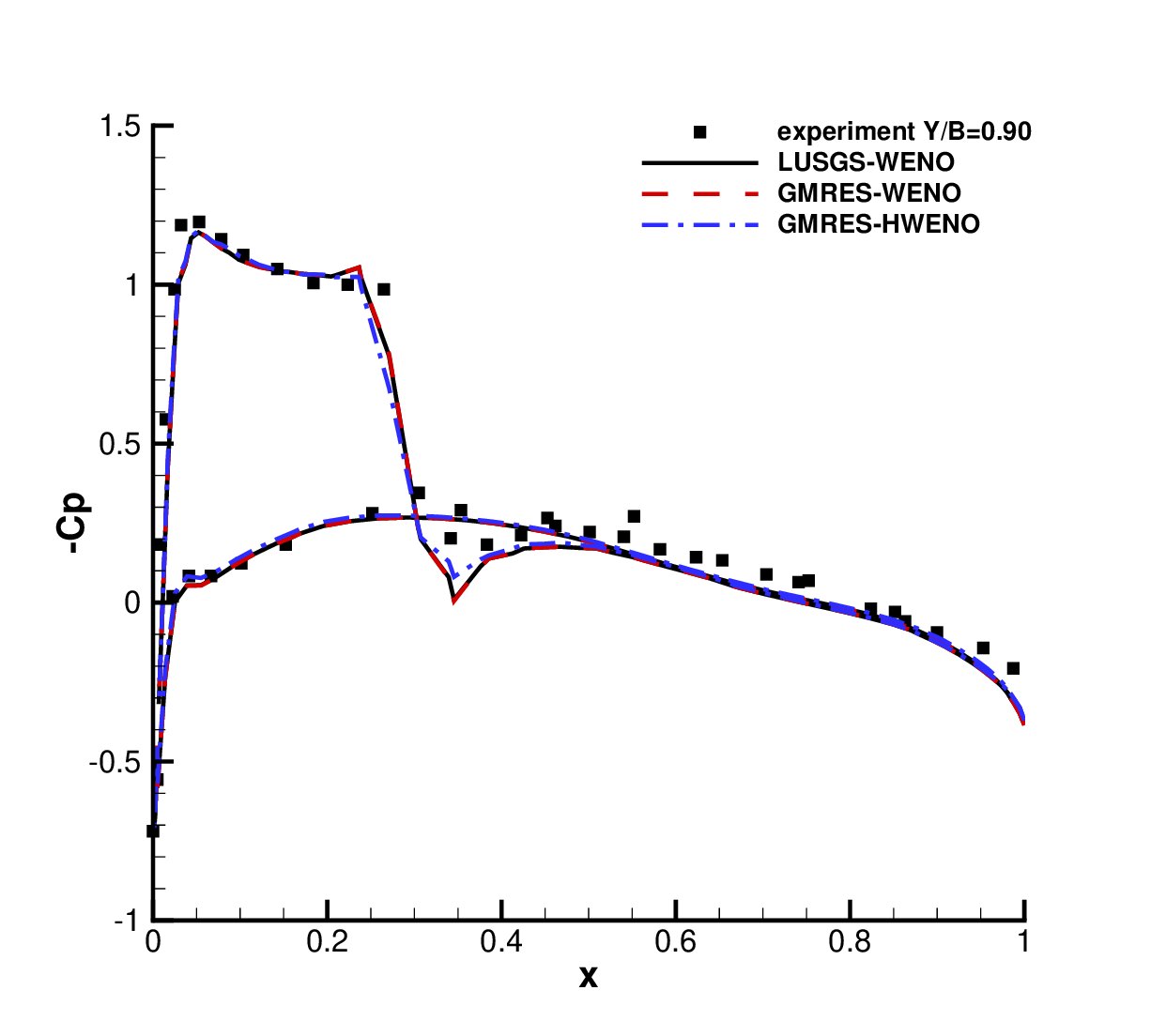}
\includegraphics[width=0.42\textwidth]{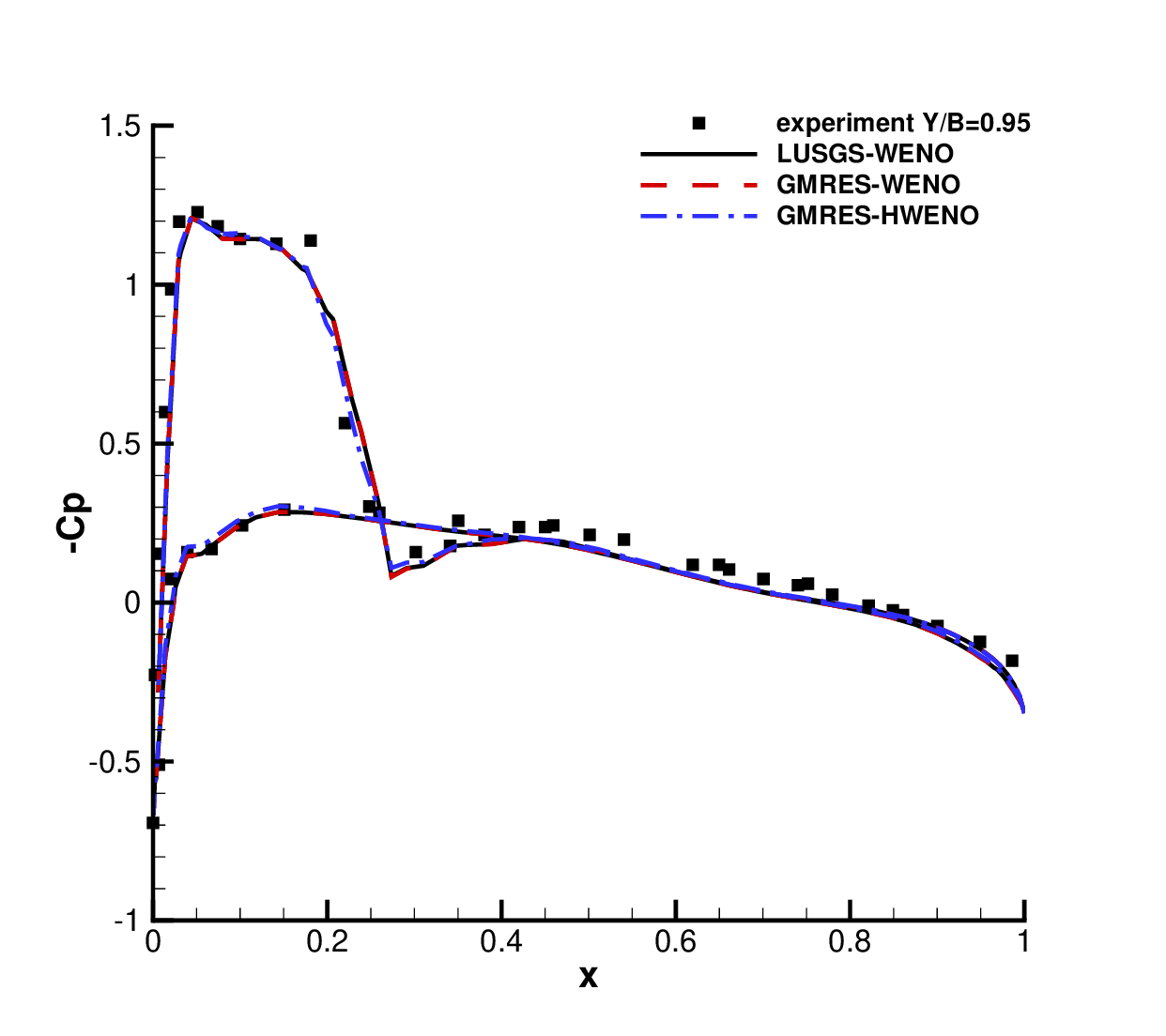}
\caption{\label{M6-wing-2} ONERA M6 wing: the pressure coefficient
distributions at $Y/B = 0.20$, $0.44$, $0.65$, $0.80$, $0.90$ and
$0.95$ for the inviscid flow with tetrahedral and hybrid meshes.}
\end{figure}

\begin{figure}[!h]
\centering
\includegraphics[width=0.48\textwidth]{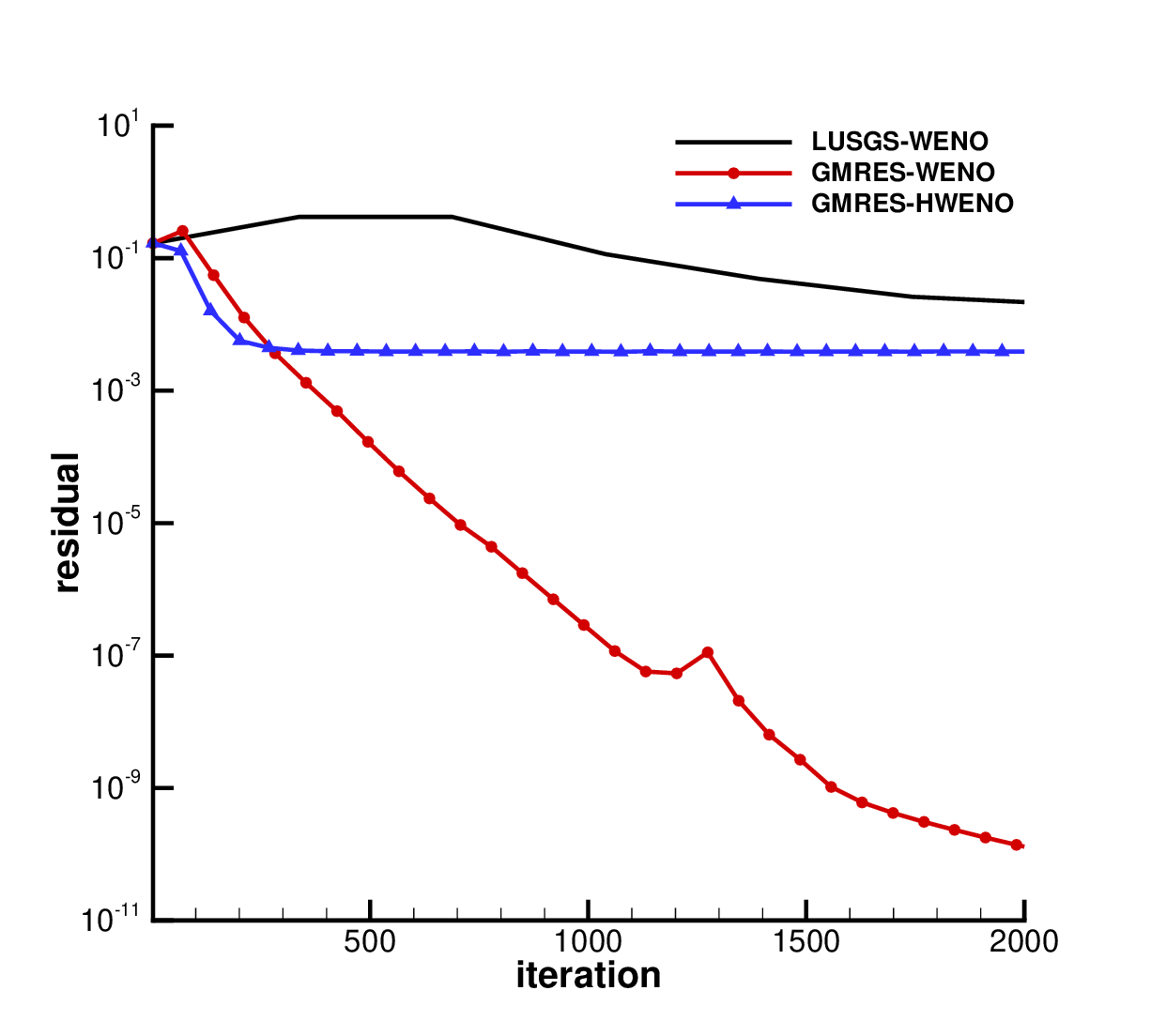}
\caption{\label{m6wing-res} Transonic flow around ONERA M6 wing: the residual convergence comparison with LUSGS-WENO, GMRES-WENO and GMRES-HWENO methods.}
\end{figure}

\subsection{Transonic flow around ONERA M6 wing}
The transonic flow around the ONERA M6 wing is a standard benchmark
for engineering simulations. Besides the three-dimensional geometry,
the flow structures are complex including the interaction of shock
and turbulent boundary. Thus, it is a good candidate to test the
performance of the extended BGK model and implicit high-order
gas-kinetic scheme. The inviscid flow around the wing is tested,
which corresponds to a rough prediction of the flow field under a
very high Reynolds number. The incoming Mach number and angle of
attack are given by
\begin{align*}
Ma_{\infty} = 0.8395, ~AoA = 3.06^{\circ}.
\end{align*}
This case is performed by the tetrahedral meshes, which includes
294216 cells. The mesh distribution is shown in
Figure.\ref{m6wing-mesh}. The stencils $S_{i_m}^{HWENO}$ are selected as tetrahedral compact sub-stencils, which are dominated by the gradient informations. 
The subsonic inflow and outflow boundaries
are all set according to the local Riemann invariants, and the
adiabatic and slip wall condition is imposed on the solid wall. The
local pressure distributions with the LUSGS method, the non-compact
GMRES method and compact GMRES method are shown in
Fig.\ref{M6-wing-1}, and the $\lambda$ shock is well resolved by all
the current implicit schemes. The comparisons on the pressure
distributions at the semi-span locations $Y/B = 0.20$, $0.44$,
$0.65$, $0.80$, $0.90$ and $0.95$ of the wing are given in
Figure.\ref{M6-wing-2}. The numerical results quantitatively agree well
with the experimental data \cite{Case-Schmitt}. The histories of residual
convergence with different implicit methods are given in
Figure.\ref{m6wing-res}. It shows that the non-compact GMERS method
converges faster than the LUSGS method. Although the residual
convergence is affected by the dominant gradient information in
tetrahedral compact sub-stencils, the compact GMRES still shows the
advantage over the LUSGS method.  
In the future, the selection of stencils and the evolution of the cell averaged gradients in time will be further investigated.

\subsection{Efficiency comparison of CPU and GPU}
The efficiency comparison of CPU and GPU codes is provided for both implicit  
non-compact and compact schemes. The CPU code is run with Intel Xeon Gold 6230R CPU
using Intel Fortran compiler with 16 OpenMP threads, while Nvidia
Quadro RTX 8000 is used for GPU code with Nvidia CUDA and NVFORTRAN
compiler. The clock rates of GPU and CPU are 1.77 GHz and 2.10 GHz
respectively, and the double precision is used in computation. The
lid-driven cavity flow with $Re=1000$ on three dimensional
unstructured tetrahedral mesh are used to test the efficiency. The
HGKS with GMRES method combined by WENO and HWENO reconstruction are
implemented with both GPU and CPU. As shown in Table.\ref{GPU-CPU},
8x speedup is achieved for non-compact and compact GMRES GPU code compared with
non-compact and compact GMRES CPU code respectively.
Taking the number of thread of CPU into account, the speedup of GPU code approximately equals to 130.
However, the compact scheme needs to store more local information than the
non-compact scheme, the computational scale is constrained by the
limited available memory of single GPU. In the future, the implicit
compact HGKS on unstructured meshes will be further upgraded with
multiple GPUs using MPI and CUDA, and more challenging problems for
compressible flows will be investigated.

\begin{table}[!h]
\begin{center}
\def\temptablewidth{1.0\textwidth}{\rule{\temptablewidth}{0.8pt}}
\begin{tabular*}{\temptablewidth}{@{\extracolsep{\fill}}c|c|c|c|c}
Scheme & Computational Mesh & CPU+OpenMP & GPU+CUDA& Speedup \\
\hline
noncompact GMRES  & 48000 Tet &  259s   & 30s &  8.63 \\
compact GMRES  & 48000 Tet &  267s  &  32s &  8.34
\end{tabular*}
{\rule{\temptablewidth}{1.0pt}}
\end{center}
\caption{\label{GPU-CPU} Efficiency comparison: the computational time per 100 steps and speedup for GPU and CPU code.
Taking the number of thread of CPU into account, the speedup of GPU code approximately equals to 130.}
\end{table}

\section{Conclusion}
In this paper, the implicit non-compact and compact HGKSs are developed on the
three-dimensional unstructured meshes. For non-compact GKS scheme, the third-order WENO reconstruction is used, 
where the stencils are selected from the neighboring cells and the neighboring cells of neighboring cells.  
Incorporate with the GMRES method based on numerical Jacobian matrix, the implicit non-compact HGKS is developed for steady problems. 
To improve the resolution and parallelism, the implicit compact HGKS is also developed with HWENO reconstruction, 
where the stencils only contain one level of neighboring cells. 
The cell averaged conservative variable is updated with GMRES method. Simultaneously, a simple strategy is used to update the cell averaged 
gradient with the spatial-temporal coupled gas-kinetic flow solver. 
To further accelerate the computation, the Jacobi iteration is chosen as preconditioner for both non-compact and compact schemes. 
Various three-dimensional numerical experiments, from the subsonic to supersonic flows, are presented to to validate the accuracy and robustness of current implicit scheme. 
To accelerate the computation, the current schemes are implemented to run on graphics processing unit (GPU) using compute unified device architecture (CUDA). 
In the future, the implicit HGKS on arbitrary unstructured meshes with multiple GPUs will be developed for the engineering turbulent flows with high-Reynolds numbers. 

\section*{Acknowledgements}

The authors would like to thank Mr. Yue Zhang for helpful discussion.
The current research of L. Pan is supported by Beijing Natural
Science Foundation (1232012), National Natural Science Foundation of
China (11701038) and the Fundamental Research Funds for the Central
Universities, China.  The work of K. Xu is supported by National Key
R$\&$D Program of China (2022YFA1004500), National Natural Science
Foundation of China (12172316), and Hong Kong research grant council
(16208021,16301222).


\end{document}